\documentclass[12pt]{amsart}
\usepackage{amsbsy,amssymb,amscd,amsfonts,latexsym,amstext,delarray,
amsmath,graphicx} 
\input xypic

\usepackage{color}

\newtheorem{thm}{Theorem}[section]
\newtheorem{prop}[thm]{Proposition}
\newtheorem{cor}[thm]{Corollary}
\newtheorem{lem}[thm]{Lemma}

\newtheorem{defn}[thm]{Definition}
\newtheorem{rem}[thm]{Remark}

\numberwithin{equation}{section}

\def\bC{{\mathbb C}}

\def\bF{{\mathbb F}}
\def\bG{{\mathbb G}}

\def\bL{{\mathbb L}}

\def\bN{{\mathbb N}}

\def\bQ{{\mathbb Q}}

\def\bS{{\mathbb S}}
\def\bT{{\mathbb T}}

\def\bW{{\mathbb W}}
\def\bX{{\mathbb X}}

\def\bZ{{\mathbb Z}}

\def\A{{\mathbb A}}
\def\C{{\mathbb C}}
\def\F{{\mathbb F}}

\def\N{{\mathbb N}}
\renewcommand{\P}{{\mathbb P}}
\def\Q{{\mathbb Q}}
\def\R{{\mathbb R}}
\def\Z{{\mathbb Z}}
\def\K{{\mathbb K}}

\def\cA{{\mathcal A}}

\def\cC{{\mathcal C}}
\def\cD{{\mathcal D}}
\def\cE{{\mathcal E}}
\def\cF{{\mathcal F}}

\def\cI{{\mathcal I}}
\def\cJ{{\mathcal J}}

\def\cL{{\mathcal L}}
\def\cM{{\mathcal M}}
\def\cN{{\mathcal N}}
\def\cO{{\mathcal O}}
\def\cP{{\mathcal P}}
\def\cQ{{\mathcal Q}}
\def\cR{{\mathcal R}}
\def\cS{{\mathcal S}}
\def\cT{{\mathcal T}}

\def\cV{{\mathcal V}}
\def\cW{{\mathcal W}}

\def\cZ{{\mathcal Z}}

\def\Aut{{\rm Aut}}

\def\Hom{{\rm Hom}}

\def\SL{{\rm SL}}
\def\Spec{{\rm Spec}}
\def\Sp{{\rm Spec}}

\def\Tr{{\rm Tr}}

\def\fm{{\mathfrak m}}

\title[Bost--Connes and $\F_1$: Grothendieck rings, spectra, Nori motives]{Bost--Connes systems and $\F_1$-structures in Grothendieck
rings, spectra, and Nori motives} 
\author{Joshua F.~Lieber, Yuri I.~Manin, and Matilde Marcolli}
\date{}
\address{California Institute of Technology, Pasadena \\ USA}
\email{jlieber@caltech.edu}
\address{Max Planck Institute for Mathematics, Bonn \\ Germany}
\email{manin@mpim-bonn.mpg.de}
\address{California Institute of Technology, Pasadena \\ USA}
\email{matilde@caltech.edu}

\begin{document}
\maketitle

\begin{abstract}
We construct geometric lifts of the Bost--Connes algebra to Grothendieck rings and
to the associated assembler categories and spectra, as well as to certain categories of
Nori motives. These categorifications are related to the integral Bost--Connes algebra 
via suitable Euler characteristic type maps and zeta functions,
and in the motivic case via fiber functors. We also discuss aspects of $\F_1$-geometry,
in the framework of torifications, that fit into this general setting.
\end{abstract}

\section{Introduction and summary}\label{IntroSec} 

This survey/research paper interweaves many different strands that recently became 
visible in the fabric of algebraic geometry, arithmetics, (higher) category theory, quantum statistics, homotopical ``brave new algebra” etc., see especially A.~Connes and C.~Consani \cite{CoCo2}
\cite{CC16}; A.~Huber, St.~M\"uller–Stach \cite{HuM-S17}, etc.

\smallskip

In this sense, our present paper can be considered as a continuation and further extension of \cite{ManMar2},
and we will be relying on much of the work in that paper for details and examples.  The motivational starting
point in \cite{ManMar2} was coming from the interpretation given in \cite{CCM} of the Bost--Connes
quantum statistical mechanical system, and in particular the integral Bost--Connes algebra, as a form
of $\F_1$-structure, or ``geometry below $\Spec(\Z)$". The main theme of \cite{ManMar2} is an exploration
of how this structure manifests itself beyond the usual constructions of $\F_1$-structures on certain classes
of varieties over $\Z$.  In particular, the results of \cite{ManMar2} focus on lifts of the
integral Bost--Connes algebra to certain Grothendieck rings and to associated homotopy-theoretic spectra
obtained via assembler categories, and also on another form of $\F_1$-structures
arising through quasi-unipotent Morse-Smale dynamical systems. 

\smallskip

The main difference between the present paper and \cite{ManMar2} consists in a change 
of the categorical environment: the unifying vision we already considered in \cite{ManMar2} was 
provided by I.~Zakharevich’s notions of assemblers and scissors congruences: cf.~\cite{Zak1}, 
\cite{Zak2}, \cite{Zak3}, and \cite{CaWoZa}. In this paper, we continue to use the
formalism of assemblers and the associated spectra, but we complement it with 
categories of Nori motives, \cite{HuM-S17}. 

\smallskip

As in \cite{ManMar2}, we 
focus primarily on various geometrizations of the Bost--Connes  
algebra(s). Some of these constructions take place in Grothendieck
rings, like the previous cases considered in \cite{ManMar2}, and are aimed at lifting
the Bost--Connes endomorphisms to the level of 
homotopy theoretic spectra through the use of Zakharevich's formalism
of assembler categories. We focus on the case of relative
Grothendieck rings, endowed with appropriate equivariant Euler characteristic
maps. For varieties that admits torifications,
we introduce zeta functions based on the counting of points
over $\F_1$ and over extensions $\F_{1^m}$. 
We present a more general construction
of Bost--Connes type systems associated to exponentiable motivic
measures and the associated zeta functions with values in Witt rings, obtained using
a lift of the Bost--Connes algebra to Witt rings via Frobenius and
Verschiebung maps. 

\smallskip

We then consider lifts of the Bost--Connes algebra to
Nori motives, where we use a (slightly generalized)
version of Nori motives, which may be of independent interest in view of
possible versions of equivariant periods. In this categorical setting we show that the fiber
functor from Nori motives maps to a categorification of the Bost--Connes
algebra previously constructed by Tabuada and the third author, compatibly with
the functors realizing the Bost--Connes structure. 

\smallskip
\subsection{Structure of the paper and main results}\label{ResultsSec}   

Below we will briefly describe the content of the subsequent Sections, and the
main results of the paper, with pointers to the specific statements where 
these are proved.

\smallskip
\subsubsection{Bost--Connes systems and relative equivariant Grothendieck rings}\label{IntroGrSec}

In \S \ref{RelGrothSec}, we show the existence of a lift of the Bost--Connes
structure to the relative equivariant Grothendieck ring $K_0^{\hat\Z}(\cV_S)$,
extending similar results previously obtained in \cite{ManMar2} for the 
equivariant Grothendieck ring $K_0^{\hat\Z}(\cV)$. The main result in this part
of the paper is Theorem~\ref{liftBCmaps}, where the existence of this lift is
proved. The rest of the section covers the preliminary results needed for this
main result.

\smallskip

In particular, we first introduce the integral Bost--Connes
algebra in \S \ref{intBCSec}, in the form in which it was introduced in \cite{CCM}. 
We recall in \S \ref{relGrK0Sec} and \ref{eqGrSec} the relative and the
equivariant relative Grothendieck rings, and in \S \ref{EqEulCharSec} the
associated equivariant Euler characteristic map. 

\smallskip

In \S \ref{VerschSec} we recall from \cite{ManMar2} the geometric form of
the Verschiebung map that is used in the lifting of the Bost--Connes structure
to varieties with suitable $\hat\Z$-actions. In \S \ref{BCreleqGrK0Sec} we
introduce the Bost--Connes maps $\sigma_n$ and $\tilde\rho_n$ on
classes in $K_0^{\hat\Z}(\cV_S)$ and Proposition~\ref{BCliftSZ} shows 
the way they transform the varieties
and the base scheme with their respective $\hat\Z$-action. 

\smallskip

In \S \ref{PrimeSec} we recall from \cite{CCM} the prime decomposition
of the integral Bost--Connes algebra, which for a finite set of primes $F$
separates out an $F$-part and an $F$-coprime part of the algebra. We
then show in \S \ref{FcoprimeSec}, and in particular Proposition~\ref{liftBC},
that, given a scheme $S$ with a good effectively finite action of $\hat\Z$, 
there is an associated finite set of primes $F$ such that the $F$-coprime
part of the Bost--Connes algebra lifts to endomorphisms of
$K_0^{\hat\Z}(\cV_S)$. 

\smallskip

Finally in \S \ref{liftSec} we show how to
lift the full Bost--Connes algebra to homomorphisms between Grothendieck
rings $K_0^{\hat\Z}(\cV_{(S,\alpha)})$ where the scheme $S$ and the
action $\alpha$ are also transformed by the Bost--Connes map. By an
analysis of the structure of periodic points in Lemma~\ref{fixedPts} we
show the compatibility with the equivariant Euler characteristic, so we
can them prove the main result in Theorem~\ref{liftBCmaps}, showing
that the equivariant Euler characteristic intertwines the Bost--Connes
maps $\sigma_n$ and $\tilde\rho_n$ on the $K_0^{\hat\Z}(\cV_{(S,\alpha)})$ 
rings with the original $\sigma_n$ and $\tilde\rho_n$ maps of the integral
Bost--Connes algebra. 

\smallskip
\subsubsection{Bost--Connes systems on assembler categories and spectra}\label{IntroAssSec} 

In \S \ref{RingSpSec} we further lift the Bost--Connes structure obtained 
at the level of Grothendieck rings $K_0^{\hat\Z}(\cV_S)$ to assembler
categories underlying these Grothendieck rings and to the homotopy-theoretic
spectra defined by these categories. Again this extends to the equivariant 
relative case results that were obtained in \cite{ManMar2} for the non-relative setting.
The main result in this part of the paper is Theorem~\ref{liftAssSZ}, where it is
shown that the maps $\sigma_n$ and $\tilde\rho_n$ on the Grothendieck
rings $K_0^{\hat\Z}(\cV_{(S,\alpha)})$ constructed in the previous section
lift to functors of the underlying assembler categories, that induce these 
maps on $K_0$. 

\smallskip

In \S \ref{AssSubsec} we recall the formalism of assembler categories
of \cite{Zak1}, underlying scissor congruence relations and Grothendieck rings.
In \S \ref{GammaSpSec} we review Segal's $\Gamma$-spaces formalism
and how one obtains the homotopy-theoretic spectrum associated to an assembler.
In \S \ref{EnrichSec} and \S \ref{BCcatSec} we then lift this formalism by endowing
the main relevant objects with an action of a finite cyclic group,
with appropriate compatibility conditions. It is this
further structure that provides a framework for the respective
lifts of the Bost--Connes algebras, as in the cases discussed
in \cite{ManMar2} and in the ones we will be discussing in the
following sections.  We give here a very general definition of
Bost--Connes systems in categories, based on endofunctors
of subcategories of the automorphism category. In the applications
considered in this paper we will be using only the special case
where the automorphisms are determined by an effectively finite
action of $\hat\Z$, but we introduce the more general framework in
anticipation of other possible applications. 

\smallskip

In \S \ref{RelativeAssSec} we construct the assembler underlying the
equivariant relative Grothendieck ring $K_0^{\hat\Z}(\cV_S)$ and 
we prove the main result in Theorem~\ref{liftAssSZ} on the lift of
the Bost--Connes structure to functors of these assemblers.

\smallskip
\subsubsection{Bost--Connes systems on Grothendieck rings and assemblers of
torified varieties}

In \S \ref{TorifiedSec} we consider the approach to $\F_1$-geometry via torifications
of varieties over $\Z$, introduced in \cite{LoLo}. The main results of this part of the
paper are Proposition~\ref{assK0Ta} and Proposition~\ref{BCliftT} where we construct
assembler categories of torified varieties and we show the existence of a lift of the
Bost--Connes algebra to these categories. 

\smallskip

In \S \ref{TorifSec} we recall the notion of torified varieties
from \cite{LoLo} and the different versions of morphisms of torified varieties from
\cite{ManMar}, and we construct Grothendieck rings of torified varieties for each
flavor of morphisms In \S \ref{TorifRelSec} we introduce a relative version of
these Grothendieck rings of torified varieties. In \S \ref{QZactSec} we describe
$\Q/\Z$ and $\hat\Z$-actions on torifications. 

\smallskip

In \S \ref{AssTorifiedSec} we construct the assembler categories underlying
these relative Grothendieck rings and in \S \ref{BCTorifSec} we prove the first
main result of this section by constructing the lift of the Bost--Connes structure.

\smallskip
\subsubsection{Torified varieties, $\F_1$-points, and zeta functions}\label{IntroF1zetaSec}

This section continues the theme of torified varieties from the previous section but
with main focus on some associated zeta functions. We consider two different kinds
of zeta function: $\F_1$-zeta functions that count $\F_1$-points of torified
varieties, in an appropriates sense that it discussed in \S \ref{CountF1Sec}, and
dynamical zeta functions associated to endomorphisms of torified varieties
that are compatible with the torification. The use of dynamical zeta functions is
motivated by a proposal made in \cite{ManMar2} for a notion of $\F_1$-structures
based on dynamical systems that induce quasi-unipotent maps in homology. 

\smallskip

The two main results of this section are Proposition~\ref{zetaF1motmeas} and
Proposition~\ref{zetadynmotmeas} where we show that the $\F_1$-zeta function,
respectively the dynamical zeta function, determine exponentiable motivic
measures from the Grothendieck rings of torified varieties introduced in the
previous section to the ring $W(\Z)$ of Witt vectors. 

\smallskip

We introduce in \S \ref{CountF1Sec} and \S \ref{BBsec} the counting of $\F_1$-points 
of a torified variety and its relation to the Grothendieck class. We in show in \S \ref{BBsec} how
the  Bialynicki--Birula decomposition can be used to determine torifications and we give in \S \ref{BBexSec}
some explicit examples of computations of Grothendieck classes in simple cases that
have physical significance in the context of BPS counting in string theory. 

\smallskip

In \S \ref{CountF1zetaSec} we introduce the $\F_1$-zeta function and  
we prove Proposition~\ref{zetaF1motmeas}. In \S \ref{F1HWSec} we explain how
the $\F_1$-zeta function can be obtained from the Hasse--Weil zeta function.

\smallskip

In \S \ref{DynZetaSec} we consider torified varieties with dynamical systems compatible
with the torification and the associated Lefschetz and Artin--Mazur dynamical zeta functions.
We recall the definition and main properties of these zeta functions in 
\S \ref{DynZetaProSec} and we prove in Proposition~\ref{zetadynmotmeas} in \S \ref{TorDynZetaSec}.

\smallskip
\subsubsection{Spectrification of Witt vectors and lifts of zeta functions}\label{IntroZetaWRSec}

In the constructions described in \S\S~3 and 4 of \cite{ManMar2}
and in \S\S~\ref{RelGrothSec}--\ref{ZetaSec} of the present paper
we obtain lifts of the integral Bost--Connes algebra to various
assembler categories and associated spectra, starting from a ring
homomorphism (motivic measure) from the relevant Grothendieck
ring to the group ring $\Z[\Q/\Z]$ of the integral Bost--Connes
algebra, that is equivariant with respect to the maps $\sigma_n$
and $\tilde\rho_n$ of \eqref{eqA} and \eqref{eqB} of the Bost--Connes
algebra and the maps (also denoted by $\sigma_n$
and $\tilde\rho_n$) on the Grothendieck
ring induced by a Bost--Connes system on the corresponding
assembler category. The motivic measure provides in this way
a map that lifts the Bost--Connes structure. 

\smallskip

This part of the paper considers then a more general class of zeta functions $\zeta_\mu$ associated 
to exponentiable motivic measures $\mu: K_0(\cV) \to R$ with values in a commutative ring $R$, 
that admit a factorization into linear factors in the subring $W_0(R)$ of the Witt ring $W(R)$.

\smallskip

Our main results in this section are Proposition~\ref{PhimuGamma}, showing that
these zeta functions lift to the level of assemblers and spectra, 
and Proposition~\ref{zetamuBC}, which shows that the Frobenius and
Verschiebung maps on the endomorphism category lift, through the lift of the zeta function,
to a Bost--Connes system on the assembler category of the Grothendieck ring of
varieties $K_0(\cV)$. 

\smallskip

The main step toward establishing the main results of this section is the
construction in \S \ref{WittSpSec} and \S \ref{zetaGammaSec} 
of a spectrification of the ring $W_0(R)$. This is obtained using
its description in terms of the $K_0$ of the endomorphism category $\cE_R$
and of $R$, and the formalism of Segal Gamma-spaces. The spectrification
we use here is not the same as the spectrification of the ring of Witt vectors
introduced in \cite{Hess}. The lifting of Bost--Connes systems via motivic measures
in discussed in \S \ref{BCMotMeasSec}, where Proposition~\ref{zetamuBC} is proved.

\smallskip

We also consider again in \S \ref{SpSpSec} the setting on dynamical $\F_1$-structures
proposed in \cite{ManMar2}, with a pair $(X,f)$ of a variety and an endomorphism
that induces a quasi-unipotent map in homology, and we associate to these data
the operator-theoretic spectrum of the quasi-unipotent map, seen as an element
in $\Z[\Q/\Z]$. This determines a spectral map $\sigma: K_0^\Z(\cV_\C)\to \Z[\Q/\Z]$
with the properties of an Euler characteristic. 

\smallskip

Another main result in this section is Proposition~\ref{AssTann}, showing
that this spectral Euler characteristic lifts to a functor from the assembler
category underlying $K_0^\Z(\cV_\C)$ to the Tannakian category ${\rm Aut}^{\bar \Q}_{\Q/\Z}(\Q)$
that categorifies the ring $\Z[\Q/\Z]$, and that the resulting functor lifts the Bost--Connes
structure on ${\rm Aut}^{\bar \Q}_{\Q/\Z}(\Q)$ described in \cite{MaTa} to a Bost--Connes
structure on the assembler of $K_0^\Z(\cV_\C)$. 

\smallskip
\subsubsection{Bost--Connes systems in categories of Nori motives}\label{IntroNoriSec}  

When we replace the formalism of assembler categories and homotopy
theoretic spectra underlying the Grothendieck rings with geometric diagrams
and associated Tannakian categories of Nori motives, with the same notion
of categorical Bost--Connes systems introduced in Definitions~\ref{BCcat} 
and \ref{BCcat2}, we can lift the Euler characteristic type motivic measures
to the level of categorifications, where, as in the previous section, the
categorification of the Bost--Connes algebra is the one introduced in \cite{MaTa},
given by a Tannakian category of $\Q/\Z$-graded
vector spaces endowed with Frobenius and Verschiebung endofunctors. 

\smallskip

In \S \ref{NoriSec} in this paper we construct Bost--Connes systems in categories of
Nori motives. The main result of this part of the paper is Theorem~\ref{BCNori}, which
shows that there is a categorical Bost--Connes system on a category of equivariant
Nori motives, and that the fiber functor to the categorification of the Bost--Connes
algebra constructed in \cite{MaTa} intertwines the respective Bost--Connes endofunctors. 

\smallskip

In \S \ref{NoriDiagSubsec} and \S \ref{DiagNoriSec} we review the construction
of Nori motives from diagrams and their representations. In \S \ref{NoriTannSec} 
we construct a category of equivariant Nori motives. In \S \ref{BCNorisubsec}   
we describe the endofunctors of this category that implement the Bost--Connes
structure and we prove the main result in Theorem~\ref{BCNori}. In \S \ref{NoriAssSec}
we generalize this result to the relative case, using Arapura's motivic sheaves
version of Nori motives. 

\smallskip

Finally, in \S \ref{NoriAssSec} we consider Nori diagrams associated to 
assemblers and we formulate the question of their ``universal cohomological representations''.
This is a contemporary embodiment of the primordial Grothendieck's dream
that motives constitute a universal cohomology theory of algebraic varieties.

\section{Bost-Connes systems in Grothendieck rings}\label{RelGrothSec}  

In \cite{ManMar2} it was shown that the integral Bost--Connes algebra of \cite{CCM} admits 
lifts to certain Grothendieck rings, via corresponding equivariant Euler characteristic maps.
The cases analyzed in \cite{ManMar2} included the cases of the equivariant Grothendieck
ring $K_0^{\hat\Z}(\cV)$ and the equivariant Konsevich--Tschinkel Burnside 
ring ${\rm Burn}^{\hat\Z}(\K)$. We treat here, in a similar way, the case of the relative
equivariant Grothendieck ring $K_0^{\hat\Z}(\cV_S)$.  This case is more 
delicate than the other cases considered in \cite{ManMar2}, because when the Bost--Connes 
maps act on the classes in $K_0^{\hat\Z}(\cV_S)$ they also change the base scheme $S$ with its 
$\hat\Z$-action. 

\smallskip

The main result in this section is the existence of a lifting of the Bost--Connes 
structure to $K_0^{\hat\Z}(\cV_S)$, proved in Theorem~\ref{liftBCmaps}.

\smallskip

We first review the definition of the integral Bost--Connes algebra in \S \ref{intBCSec} 
and the equivariant relative Grothendieck ring in \S \ref{eqGrSec}. The rest of the section
then develops the intermediate steps leading to the proof of the main 
results of Theorem~\ref{liftBCmaps}.

\subsection{Bost-Connes algebra}\label{intBCSec}  

The Bost--Connes algebra was introduced in \cite{BC} as a quantum statistical mechanical
system that exhibit the Riemann zeta function as partition function, the generators of
the cyclotomic extensions of $\Q$ as values of zero-temperature KMS equilibrium states
on arithmetic elements in the algebra, and the abelianized Galois group $\hat\Z^*\simeq {\rm Gal}(\bar\Q/\Q)^{ab}$ as group of quantum symmetries. In particular, the arithmetic subalgebra of the 
 Bost--Connes system is given by the semigroup crossed product
\begin{equation}\label{QBC}
\Q[\Q/\Z] \rtimes \N
\end{equation} 
of the multiplicative semigroup $\N$ of positive integers acting on the group algebra of the group $\Q/\Z$.

\smallskip

The additive group $\bQ / \bZ$ can be identified with the multiplicative group 
$\nu^*$ of {\it roots of unity embedded into}  $\bC^*$: namely,
$r\in \bQ / \bZ$ corresponds to $e(r):=\exp (2\pi i\,r)$. More generally, the choice of the 
embedding can be modified by an arbitrary choice of an element in $\hat\Z^*=\Hom(\Q/\Z,\Q/\Z)$,
as is usually done in representations of the Bost--Connes algebra, see \cite{BC}. Thus,
we will use here interchangeably the notation $\zeta$ or $r$ for elements of $\Q/\Z$
assuming a choice of embedding as above. 
The group algebra $\bQ [\nu^*]$ consists of formal finite linear combinations 
$\sum_{a_{\zeta}\in\bQ} a_{\zeta} \zeta$
of roots of unity $\zeta\in \nu^*$. Formality means here that the sum is {\it not}
related to the additive structure of $\bC$.

\smallskip

The action of the semigroup $\N$ on $\Q[\Q/\Z]$ that defines the crossed product \eqref{QBC} is given by the endomorphisms
\begin{equation}\label{rhonQBC}
\rho_n (\sum a_{\zeta} \zeta) := \sum a_{\zeta}\,\, \frac{1}{n} \sum_{  \zeta^{\prime\, n} =\zeta }    \zeta^{\prime} .
\end{equation}
Equivalently, the algebra \eqref{QBC} is generated by elements $e(r)$ with the
relations $e(0)=1$, $e(r+r')=e(r)e(r')$, and elements $\mu_n$ and $\mu_n^*$ satisfying the relations
\begin{equation}\label{QBCrels1}
\begin{array}{ccc}
\mu_n^* \mu_n=1,\, \forall n;  & \mu_n \mu_n^* =\pi_n,  \, \forall n & \text{ with } \pi_n = \frac{1}{n} \sum_{  nr =0 } e(r); \\[3mm]
\mu_{nm}=\mu_n \mu_m, \, \forall n,m; &  \mu_{nm}^*=\mu_n^* \mu_m^* , \, \forall n,m; & \mu_n^* \mu_m=\mu_m\mu_n^* \text{ if } (n,m)=1.
\end{array}
\end{equation}
The semigroup action \eqref{rhonQBC} is then equivalently written as $\rho_n(a)=\mu_n\, a\, \mu_n^*$,
for all $a=\sum a_{\zeta} \zeta$ in $\Q[\Q/\Z]$.
The element $\pi_n\in \Q[\Q/\Z]$ is an idempotent, hence the generators $\mu_n$ are isometries
but not unitaries. 
See \cite{BC} and \S 3 of \cite{CoMa} for a detailed discussion of the
Bost--Connes system and the role of the arithmetic subalgebra \eqref{QBC}. 

\smallskip

In \cite{CCM} an integral model of the Bost--Connes algebra was constructed in order to
develop a model of $\F_1$-geometry in which the Bost--Connes system encodes the
extensions $\F_{1^m}$, in the sense of \cite{KapSmi}, of the ``field with one element" $\F_1$.

\smallskip

The integral Bost--Connes algebra is obtained by considering the group ring
$\Z[\Q/\Z]$, which we can again implicitly identify with $\Z[\nu^*]$ for a choice
of embedding $\Q/\Z\hookrightarrow \C$ as roots of unity. 

\smallskip

Define its {\it ring endomorphisms} $\sigma_n$:
\begin{equation}\label{eqA}
 \sigma_n (\sum a_{\zeta} \zeta):=
\sum a_{\zeta} \zeta^n.
\end{equation}

\smallskip

Define {\it additive maps} $\tilde{\rho}_n$: $\bZ [\nu^*] \to \bZ [\nu^*]$:
\begin{equation}\label{eqB}
\tilde{\rho}_n (\sum a_{\zeta} \zeta) := \sum a_{\zeta} \sum_{  \zeta^{\prime\, n} =\zeta }    \zeta^{\prime} .
\end{equation}

\smallskip

The maps $\sigma_n$ and $\tilde\rho_n$ satisfy the relations
\begin{equation}\label{sigmarhorels}
\sigma_n \circ \tilde\rho_n = n\, {\rm id}, \ \ \  \tilde\rho_n \circ \sigma_n = n \, \pi_n. 
\end{equation}

\smallskip

The integral Bost--Connes algebra is then defined as the algebra generated by
the group ring $\Z[\Q/\Z]$ and generators $\tilde\mu_n$ and $\mu_n^*$ with the
relations 
\begin{equation}\label{ZBCrel}
\begin{array}{ccc}
\tilde\mu_n \, a \, \mu^*_n =\tilde\rho_n(a), \, \forall n; &  \mu^*_n \, a = \sigma_n(a)\, \mu^*_n, \, \forall n; & a\, \tilde\mu_n=
\tilde\mu_n \, \sigma_n(a), \, \forall n; \\[3mm]
\tilde\mu_{nm}=\tilde\mu_n \tilde\mu_m, \, \forall n,m; & \mu^*_{nm}=\mu^*_n \mu^*_m, \, \forall n,m; &
\tilde\mu_n \mu^*_m =\mu^*_m    \tilde\mu_n \, \text{if } (n,m)=1.
\end{array}
\end{equation}
where the relations in the first line hold for all $a=\sum a_{\zeta} \zeta\in \Z[\Q/\Z]$, with $\sigma_n$
and $\tilde\rho_n$ as in \eqref{eqA} and \eqref{eqB}.

\smallskip

The maps $\tilde\rho_n$ of the integral Bost--Connes algebra and the
semigroup action $\rho_n$ in the rational Bost--Connes algebra \eqref{QBC}
are related by
$$
\rho_n = \frac{1}{n} \tilde{\rho}_n
$$
with $\tilde{\rho}_n$ defined as in \eqref{eqB}.

\smallskip
\subsection{Relative Grothendieck ring}\label{relGrK0Sec}   

We describe here a variant of construction of \cite{ManMar2}, where we work with
relative Grothendieck rings and with an Euler characteristic with values in a
Grothendieck ring of locally constant sheaves. We show that this relative
setting provides ways of lifting to the level of Grothendieck classes certain 
subalgebras of the integral Bost--Connes algebras associated to the choice of
a finite set of non-archimedean places. 

\smallskip

\begin{defn}\label{RelGrK0}  {\rm 
The relative Grothendieck ring $K_0(\cV_S)$ of varieties over a base variety 
$S$ over a field $\K$ is generated by the isomorphism classes of data 
$f: X \to S$ of a variety $X$ over $S$ with the relations
$$  [ f: X \to S ] = [f|_Y : Y \to S ] + [ f|_{X\smallsetminus Y}: X\smallsetminus Y \to S]  $$
as in \eqref{relK0S}
for a closed embedding $Y\hookrightarrow X$ of varieties over $S$. The product
is given by the fibered product $X\times_S Y$. We will write $[X]_S$ as shorthand notation 
for the class $[f: X\to S]$ in $K_0(\cV_S)$.  }
\end{defn}

\smallskip

A morphism $\phi: S \to S'$ induces a base change ring homomorphism
$\phi^* : K_0(\cV_{S'}) \to K_0(\cV_S)$ and a direct image map 
$\phi_*: K_0(\cV_S) \to K_0(\cV_{S'})$
which is a group homomorphisms and a morphism of $K_0(\cV_{S'})$-modules,
but not a ring homomorphism. The class $[ \phi: S \to S' ]$ as an element in 
$K_0(\cV_{S'})$ is the image of $1\in K_0(\cV_S)$ under $\phi_*$. 

\smallskip

When $S={\rm Spec}(\K)$ one recovers the ordinary Grothendieck ring $K_0(\cV_\K)$. 

\smallskip
\subsection{Equivariant relative Grothendieck ring}\label{eqGrSec}  

Let $X$ be a variety with a good action $\alpha: G\times X \to X$ by a finite 
group $G$ and $X'$ a variety with a good action $\alpha'$ by $G'$. As morphisms we
then consider pairs $(\phi,\varphi)$ of a morphism $\phi: X \to X'$
and a group homomorphism $\varphi: G \to G'$ such that
$\phi(\alpha(g,x))=\alpha'(\varphi(g),\phi(x))$, for all $g\in G$ and $x\in X$. 
Thus, isomorphisms of varieties with good $G$-actions are pairs of an
isomorphism $\phi: X \to X'$ of varieties and a group automorphism 
$\varphi\in \Aut(G)$ with the compatibility condition as above.

\smallskip

Given a base variety (or scheme) $S$ with a given good action $\alpha_S$ of a finite group $G$,
and varieties $X, X'$ over $S$, with good $G$-actions $\alpha_X, \alpha_{X'}$ and 
$G$-equivariant maps $f: X \to S$ and $f':X'\to S$, we consider morphisms given by a triple
$(\phi,\varphi, \phi_S)$ of a morphism $\phi: X \to X'$, a group homomorphism 
$\varphi: G \to G$ with the compatibility as above, and an endomorphism $\phi_S: S \to S$
such that $f'\circ \phi = \phi_S \circ f$. Then these maps also satisfy 
$\phi_S (\alpha_S(g,f(x))=\alpha_S (\varphi(g), \phi_S (f(x)))$. 

\smallskip

\begin{defn}\label{equivRelGrK0} {\rm 
The relative equivariant Grothendieck ring $K^G_0(\cV_S)$ is obtained as follows.
Consider the abelian group generated by isomorphism classes $[f: X \to S]$
of varieties over $S$ with compatible good $G$-actions, with respect to isomorphisms
$(\phi,\varphi, \phi_S)$ as above, with the inclusion-exclusion relations generated by
equivariant embeddings with
compatible $G$-equivariant maps
\begin{equation}\label{overSmaps}
 \xymatrix{ 
Y \ar@{^{(}->}[r] \ar[rd]_{f|_Y} & X \ar[d]^f & X\smallsetminus Y \ar@{_{(}->}[l] \ar[ld]^{f|_{X\smallsetminus Y}} \\
& S  & 
} \end{equation}
and isomorphisms. This means that we have
$[f: X \to S]=[f_Y: Y \to S] + [f_{X\smallsetminus Y}: X\smallsetminus Y \to S]$ if there are
isomorphisms $(\phi_Y,\varphi_Y,\phi_{S,Y})$ and $(\phi_{X\smallsetminus Y}, \varphi_{X\smallsetminus Y}, \phi_{S,X\smallsetminus Y})$, such that the diagram commutes
\begin{equation}\label{overSmaps2}
 \xymatrix{ Y \ar[r]^{\phi_Y} \ar[d]^{f_Y} & 
Y \ar@{^{(}->}[r] \ar[rd]_{f|_Y} & 
X \ar[d]^f 
& X\smallsetminus Y \ar@{_{(}->}[l] \ar[ld]^{f|_{X\smallsetminus Y}} & 
X\smallsetminus Y \ar[l]_{\phi'_{X\smallsetminus Y}} \ar[d]_{f_{X\smallsetminus Y}}
 \\ S \ar[rr]_{\phi_{S,Y}} 
& &  S  & & S \ar[ll]^{\phi_{S,X\smallsetminus Y}}
} \end{equation} 
The product $[f: X \to S]\cdot [f': X' \to S]$ given by $[\tilde f : X\times_S X' \to S]$
with $\tilde f= f\circ \pi_X =f'\circ \pi_{X'}$ is well defined on isomorphism classes,
with the diagonal action $\tilde\alpha(g,(x,x'))=(\alpha_X(g,x),\alpha_{X'}(g,x'))$
satisfying $f(\alpha_X(g,x))=\alpha_S(g,f(x))=\alpha_S(g,f'(x'))=f'(\alpha_{X'}(g,x'))$. 
}\end{defn}

\smallskip

We will use the following terminology for the $\hat\Z$-actions we consider.

\begin{defn}\label{efinact} {\rm 
A good effectively finite action of $\hat\Z$ on a variety $X$ is a good action
that factors through an action of some quotient $\Z/N\Z$. We will
write $\Z/N\Z$-effectively finite when we need to explicitly keep track 
of the level $N$. }
\end{defn}

\smallskip

In the case of the equivariant Grothendieck ring $K_0^{\hat\Z}(\cV)$ 
considered in \cite{ManMar2}, we can then also consider a relative version 
$K_0^{\hat\Z}(\cV_S)$, with $S$ a variety with a good effectively finite $\hat\Z$-action as above.
We consider the Grothendieck ring $K_0^{\hat\Z}(\cV_S)$ given by the isomorphism
classes of $S$-varieties $f: X\to S$ with good effectively finite $\hat\Z$-actions with respect to
which $f$ is equivariant, with the notion of isomorphism described above. 
The product is given by the fibered product over $S$ with the diagonal $\hat\Z$-action.
The inclusion-exclusion relations are as in \eqref{relK0S} where $Y\hookrightarrow X$ 
and $X\smallsetminus Y \hookrightarrow X$ are equivariant embeddings with
compatible $\hat\Z$-equivariant maps as in \eqref{overSmaps2}.

\smallskip
\subsection{Equivariant Euler characteristic}\label{EqEulCharSec}    

There is an Euler characteristic map given by a ring homomorphism
\begin{equation}\label{EulGS}
\chi^{\hat\Z}_S: K_0^{\hat\Z}(\cV_S) \to K_0^{\hat\Z}(\Q_S)
\end{equation}
to the Grothendieck ring of constructible sheaves over $S$ with $\hat\Z$-action,
\cite{GuZa}, \cite{Looij}, \cite{MaxSch}, \cite{Verdier}. 

\smallskip

\begin{lem}\label{fixedptsF}
Let $S$ be a variety with a good $\Z/N\Z$-effectively finite
$\hat\Z$-action.
Given a constructible sheaf $[\cF]$ in $K_0^{\hat\Z}(\Q_S)$, let $\cF|_{S^g}$ denote
the restrictions to the fixed point sets $S^g$, for $g\in\Z/N\Z$. These determine classes in
$K_0(\Q_{S^g})\otimes \Z[\Q/\Z]$. One obtains in this way a map
\begin{equation}\label{delochi}
 \chi: K_0^{\hat\Z}(\cV_S) \to \bigoplus_{g\in \Z/N\Z} K_0(\Q_{S^g})\otimes \Z[\Q/\Z].
\end{equation} 
\end{lem}
 
\proof The $\hat\Z$ action on $S$ factors through some $\Z/N\Z$,
hence the fixed point sets are given by $S^g$ for $g\in \Z/N\Z$.
Given a constructible sheaves $\cF$ over $S$ with $\hat\Z$-action, consider 
the restrictions $\cF|_{S^g}$. The subgroup $\langle g \rangle$ generated by $g$
acts trivially on $S^g$, hence for each $s\in S^g$ it acts on the stalk
$\cF_s$. Thus, these restrictions define classes $[\cF|_{S^g}]\in K_0(\Q_{S^g})\otimes 
R({\langle g \rangle}) \subset K_0(\Q_{S^g})\otimes\Z[\Q/\Z]$.
By precomposing with the Euler characteristic \eqref{EulGS} one then obtains
the map \eqref{delochi}. 
\endproof

\smallskip

We will also consider the map $K_0^{\hat\Z}(\cV_S) \to K_0(\Q_{S^G})\otimes \Z[\Q/\Z]$
given by the Euler characteristic followed by restriction of sheaves to the fixed point
set $S^G$ of the group action. 

\smallskip
\subsubsection{Fixed points and delocalized homology}\label{delocSec}  

Equivariant characteristic
classes from constructible sheaves to delocalized homology are
constructed in  \cite{MaxSch}.  

\smallskip

For a variety $S$ with a good action by a finite group $G$, and a (generalized)
homology theory $H$, the associated delocalized equivariant theory is given by 
$$ H^G(S) = (\oplus_{g\in G} H(S^g))^G $$
where the disjoint union $\sqcup_g S^g$ of the fixed point sets $S^g$ has 
an induced $G$-action $h: S^g\to S^{hgh^{-1}}$. In the case of an
abelian group we have $H^G(S) = (\oplus_{g\in G} H(S^g))^G$. 

\smallskip

As an observation, we can see explicitly the relation of delocalized homology 
to the integral Bost--Connes algebra, by considering the following cases (see
Remark~\ref{remHequiv}).
Let $S$ be a variety with a good $\Z/N\Z$-effectively finite $\hat\Z$-action. 
If $S$ has the trivial $\Z/N\Z$-action we have
$H^{\Z/N\Z}(S)=H(S)\otimes \Z[\Z/N\Z]$. In particular, if
$S$ is just a point, then this is $\Z[\Z/N\Z]$.
More generally, there is a morphism $$ \Z[\Z/N\Z] \times H^{\Z/N\Z}(S)\to H^{\Z/N\Z}(S) $$
induced by $H^{\Z/N\Z}(pt)\times H^{\Z/N\Z}(S)\to H^{\Z/N\Z\times \Z/N\Z}(pt\times S)\to H^{\Z/N\Z}(S)$
with the restriction to the diagonal subgroup as the last map.

\smallskip
\subsection{The geometric Verschiebung action}\label{VerschSec}  

We recall here how to construct the geometric Verschiebung action used
in \cite{ManMar2} to lift the Bost--Connes maps to the level of Grothendieck
rings. This has the effect of transforming an action of $\hat\Z$ on $X$
that factors through some $\Z/N\Z$ into an action of $\hat\Z$ on
$X\times Z_n$, with $Z_n = \{ 1, \ldots, n \}$, that factors through $\Z/Nn\Z$.
For $x\in X$, let $\underline x=(x,a_i)_{a_i \in Z_n}=(x_i)_{i=1}^n$
be the subset $\{ x \}\times Z_n$. For $\zeta_N$ a primmitive $N$-th root
of unity, we write in matrix form
$$ V_n(\zeta_{N n})=\left( \begin{array}{cccccc} 0 & 0 & \cdots & 0 & \alpha(\zeta_N) \\
1 & 0 & \cdots & 0 & 0 \\
0 & 1 & \cdots & 0 & 0  \\
\vdots & & \vdots & & \vdots \\
0 & 0 & \cdots & 1 & 0 \\
\end{array}\right)  $$
so that we can write
\begin{equation}\label{Verschiebung}
 V_n(\zeta_{Nn}) \cdot \underline{x} =\left\{ \begin{array}{ll}
(x, a_{i+1}) & i=1,\ldots, n-1 \\
(\alpha(\zeta_N)\cdot x, a_1) & i=n
\end{array}\right. 
\end{equation}
which satisfies $V_n(\zeta_{N n})^n=\alpha(\zeta_N) \times {\rm Id}_{Z_n}$.
The resulting action $\Phi_n(\alpha)$ of $\hat\Z$ on $X\times Z_n$ 
that factors through $\Z/Nn\Z$ is specified by setting
\begin{equation}\label{Phinact}
 \Phi_n(\alpha)(\zeta_{Nn})\cdot (x,a)= (V_n(\alpha(\zeta_N)) \cdot \underline{x})_a. 
 \end{equation}

\smallskip
\subsection{Lifting the Bost--Connes endomorphisms} \label{BCreleqGrK0Sec}   

Consider a base scheme $S$ with a good effectively finite action of $\hat\Z$. 
Let $f: X\to S$ be a variety over $S$ with a good effectively finite $\hat\Z$ action
such that the map is $\hat\Z$-equivariant. We denote by $\alpha_S: \hat\Z \times S \to S$
the action on $S$ and by $\alpha_X: \hat\Z \times X \to X$ the action on $X$. We write the equivariant
relative Grothendieck ring as $K_0^{\hat\Z}(\cV_{(S,\alpha_S)})$ to explicitly 
remember the fixed (up to isomorphisms as in \S \ref{eqGrSec}) action on $S$.

\smallskip

\begin{defn}\label{defBCliftSZ}{\rm
Let $(S,\alpha_S)$ be a scheme with a good effectively finite action of $\hat\Z$. 
Let $Z_n=\Spec(\Q^n)$ and let $\Phi_n(\alpha_S)$ be the action of $\hat\Z$ on $S\times Z_n$
as in \eqref{Verschiebung} and \eqref{Phinact}. Given a class $[f: (X,\alpha_X) \to (S,\alpha_S)]$ 
in $K_0^{\hat\Z}(\cV_{(S,\alpha_S)})$, with $\alpha_X$ the
compatible $\hat\Z$-action on $X$, let
\begin{equation}\label{sigmanK0S}
 \sigma_n [f: (X,\alpha_X) \to (S,\alpha_S)]= [f:(X,\alpha_X\circ \sigma_n)\to (S,\alpha_S\circ \sigma_n)] 
\end{equation}
\begin{equation}\label{rhonK0S}
 \tilde\rho_n[f: (X,\alpha_X) \to (S,\alpha_S)]=[f\times {\rm id}: (X\times Z_n,\Phi_n(\alpha_X))\to 
(S\times Z_n,\Phi_n(\alpha_S))]. 
\end{equation}
}\end{defn}

\smallskip

\begin{prop}\label{BCliftSZ}
For all $n\in\N$ the $\sigma_n$ defined in \eqref{sigmanK0S} are ring homomorphisms
\begin{equation}\label{sigmaSn}
\sigma_n : K_0^{\hat\Z}(\cV_{(S,\alpha_S)}) \to K_0^{\hat\Z}(\cV_{(S,\alpha_S\circ \sigma_n)})
\end{equation}
and the $\tilde\rho_n$ defined in \eqref{rhonK0S} are group homomorphisms 
\begin{equation}\label{rhoSn}
\tilde\rho_n: K_0^{\hat\Z}(\cV_{(S,\alpha_S)}) \to K_0^{\hat\Z}(\cV_{(S\times Z_n,\Phi_n(\alpha_S))}) ,
\end{equation}
with compositions satisfying
$$ \tilde\rho_n \circ \sigma_n: K_0^{\hat\Z}(\cV_{(S,\alpha_S)}) \to K_0^{\hat\Z}(\cV_{(S\times Z_n,\alpha_S\times \alpha_n)}) \to K_0^{\hat\Z}(\cV_{(S,\alpha_S)}) $$
$$ \sigma_n\circ \tilde\rho_n: K_0^{\hat\Z}(\cV_{(S,\alpha_S)}) \to K_0^{\hat\Z}(\cV_{(S,\alpha_S)^{\oplus n}})
\to K_0^{\hat\Z}(\cV_{(S,\alpha_S)}), $$
with $\sigma_n \circ \tilde\rho_n =n \, id$ and $\tilde\rho_n\circ \sigma_n$ is the product
by $(Z_n,\alpha_n)$.
\end{prop}

\proof Consider the $\sigma_n$ defined in \eqref{sigmanK0S}.
Since the group $\hat\Z$ is commutative and so is its endomorphism ring, these transformations $\sigma_n$
respect isomorphism classes since for an isomorphism $(\phi,\varphi,\phi_S)$ the actions satisfy
$$ \phi_X (\alpha_X(\sigma_n(g),x))=\alpha'_X(\varphi(\sigma_n(g)),\phi(x))
=\alpha'_X(\sigma_n(\varphi(g)),\phi(x)), $$ 
and similarly for the actions $\alpha_S, \alpha'_S$, so that $(\phi,\varphi,\phi_S)$ is also an
isomorphism of the images under $\sigma_n$. Similarly, the $\tilde\rho_n$ defined in \eqref{rhonK0S}
are well defined on the isomorphism classes.

As in \cite{ManMar2} we see that $\sigma_n \circ \tilde\rho_n[f: (X,\alpha_X) \to (S,\alpha_S)]=
[f: (X,\alpha_X) \to (S,\alpha_S)]^{\oplus n}$ and $\tilde\rho_n\circ \sigma_n
[f: (X,\alpha_X) \to (S,\alpha_S)]=[f\times {\rm id}: (X\times Z_n,\alpha_X\times \alpha_n)\to 
(S\times Z_n,\alpha_S\times \alpha_n)]$. The Grothendieck groups $K_0^{\hat\Z}(\cV_{(S\times Z_n,\alpha_S\times \alpha_n)})$ and $K_0^{\hat\Z}(\cV_{(S,\alpha_S)^{\oplus n}})$ map to 
$K_0^{\hat\Z}(\cV_{(S,\alpha_S)})$ via the morphism induced by composition with the natural maps
of the respective base varieties to $(S,\alpha_S)$.
\endproof

\smallskip

The fact that the ring homomorphisms \eqref{sigmaSn} and \eqref{rhoSn}  
determine a lift of the ring endomorphism $\sigma_n: \Z [\Q/\Z] \to \Z[\Q/\Z]$ and
group homomorphisms $\tilde\rho_n: \Z [\Q/\Z] \to \Z[\Q/\Z]$
of the integral Bost--Connes algebra is discussed in Proposition~\ref{liftBC}
and \S \ref{liftSec}.

\smallskip

We know from \cite{ManMar2} that the integral Bost--Connes algebra lifts to
the equivariant Grothendieck ring $K^{\hat\Z}(\cV_\Q)$ via maps $\sigma_n$
and $\tilde\rho_n$ that, respectively, precompose the action with the 
Bost--Connes endomorphism $\sigma_n$ and apply a geometric form
of the Verschiebung map. The main difference with the relative case considered here
lies in the fact that the lifts to the equivariant relative Grothendieck rings 
given by the maps \eqref{sigmaSn} and \eqref{rhoSn} need to transform in a compatible 
way the actions on both $X$ and $S$.

\begin{rem}\label{actXSrem} {\rm
Because the maps $\sigma_n$ and $\tilde\rho_n$ of \eqref{sigmaSn} and \eqref{rhoSn} simultaneously
modify the action on the varieties and on the base scheme $S$, they do not give endomorphisms
of the same $K_0^{\hat\Z}(\cV_{(S,\alpha_S)})$. However, given $(S,\alpha_S)$, it is possible
to identify a subalgebra of the integral Bost--Connes algebra that lift to 
endomorphisms of a corresponding subring of $K_0^{\hat\Z}(\cV_{(S,\alpha_S)})$, using the
notion of ``prime decomposition" of the Bost--Connes algebra. 
We discuss this more carefully in \S \ref{PrimeSec}, \S \ref{FcoprimeSec} and \S \ref{liftSec}.}
\end{rem}

\smallskip
\subsection{Prime decomposition of the Bost--Connes algebra}\label{PrimeSec}   

As in \cite{CCM}, for each prime $p$, we can decompose the group $\Q/\Z$ into a 
product $\Q_p/\Z_p \times (\Q/\Z)^{(p)}$, where $\Q_p/\Z_p$ is the Pr\"ufer group,
namely the subgroup of elements of $\Q/\Z$ where the denominator is a power of $p$,
isomorphic to $\Z[\frac{1}{p}]/\Z$, while $(\Q/\Z)^{(p)}$ consists of the elements with denominator prime to $p$.

\smallskip

Similarly, given a finite set $F$ of primes, we can decompose $\Q/\Z=(\Q/\Z)_F\times (\Q/\Z)^{F}$,
where the first term $(\Q/\Z)_F$ is identified with
fractions in $\Q/\Z$ whose denominator has prime factor decomposition consisting only
of primes in $F$, while elements in $(\Q/\Z)^{F}$ have denominators prime to all $p\in F$. 
The group ring decomposes accordingly as $\Z[(\Q/\Z)_F]\otimes \Z[(\Q/\Z)^{F}]$. 

\smallskip

The subsemigroup
$\N_F \subset \N$ generated multiplicatively by the primes $p\in F$ acts on
$\Z[(\Q/\Z)_F]\otimes_\Z \Q=\Q[(\Q/\Z)_F]$ by endomorphisms
$$ \rho_n (e(r))= \frac{1}{n} \sum_{nr' =r} e(r'), \ \ \  n\in \N_F, \,\, r\in (\Q/\Z)_F. $$
The corresponding morphisms $\sigma_n (e(r))=e(nr)$ and maps $\tilde\rho_n (e(r))=\sum_{nr' =r} e(r')$
act on $\Z[(\Q/\Z)_F]$ and we can consider the associated algebra $\cA_{\Z,F}$ generated by
$\Z[(\Q/\Z)_F]$ and $\tilde\mu_n$, $\mu_n^*$ with $n\in \N_F$, with the relations
\begin{equation}\label{intBCrels1}
\tilde\mu_{nm}=\tilde\mu_n \tilde\mu_m, \ \ \ \ \mu_{nm}^*=\mu_n^* \mu_m^*, \ \ \ \ 
\mu_n^* \tilde \mu_n =n, \ \ \ \ \tilde\mu_n \mu_m^* =\mu_m^* \tilde\mu_n, 
\end{equation}
 where the first two relations hold for arbitrary $n,m\in \N$, the third for arbitrary $n\in \N$ and
the forth for $n,m\in \N$ satisfying $ (n,m)=1$, and the relations
\begin{equation}\label{intBCrels2}
x \tilde\mu_n = \tilde\mu_n \sigma_n(x) \ \ \ \ \mu_n^* x = \sigma_n(x) \mu_n^*, \ \ \ \ 
\tilde\mu_n  x \mu_n^* = \tilde\rho_n(x),
\end{equation}
for any $x\in \Z[\Q/\Z]$, where $\tilde\rho_n(e(r))=\sum_{nr'=r} e(r')$, and 
with $$\cA_{\Z,F}\otimes_\Z\Q=\Q[(\Q/\Z)_F]\rtimes \N_F .$$
We refer to $\cA_{\Z,F}$ as the $F$-part of the  integral Bost--Connes algebra. 

\smallskip

The decomposition $\N = \N_F \times \N^{(F)}$, where $\N^{(F)}$ is generated by all
primes $p \notin F$, gives also an algebra 
$\cA_\Z^{(F)}$ generated by $\Z[(\Q/\Z)^F]$ and the $\tilde \mu_n$ and $\mu_n^*$ as
in \eqref{intBCrels2} with $p\notin F$ with
$$ \cA_\Z^{(F)} \otimes_\Z\Q=\Q[(\Q/\Z)^F]\rtimes \N^{(F)} .$$
We refer to $\cA_\Z^{(F)}$ as the
$F$-coprime part of the integral Bost--Connes algebra.

\smallskip
\subsection{Lifting the $F_N$-coprime Bost--Connes algebra}\label{FcoprimeSec}   

Let $F=F_N$ be the set of prime factors of $N$ and let $\Z[(\Q/\Z)^F]$ denote, as before,
the part of the group ring of $\Q/\Z$ involving only denominators relatively prime to $N$.
The semigroup $\N^{(F)}$ is generated by primes $p\!\!\not | N$ and we consider the
morphisms $\sigma_n (e(r))=e(nr)$ and maps $\tilde\rho_n (e(r))=\sum_{nr' =r} e(r')$ 
with $n\in \N^{(F)}$ and $r\in (\Q/\Z)^F$ as discussed above. 

\smallskip

\begin{prop}\label{liftBC}
Let $S$ be a base scheme with a good $\Z/N\Z$-effectively finite action of $\hat\Z$. 
Let $\cZ_{n,S}$ be defined as $\cZ_{n,S}=S \times Z_n$, with $Z_n =\Spec(\Q^n)$, with the 
action $\Phi_n(\alpha_S)$ obtained as in \eqref{Verschiebung} and \eqref{Phinact}.
The endomorphisms $\sigma_n: \Z [(\Q/\Z)^{F_N}] \to \Z[(\Q/\Z)^{F_N}]$ with $n\in \N^{(F_N)}$
of the $F_N$-coprime part of the integral Bost--Connes algebra 
lift to endomorphisms
$\sigma_n : K_0^{\hat\Z}(\cV_S) \to K_0^{\hat\Z}(\cV_S)$, as in \eqref{sigmanK0S}, which 
define a semigroup action of the multiplicative group $\N^{(F_N)}$ on the
Grothendieck ring $K_0^{\hat\Z}(\cV_S)$. 
The maps $\tilde\rho_n$, for $n\in \N^{(F_N)}$, 
lift to group homomorphisms $\tilde\rho_n : K_0^{\hat\Z}(\cV_S) \to K_0^{\hat\Z}(\cV_S)$, 
as in \eqref{rhonK0S}, so that $\sigma_n\circ \tilde\rho_n[f: X\to S]=[f: X\to S]^{\oplus n}$ and
$\tilde\rho_n\circ \sigma_n [f: X\to S]=[f: X\to S] \cdot \cZ_{n,S}$.
\end{prop}

\proof Given the base variety $S$ with a good $\Z/N\Z$-effectively finite $\hat\Z$-action, 
let $F=F_N$ denote
the set of prime factors of $N$. Let $X$ be a variety over $S$,
with a $\hat\Z$-equivariant map $f: (X,\alpha_X) \to (S,\alpha_S)$,
where we explicitly write the actions, satisfying $f(\alpha_X(\zeta, x))=\alpha_S (\zeta,f(x))$.
For $(N,n)=1$, the maps
$\sigma_n: [f: (X,\alpha_X) \to (S,\alpha_S)]=[f: (X,\alpha_X\circ \sigma_n) \to (S,\alpha_S\circ \sigma_n)]$,
as in \eqref{sigmanK0S}, satisfy
$(S,\alpha_S\circ \sigma_n)\simeq (S,\alpha_S)$
with the notion of isomorphism discussed in \S \ref{eqGrSec}, since $\zeta\mapsto \sigma_n(\zeta)$ is an automorphism of $\Z/N\Z$. Thus, the maps $\sigma_n$,
for $n\in \N^{(F_N)}$ determine a semigroup action of $\N^{(F_N)}$ by endomorphisms of 
$K_0^{\hat\Z}(\cV_S)$. 

Consider then $(\cZ_{n,N},\Phi_n(\alpha_S))$ as above, which we write equivalently
as $\tilde\rho_n(S,\alpha_S)$ where $\tilde\rho_n$ is the lift of the Bost--Connes
map to $K^{\hat\Z}(\cV)$ as in Proposition~3.5 of \cite{ManMar2}. We know that 
$\tilde\rho_n\circ\sigma_n [S,\alpha_S] = [S,\alpha_S] \cdot [Z_n,\alpha_n]$ in 
$K^{\hat\Z}(\cV)$. Since for
$(n,N)=1$ we have $(S,\alpha_S\circ \sigma_n) \simeq (S,\alpha_S)$, this gives
$(\cZ_{n,N},\Phi_n(\alpha_S))\simeq (S\times Z_n,\alpha_S \times \gamma_n)$. 
Then setting $\tilde\rho_n (f: X \to S)= (\tilde f: X\times_S \cZ_{n,S} \to S)$ with
$\tilde f= f\circ \pi_X$ gives $X\times_S \cZ_{n,S}\simeq X\times Z_n$, and the
composition properties for $\tilde\rho_n\circ \sigma_n$ and $\sigma_n\circ \tilde\rho_n$
are satisfied.

Given a class $[f: X \to S]$, let $[\cF_{X,S}]$ be the class in $K_0^{\hat\Z}(\Q_S)$
of the constructible sheaf given by the Euler characteristic \eqref{EulGS} of $[f: X \to S]$.
Let $[\cF_{X,S}|_{S^{\Z/N\Z}}]$ be the resulting class in $K_0(S^{\Z/N\Z})\otimes \Z[\Q/\Z]$
obtained by restriction to the fixed point set $S^{\Z/N\Z}$ with the element in $\Z[\Q/\Z]$
specifying the representation of $\hat\Z$ on the stalks of the sheaf $\cF_{X,S}|_{S^{\Z/N\Z}}$.
For $(N,n)=1$, the action of $\sigma_n$ by automorphisms of $\Z/N\Z$ wih the resulting
action by isomorphisms of $S$ induces an action by isomorphisms on the $K_0(S^{\Z/N\Z})$
part and the usual Bost--Connes action on $\Z[\Q/\Z]$. The restriction of the semigroup
action of $\N^{(F_N)}$ to the subring $\Z [(\Q/\Z)^{F_N}]$ is then the image of the
action of the maps $\sigma_n$ and $\tilde\rho_n$ on the preimage of this subring
under the morphism $K_0^{\hat\Z}(\cV_S) \to K_0(\Q_{S^G})\otimes \Z[\Q/\Z]$. 
\endproof

\smallskip

While this construction captures a lift of the $\Z [(\Q/\Z)^{F_N}]$ part of
the Bost--Connes algebra with the semigroup action of $\N^{(F_N)}$, 
the fact that the endomorphisms $\sigma_n$
acting on the roots of unity in $\Z/N\Z$ are automorphisms when $(N,n)=1$
loses some of the interesting structure of the Bost--Connes algebra, which
stems from the partial invertibility of these morphisms.  Thus, one also
wants to recover the structure of the complementary part of the
Bost--Connes algebra with the group ring $\Z[(\Q/\Z)_{F_N}]$ and the
semigroup $\N_{F_N}$. 

\smallskip
\subsection{Lifting the full Bost--Connes algebra}\label{liftSec}   

Unlike the $\Z [(\Q/\Z)^{F_N}]$ part of the Bost--Connes algebra
described above, when one considers the full Bost--Connes algebra,
including the $F_N$-part, the lift to the Grothendieck ring no longer
consists of endomorphisms of a fixed $K_0^{\hat\Z}(\cV_{(S,\alpha)})$,
but is given as in Proposition~\ref{BCliftSZ} by homomorphisms as in \eqref{sigmanK0S},
\eqref{sigmaSn} and \eqref{rhonK0S}, \eqref{rhoSn},
$$ \sigma_n : K_0^{\hat\Z}(\cV_{(S,\alpha_S)}) \to K_0^{\hat\Z}(\cV_{(S,\alpha_S\circ \sigma_n)}), $$
$$ \tilde\rho_n: K_0^{\hat\Z}(\cV_{(S,\alpha_S)}) \to K_0^{\hat\Z}(\cV_{(S\times Z_n,\Phi_n(\alpha_S))}). $$

\smallskip

For $G$ a finite abelian group with a good action $\alpha: G \times S \to S$ on a variety $S$, 
let $(S,\alpha)^G_k=\{ s\in S\,:\, \alpha(g^k,s)=s , \, \forall g\in G \}$ denote the
set of periodic points of period $k$, with $(S,\alpha)^G_1=(S,\alpha)^G$ the set of fixed points.
We always have $(S,\alpha)^G_k \subseteq (S,\alpha)_{km}^G$ for all $m\in \N$,
hence in particular a copy of the fixed point set $(S,\alpha)^G$ is contained in all $(S,\alpha)^G_k$.
For $G=\Z/N\Z$, with $\zeta_N$ a primitive $N$-th root of unity generator, the set of $k$-periodic
points is given by $(S,\alpha)^{\Z/N\Z}_k=\{ s\in S\,:\, \alpha(\zeta_N^k,s)=s \}$.

\begin{lem}\label{fixedPts}
The sets of periodic points satisfy $(S,\alpha\circ \sigma_n)^G_k=(S,\alpha)^G_{nk}$.
The sets $(S\times Z_n,\Phi_n(\alpha))^G_k$ can be non-empty only when $n|k$ with
$(S\times Z_n,\Phi_n(\alpha))^G_k = ((S,\alpha)_{k/n}^G)^n$.
\end{lem}

\proof Under the action $\alpha\circ \sigma_n$ the periodicity condition
means $\alpha\circ \sigma_n(\zeta^k,s)=\alpha(\zeta^{nk},s)=s$ for all 
$\zeta\in G$ hence the identification $(S,\alpha\circ \sigma_n)^G_k=(S,\alpha)^G_{nk}$.
In the case of the geometric Verschiebung action $\Phi_n(\alpha)$ on $S\times Z_n$,
the $k$-periodicity condition $\Phi_n(\alpha)(\zeta^k, (s,z))=(s,z)$ implies that $n|k$ for the
$k$-periodicity in the $z\in Z_n$ variable and that $\alpha(\zeta^{k/n},s)=s$.  
\endproof

The identification $(S,\alpha\circ \sigma_n)^G_k=(S,\alpha)^G_{nk}$ implies
the inclusion $(S,\alpha)^G_k \subseteq (S,\alpha\circ \sigma_n)^G_k$ and 
in particular the inclusion of the fixed point sets $(S,\alpha)^G \subseteq 
(S,\alpha\circ \sigma_n)^G$. Similarly, $(S\times Z_n,\Phi_n(\alpha))^G\subseteq ((S,\alpha)^G)^n$.
Since these inclusions will in general be strict, due to the fact that the endomorphisms
$\sigma_n$ are not automorphisms, one cannot simply use the map 
given by the equivariant Euler  characteristic followed by the restriction to the fixed point set 
$$ K_0^{\hat\Z}(\cV_S) \to K_0(\Q_{S^{\hat\Z}}) \otimes \Z[\Q/\Z] $$
to lift the Bost--Connes endomorphisms to the maps 
\eqref{sigmaSn} and \eqref{rhoSn} of Proposition~\ref{BCliftSZ}. However,
a simple variant of the same idea, where we consider sets of periodic points,
gives the lift of the full Bost--Connes algebra to the equivariant relative Grothendieck rings
$K_0^G(\cV_{(S,\alpha)})$.

\smallskip

Consider the equivariant Euler characteristic map followed by the restrictions to
the sets of periodic points
\begin{equation}\label{mapK0perPts}
K_0^G(\cV_{(S,\alpha)}) \stackrel{\chi_S^G}{\to} K_0^G(\Q_{(S,\alpha)}) \to \bigoplus_{k\geq 1} K_0^G(\Q_{(S,\alpha)^G_k}). 
\end{equation}
Also, for a given $n\in \N$, consider the same map composed with the projection to
the summands with $n|k$
\begin{equation}\label{mapK0perPtsn}
\chi_{S,n}^G: K_0^G(\cV_{(S,\alpha)}) \stackrel{\chi_S^G}{\to} K_0^G(\Q_{(S,\alpha)}) \to \bigoplus_{k\geq 1\,:\,n|k} K_0^G(\Q_{(S,\alpha)^G_k}). 
\end{equation}
For simplicity we consider the case where the fixed point set and
periodic points sets of the action $(S,\alpha)$ are all finite sets.

\smallskip

\begin{defn}\label{intertwinedef}  {\rm
Let $(S,\alpha)$ be a variety with a good effectively finite $\hat\Z$-action.
Consider data $(A_{(S,\alpha),n}, f_{(S,\alpha),n})$ and $(B_{(S,\alpha)},h_{(S,\alpha)})$ of a family of rings
$A_{(S,\alpha),n}$ with $n\in \N$ and $B_{(S,\alpha)}$ and ring homomorphisms
$f_{(S,\alpha),n}: K_0^G(\cV_{(S,\alpha)}) \to A_{(S,\alpha),n}\otimes \Z[\Q/\Z]$ and $h_{(S,\alpha)}: K_0^G(\cV_{(S,\alpha)}) \to B_{(S,\alpha)}\otimes \Z[\Q/\Z]$. The maps $f_{(S,\alpha),n}$ and $h_{(S,\alpha)}$ 
are said to {\rm intertwine the Bost--Connes structure} if 
there are ring isomorphisms $J_n: A_{(S,\alpha),n}\to B_{(S,\alpha\circ\sigma_n)}$ 
and isomorphisms of abelian groups $\tilde J_n: B_{(S,\alpha)}\to A_{(S\times Z_n, \Phi_n(\alpha))}$,
such that the following holds.
\smallskip
\begin{enumerate}
\item There is a commutative diagram of ring homomorphisms
$$ \xymatrix{ K_0^{\hat\Z}(\cV_{(S,\alpha)}) \ar[r]^{ f_{(S,\alpha),n} \quad \quad} \ar[d]^{\sigma_n} & A_{(S,\alpha),n} \otimes \Z[\Q/\Z]  \ar[d]^{J_n \otimes \sigma_n} \\ 
K_0^{\hat\Z}(\cV_{(S,\alpha \circ \sigma_n)}) \ar[r]^{\quad h_{(S,\alpha)} \quad \quad} & 
B_{(S,\alpha\circ \sigma_n)} \otimes \Z[\Q/\Z] 
} $$
where the maps $\sigma_n:K_0^{\hat\Z}(\cV_{(S,\alpha)})\to K_0^{\hat\Z}(\cV_{(S,\alpha \circ \sigma_n)})$
are as in \eqref{sigmaSn} and the maps 
$\sigma_n: \Z[\Q/\Z] \to \Z[\Q/\Z]$ are the endomorphisms of the integral Bost--Connes algebra.
\smallskip
\item There is a commutative diagram of group homomorphisms 
$$ \xymatrix{ K_0^{\hat\Z}(\cV_{(S,\alpha)}) \ar[rr]^{ h_{(S,\alpha)} \quad \quad} \ar[d]^{\tilde\rho_n} & & B_{(S,\alpha)} \otimes \Z[\Q/\Z]  \ar[d]^{\tilde J_n \otimes \tilde\rho_n} \\ 
K_0^{\hat\Z}(\cV_{(S\times Z_n,\Phi_n(\alpha))}) \ar[rr]^{\quad f_{(S\times Z_n, \Phi_n(\alpha)),n} \quad \quad} & & A_{(S\times Z_n,\Phi_n(\alpha)} \otimes \Z[\Q/\Z] 
} $$
where the maps $\tilde\rho_n: K_0^{\hat\Z}(\cV_{(S,\alpha)})\to K_0^{\hat\Z}(\cV_{(S\times Z_n,\Phi_n(\alpha))})$
are as in \eqref{rhoSn} and the $\tilde\rho_n: \Z[\Q/\Z]\to \Z[\Q/\Z]$ are the maps \eqref{eqB} of the integral
Bost--Connes algebra.
\end{enumerate} }
\end{defn}

\smallskip

\begin{thm}\label{liftBCmaps}
Let $(S,\alpha)$ be a variety with a good effectively finite $\hat\Z$-action,
such that the set $(S,\alpha)^{\hat\Z}_k$ of $k$-periodic points for this action is finite, for all $k\geq 1$.
Then the maps \eqref{mapK0perPts} and \eqref{mapK0perPtsn} intertwine the Bost--Connes structure
in the sense of Definition~\ref{intertwinedef}.
\end{thm}

\proof
Under the assumptions that all the $(S,\alpha)^G_k$ for $k\geq 0$ are finite sets, we can
identify the target of the map with
$\oplus_k K_0(\Q_{(S,\alpha)^G_k}) \otimes R(G)$.
In the case of varieties with good effectively finite $\hat\Z$ actions,
we obtain in this way ring homomorphisms
$$ \chi^{\hat\Z}_{(S,\alpha)}: K_0^{\hat\Z}(\cV_{(S,\alpha)})\to 
\bigoplus_{k\geq 1} K_0(\Q_{(S,\alpha)^{\hat\Z}_k}) \otimes \Z[\Q/\Z] $$
$$ \chi^{\hat\Z}_{(S,\alpha),n} : K_0^{\hat\Z}(\cV_{(S,\alpha)})\to 
\bigoplus_{k\geq 1\,:\,n|k} K_0(\Q_{(S,\alpha)^{\hat\Z}_k}) \otimes \Z[\Q/\Z]. $$
These maps fit in the following commutative diagrams of ring homomorphisms
$$ \xymatrix{ K_0^{\hat\Z}(\cV_{(S,\alpha)}) \ar[r]^{ \chi^{\hat\Z}_{(S,\alpha),n} \quad \quad} \ar[d]^{\sigma_n} & \bigoplus_{n|k} K_0(\Q_{(S,\alpha)^{\hat\Z}_k}) \otimes \Z[\Q/\Z]  \ar[d]^{J_n \otimes \sigma_n} \\ 
K_0^{\hat\Z}(\cV_{(S,\alpha \circ \sigma_n)}) \ar[r]_{\bar\chi^{\hat\Z}_{(S,\alpha\circ \sigma_n)} \quad \quad} & 
\bigoplus_\ell K_0(\Q_{(S,\alpha \circ \sigma_n)^{\hat\Z}_\ell}) \otimes \Z[\Q/\Z] 
} $$
where the map $(J_n)_{k,\ell}$ is non-trivial for $k=\ell n$ and identifies $K_0(\Q_{(S,\alpha)^{\hat\Z}_\ell})$ with $K_0(\Q_{(S,\alpha \circ \sigma_n)^{\hat\Z}_k})$, while the maps $\sigma_n: \Z[\Q/\Z] \to \Z[\Q/\Z]$ are the Bost--Connes endomorphisms. Similarly, we obtain commutative diagrams of group homomorphisms
$$ \xymatrix{ K_0^{\hat\Z}(\cV_{(S,\alpha)}) \ar[r]^{ \chi^{\hat\Z}_{(S,\alpha)} \quad \quad} \ar[d]^{\tilde\rho_n} & \bigoplus_\ell K_0(\Q_{(S,\alpha)^{\hat\Z}_\ell}) \otimes \Z[\Q/\Z]  \ar[d]^{\tilde J_n \otimes \tilde\rho_n} \\ 
K_0^{\hat\Z}(\cV_{(S\times Z_n,\Phi_n(\alpha))}) \ar[r]_{\chi^{\hat\Z}_{(S\times Z_n, \Phi_n(\alpha)),n} \quad \quad} & \bigoplus_{n|k} K_0(\Q_{(S\times Z_n,\Phi_n(\alpha))^{\hat\Z}_k}) \otimes \Z[\Q/\Z] 
} $$
where $(\tilde J_n)_{\ell,k}$ is non-trivial for $k=\ell n$ and maps $K_0(\Q_{(S,\alpha)^{\hat\Z}_k})$  to $K_0(\Q_{(S,\alpha)^{\hat\Z}_k})^{\oplus n}$
and identifies the latter with $K_0(\Q_{(S\times Z_n,\Phi_n(\alpha))^{\hat\Z}_\ell})$.
\endproof

\smallskip

\begin{rem}\label{remHequiv}{\rm
A similar argument can be given using a map obtained by composing the equivariant Euler
characteristic considered here with values in $K_0^{\hat\Z}(\Q_S)$ with equivariant characteristic
classes from constructible sheaves to delocalized equivariant homology as in \cite{MaxSch}, see
\S \ref{delocSec}. }
\end{rem}

\medskip

\section{From Grothendieck Rings to Spectra}\label{RingSpSec}   

In this section we show that the Bost--Connes structure can be lifted further from
the level of the relative Grothendieck ring to the level of spectra, using the assembler
category construction of \cite{Zak1}.

\smallskip

The results of this section are a natural continuation of the results in \cite{ManMar2}.
The general theme considered there consisted of the following steps:
\begin{itemize}
\item Appropriate
equivariant Euler characteristic maps from certain $\hat\Z$-equivariant Grothendieck
rings to the group ring $\Z[\Q/\Z]$ are constructed.
\item These Euler characteristic maps are then used to lift the
Bost--Connes operations $\sigma_n$ and $\tilde \rho_n$ from $\Z[\Q/\Z]$ to
corresponding operations in the equivariant Grothendieck ring.
\item Assembler categories with $K_0$ given by the equivariant Grothendieck ring are constructed.
\item Endofunctors $\sigma_n$ and $\tilde \rho_n$ of these assembler categories 
are constructed so that they
induce the Bost--Connes structure in the Grothendieck ring.
\item Induced maps of spectra are obtained from these endofunctors through the
Gamma-space construction that associated a spectrum to an assembler category.
\end{itemize}
The construction of Bost--Connes operations $\sigma_n$ and $\tilde \rho_n$ on the
equivariant Grothendieck rings was generalized in the previous section to the
case of relative Grothendieck rings. This section deals with the corresponding
generalization of the remaining steps. 

\smallskip 

We start this section by a brief survey in \S \ref{AssSubsec} of Zakharevich's formalism of
{\it assemblers} which axiomatizes the ``scissors congruence'' relations
\eqref{relK0S}.

\smallskip

A general framework for categorical Bost--Connes systems is introduced in \S \ref{EnrichSec}
and \S \ref{BCcatSec} in terms of subcategories of the automorphism category (in our
examples encoding the $\hat\Z$-actions) and endofunctors $\sigma_n$ and $\tilde\rho_n$
implementing the Bost--Connes structure. 

\smallskip

In \S \ref{RelativeAssSec} we construct an assembler category for the
equivariant relative Grothendieck ring, and we prove the main result
of this section, Theorem~\ref{liftAssSZ}, on the lifting of the Bost--Connes
structure to this assembler category. 

\smallskip
\subsection{Assemblers}\label{AssSubsec}    

Below we will recall the basics of        
a general formalism for scissors congruence relations applicable in algebraic geometric
contexts defined by
I.~Zakharevich  in \cite{Zak1} and \cite{Zak2}.
The abstract form of scissors congruences consists of
categorical data called {\it assemblers}, which in turn determine a homotopy--theoretic
{\it spectrum}, whose homotopy groups embody scissors congruence relations.
This formalism is applied in \cite{Zak3} in the framework producing an assembler
and a spectrum whose $\pi_0$ recovers the Grothendieck ring of varieties.
This is used to obtain a characterisation  of the kernel of multiplication by
the Lefschetz motive, which provides a general explanation for the 
observations of \cite{Bor14}, \cite{Mart16} on the fact that the Lefschetz motive
is a zero divisor in the Grothendieck ring of varieties.

\smallskip

Consider  a (small) category $\cC$ and an object $X$  in $\cC$.
\smallskip

\begin{defn}\label{sievedef} {\rm
 A {\it sieve} $\cS$ over $X$ in $\cC$ is a family of
morphisms $f_i: X^{\prime}_i\to X$  (also called ``objects over $X$'') satisfying
the following conditions:

\begin{itemize}
\item[a)]
     Any isomorphism with target $X$ belongs to $\cS$ (as a family with one element).
    
    \item[b)] If a morphism $X^{\prime} \to X$ belongs to $\cS$, then its precomposition
with any other morphism in $\cC$ with target $X^{\prime}$    
$$
X^{\prime\prime}\to X^{\prime}\to  X
$$
also belongs to $\cS$.
\end{itemize} }
\end{defn}

\medskip

It follows that composition of any two morphisms in $\cS$ composable in $\cC$
itself belongs to $\cS$ so that any sieve is a category in its own right.

\smallskip

\begin{defn}\label{grotop} {\rm 
A Grothendieck
topology on a category $\cC$ consists of the assignment of a 
collection of sieves $\cJ(X)$  given for all 
objects $X$ in $\cC$, with the following  properties: 

\begin{itemize}
\item[a)] the total overcategory $\cC /X$ of morphisms with target $X$ is a member of
the collection $\cJ(X)$.

\item[b)] The pullback of any sieve
in $\cJ(X)$ under a morphism $f: Y\to X$ exists and is a sieve in $\cJ(Y)$.
Here pullback of a sieve is defined as the family of pullbacks of its objects, $X^{\prime}\to X$,
whereas pullback of  such an object w.r.t.~$Y\to X$ is defined as $pr_{Y}:\,Y\times_{X} X^{\prime} \to Y$.

\item[c)] given $\cC'\in \cJ(X)$ and a sieve   $\cT$ in $\cC/X$, if for every
$f: Y \to X$ in $\cC'$ the pullback $f^*\cT$ is in $\cJ(Y)$ then
$\cT$ is in $\cJ(X)$. 
\end{itemize} }
\end{defn}

For more details, see \cite{KSch06}, Chapters~16 and 17, or \cite{HuM-S17}, pp.~20--22.

\smallskip

Let $\cC$ be a category with a Grothendieck topology. Zakharevich's notion
of an assembler category is then defined as follows. 

\smallskip

\begin{defn}\label{covfam} {\rm
A collection of morphisms $\{ f_i: X_i \to X \}_{i\in I}$ in $\cC$  
is a {\it covering family} if the full subcategory of $\cC/X$ 
that contains all the morphisms of $\cC$ that
factor through the $f_i$,
$$
\{ g: Y \to X\,|\, \exists i\in I\, \, h: Y\to X_i\,\, \text{ such that } f_i \circ h = g \}, 
$$
belongs to the sieve collection $\cJ(X)$. }
\end{defn}

\smallskip

In a category $\cC$ with an initial object $\emptyset$ two morphisms $f: Y \to X$
and $g: W\to X$ are called {\it disjoint} if the pullback $Y\times_X W$ exists and is equal to $\emptyset$.
A collection $\{ f_i: X_i \to X \}_{i\in I}$ in $\cC$ is disjoint if $f_i$ and $f_j$ are disjoint
for all $i\neq j \in I$.

\smallskip

\begin{defn}\label{assdef} {\rm 
An assembler category $\cC$ is a small category endowed with
a Groth\-endieck topology, which has an initial object
$\emptyset$ (with the empty family as covering family), and where all 
morphisms are monomorphisms, with the property that any two finite 
disjoint covering families of $X$ in $\cC$ have a common refinement 
that is also a finite disjoint covering family. }
\end{defn}

\smallskip

A morphism of assemblers is a functor $F: \cC \to \cC'$ that is continuous
for the Grothendieck topologies and preserves the initial object and the
disjointness property, that is, if two morphisms are disjoint in $\cC$ their images
are disjoint in $\cC'$.

\smallskip

For $X$ a finite set, the coproduct of assemblers $\bigvee_{x\in X} \cC_x$
is a category whose 
objects are the initial object $\emptyset$ and all the non--initial objects of 
the assemblers $\cC_x$. Morphisms of non--initial objects are induced by those 
of $\cC_x$. 

\smallskip

Consider a pair $(\cC, \cD)$ where $\cC$
is an assembler category, and $\cD$ is a sieve in $\cC$.

\smallskip

One has then an associated
assembler $\cC \smallsetminus \cD$ defined as the full subcategory of $\cC$
containing all the objects that are not initial objects of $\cD$. The assembler
structure on $\cC \smallsetminus \cD$ is determined by taking as covering families
in $\cC \smallsetminus \cD$ those collections $\{ f_i : X_i \to X \}_{i\in I}$ with $X_i, X$
objects in $\cC \smallsetminus \cD$  that can be completed to a covering family in $\cC$,
namely such that there exists $\{ f_j: X_j \to X \}_{j\in J}$ with $X_j$ in $\cD$ such that
$\{ f_i : X_i \to X \}_{i\in I} \cup \{ f_j: X_j \to X \}_{j\in J}$ is a covering family in $\cC$.

\smallskip

Moreover,
there is a morphism of assemblers $\cC \to \cC \smallsetminus \cD$ that maps objects
of $\cD$ to $\emptyset$ and objects of $\cC \smallsetminus \cD$ to themselves and
morphisms with source in $\cC \smallsetminus \cD$ to themselves and morphisms
with source in $\cD$ to the unique morphism to the same target with source $\emptyset$. 
The data $(\cC, \cD, \cC \smallsetminus \cD)$ are called the abstract scissors congruences.

\smallskip

The construction of $\Gamma$-spaces, which we review more in 
detail in \S \ref{GammaSpSec}, then provides the homotopy theoretic 
spectra associated to assembler categories as in \cite{Zak1}. 
This construction of assembler categories and spectra provides
the formalism  we use here and in the previous paper \cite{ManMar2} 
to lift Bost--Connes type algebras to the level of Grothendieck rings and spectra.

\smallskip
\subsection{From categories to $\Gamma$-spaces and spectra} \label{GammaSpSec} 

The Segal construction \cite{Segal} associates a $\Gamma$-space (hence a spectrum)
to a category $\cC$ with a zero object and a categorical sum. Let $\Gamma^0$ be
the category of finite pointed sets with objects $\underline{n}=\{0,1,\ldots, n\}$ and
morphisms $f: \underline{n}\to \underline{m}$ the functions with $f(0)=0$. Let $\Delta_*$
denote the category of pointed simplicial sets. 
The construction can be generalized to symmetric monoidal categories, \cite{Thoma}.
The associated $\Gamma$-space
$F_\cC: \Gamma^0 \to \Delta_*$ is constructed as follows. First assign to a finite pointed set $X$
the category $P(X)$ with objects all the pointed subsets of $X$ with morphisms given
by inclusions. A functor $\Phi_X: P(X) \to \cC$ is summing if it maps $\emptyset \in P(X)$
to the zero object of $\cC$ and given $S,S'\in P(X)$ with $S\cap S'=\{ \star \}$ the base point of $X$,
the morphism $\Phi_X(S)\oplus \Phi_X(S') \to \Phi_X(S\cup S')$ is an isomorphism. Let
$\Sigma_\cC(X)$ be the category whose objects are the summing functors $\Phi_X$ with
morphisms the natural transformations that are isomorphisms on all objects of $P(X)$. Setting
$$ \Sigma_\cC(f)(\Phi_X)(S)=\Phi_X(f^{-1}(S)), $$
for a morphisms $f: X\to Y$ of pointed sets and $S\in P(Y)$ gives a functor
$\Sigma_\cC: \Gamma^0 \to {\rm Cat}$ to the category of small categories. 
Composing with the nerve $\cN$ gives a functor
$$ F_\cC = \cN\circ \Sigma_\cC : \Gamma^0 \to \Delta_* $$
which is the $\Gamma$-space associated to the category $\cC$. 
The functor $F_\cC: \Gamma^0 \to \Delta_* $ obtained in this way
is extended to an endofunctor $F_\cC: \Delta_* \to \Delta_*$ via the coend 
$$ F_{\cC}(K)=\int^{\underline{n}} K^n \wedge F_{\cC}(\underline{n}). $$
One obtains the spectrum $\bX=F_\cC(\bS)$ associated to the $\Gamma$-space $F_\cC$
by setting $\bX_n =F_\cC(S^n)$ with maps $S^1\wedge F_\cC(S^n) \to F_\cC(S^{n+1})$.
The construction is functorial in $\cC$, 
with respect to functors preserving sums and the zero object.

\smallskip

When $\cC$ is the category of finite sets, $F_\cC(\bS)$ is the sphere
spectrum $\bS$, and when $\cC=\cP_R$ is the category of finite projective
modules over a commutative ring $R$, the spectrum $F_{\cP_R}(\bS)=K(R)$
is the $K$-theory spectrum of $R$.

\smallskip

The Segal construction determines a functor from the category of small symmetric 
monoidal categories to the category of $-1$-connected spectra. It is shown in
\cite{Thoma} that this functor determines an equivalence of categories between the 
stable homotopy category of $-1$-connected spectra and a localization of the category 
of small symmetric monoidal categories, obtained by inverting morphisms sent to weak 
homotopy equivalences by the functor. 

\smallskip

Given an assembler category $\cC$, one considers a category $\cW(\cC)$ with
objects $\{ A_i \}_{i\in I}$ given by collections of non-initial objects $A_i$ in $\cC$
indexed by finite sets and morphisms $\phi: \{ A_i \}_{i\in I} \to \{ B_j \}_{j\in J}$ consisting of
a map of finite sets $f: I \to J$ and morphisms $\phi_i: A_i \to B_{f(i)}$ that form
disjoint covering families $\{ \phi_i\,|\, i\in f^{-1}(j)\}$, for all $j\in J$. One then obtains
a $\Gamma$-space as the functor that assigns to a finite pointed set $(X,x_0)$ the
simplicial set $\cN\cW(X\wedge \cC)$, the nerve of the category $\cW(X\wedge \cC)$
where $X\wedge \cC$ is the assembler $X\wedge \cC=\bigvee_{x\in X \smallsetminus \{ x_0\}} \cC$.
The spectrum associated to the assembler $\cC$ is the spectrum defined by this $\Gamma$ space,
namely $X_n=\cN\cW(S^n\wedge \cC)$.

\smallskip

For another occurrence of $\Gamma$-spaces in the context of $\F_1$-geometry,
see \cite{CC16}.

\medskip
\subsection{Automorphism category and enhanced assemblers}\label{EnrichSec}  

We describe
in this and the next subsection a general formalism of ``enhanced assemblers" underlying all
the explicit cases of Bost--Connes structures in Grothendieck rings discussed in \cite{ManMar2} and 
in some of the later sections of this paper.  

\smallskip

We first recall the definition of the automorphism category.

\begin{defn}\label{defAutoCat}{\rm  The automorphism category ${\rm Aut}(C)$
of a category $C$ is given by:

\begin{itemize}

\item[(i)] Objects of ${\rm Aut}(C)$ are pairs $\hat X=(X,v_X)$ where $X\in {\rm Obj}(C)$
and $v_X: X\to X$ is an automorphism of $X$.

\item[(ii)] Morphisms $\hat{f}:\, (X,v_X)\to (Y,v_Y)$ in ${\rm Aut}(C)$ are
morphisms $f:\, X\to Y$ such that $f\circ v_X=v_Y\circ f:\, X\to Y$ in
$C$.

\item[(iii)] The forgetful functor sends $\hat{X}$ to $X$ and $\hat{f}$ to $f$.
\end{itemize}}\end{defn} 

\smallskip

We use here a standard categorical notation according to which,
say, $f\circ v_X$ is the precomposition of $f$ with $v_X$.

\smallskip

Thus, we can make the following general definition. In the following we will
be especially interested in the case where the chosen subcategory is
determined by a group action, see Remark~\ref{AutCrem}.

\begin{defn}\label{enrichdef} {\rm 
Let  $C$ be a category. We will call here an {\it enhancement} of $C$ 
a pair consisting of a choice of a subcategory $\hat{C}$ of the automorphism 
category ${\rm Aut}(C)$ and the forgetful functor $\hat{C}\to C$, where objects 
$(X,v_X)$ of $\hat{C}$ have automorphisms
$v_X: X\to X$ of finite order, . }
\end{defn}

\smallskip

The main idea here is that a subcategory category $\hat C$ of the 
automorphism category of $C$ is where 
the endofunctors defining the lifts of the Bost--Connes structure are defined,
as we make more precise in Definitions~\ref{BCcat} and \ref{BCcat2}.

\smallskip

\begin{rem}\label{AutCrem}{\rm 
In the cases considered in \cite{ManMar2}
and in this paper, the subcategory of $\hat{C}$ of ${\rm Aut}(C)$ is usually determined by
a finite group action, so that elements of
$\hat{C}$ are of the form $(X,\alpha_{X}(g))$ with $\alpha_X: G \times X \to X$ the group
action. However, one expects other interesting examples that are
not necessarily given by group actions, hence it is worth considering this more general formulation.}
\end{rem}

\smallskip

\begin{rem}\label{remLiftAss}{\rm 
Assume that  $C$ is endowed with a structure of assembler. 
Then a series of constructions presented in \S\S~3 and 4 of \cite{ManMar2} 
and in \S\S~\ref{RelGrothSec}--\ref{ZetaSec} of this paper,
and restricted there to  various categories of schemes,
show in fact how  this structure of assembler can be lifted from
$C$ to $\hat{C}$.}
\end{rem}

\smallskip

In particular the Bost--Connes type structures we are investigating can be
formulated broadly in this setting of enhanced assemblers as follows.

\smallskip

\smallskip
\subsection{Bost--Connes systems on categories}\label{BCcatSec}   

Let $\hat C$ be an enhancement of a category $C$, in the sense of
Definition~\ref{enrichdef}. 

\begin{defn}\label{BCcat}  {\rm 
We assume here that $C$ is an 
additive (symmetric) monoidal category and that the enhancement
$\hat C$ is compatible with this structure. 
A Bost--Connes system in an enhancement $\hat C$ of $C$ consists of 
two families of endofunctors $\{ \sigma_n \}_{n\in \N}$ and
$\{ \tilde\rho_n \}_{n\in \N}$ of $\hat C$  with the following
properties:
\begin{enumerate}
\item The functors $\sigma_n$ are compatible with both the
additive and the (symmetric) monoidal structure, while the
functors $\tilde\rho_n$ are functors of additive categories.
\item For all $n,m\in \N$ these endofunctors satisfy
$$ \sigma_{nm}=\sigma_n \circ \sigma_m, \ \ \ \  \tilde\rho_{nm}=\tilde\rho_n \circ \tilde\rho_m. $$
\item The compositions satisfy
\begin{equation}\label{endofunctorBCrels}
 \sigma_n\circ \tilde\rho_n (X,v_X)=(X,v_X)^{\oplus n} \ \ \ \text{ and } \ \ \   \tilde\rho_n \circ \sigma_n (X,v_X)
= (X,v_X) \otimes (Z_n,v_n), 
\end{equation}
for an object $(Z_n,v_n)$ in $\hat C$ that depends on $n$ but not on $(X,v_X)$, 
and similarly on morphisms. 
\end{enumerate}
Here $\oplus$ refers to the additive structure of $C$ and $\otimes$ to the monoidal structure. }
\end{defn}

\smallskip

\begin{rem}\label{remBCCat}{\rm 
In all the explicit cases considered in \cite{ManMar2} and in this paper, 
the endofunctors $\sigma_n$ and $\tilde\rho_n$ of Definition~\ref{BCcat} have the form
$$ \sigma_n (X,v_X)= (X,v_X \circ \sigma_n) \ \ \  \text{ and } \ \ \ 
\tilde\rho_n (X,v_X)= (X\times Z_n, \Phi_n(v_X)), $$
where the endomorphism $v_X$ is the action of a generator of some finite cyclic group
$\Z/N\Z$ quotient of $\hat\Z$ and the action satisfies $v_X\circ \sigma_n (\zeta,x)=v_X( \sigma_n(\zeta),x)$,
where $\sigma_n(\zeta)=\zeta^n$ is the Bost--Connes map of \eqref{eqA}, while
the action $\Phi_n(v_X)$ on $X\times Z_n$ is a geometric form of the Verschiebung, as
will be discussed more explicitly in \S \ref{VerschSec}. 
The object $(Z_n,v_n)$ in Definition~\ref{BCcat} plays the role of the element $n\pi_n$
in the integral Bost--Connes algebra and the relations \eqref{endofunctorBCrels} play the
role of the relations \eqref{sigmarhorels}.
}\end{rem}

\smallskip

This definition covers the main examples considered in \S\S~3 and 4 of \cite{ManMar2} 
obtained using the assembler categories associated to the equivariant Grothendieck
ring $K_0^{\hat \Z}(\cV)$ of varieties with a good $\hat\Z$-action factoring through some
finite cyclic quotient and to the equivariant version ${\rm Burn}^{\hat\Z}$ of the 
Kontsevich--Tschinkel Burnside ring. This same definition also accounts for the construction
we will discuss in \S \ref{TorifiedSec} of this paper, 
based on assembler categories associated to torified varieties (see Remark~\ref{BCliftTcat}). 

\smallskip

The more
general formulation given in Definition~\ref{BCcat} is motivated by the fact that one
expects other significant examples of categorical Bost--Connes structures where
the choice of the subcategory $\hat C$ of the automorphism category ${\rm Aut}(C)$
is not determined by the action of a cyclic group as in the cases discussed here.
Such more general classes of categorical Bost--Connes systems are not discussed
in the present paper, but they are a motivation for future work, for which we just set the
general framework in this section.

\smallskip

A generalization of Definition~\ref{BCcat} is needed when considering relative cases,
in particular the lift to assemblers of the 
construction presented in \S \ref{RelGrothSec} for
relative equivariant Grothendieck rings $K_0^{\hat\Z}(\cV_S)$.
The reason why we need the following modification of 
Definition~\ref{BCcat} is the fact that, in the relative setting, the base scheme 
$S$ itself has its enhancement structure (the group action, in the specific examples) 
modified by the endofunctors implementing
the Bost--Connes structure and this needs to be taken into account. We will see this
additional structure more explicitly applied in \S \ref{RelativeAssSec}, in the specific
case where the automorphisms are determined by a group action (see Remark~\ref{asshatZrem}).

\begin{defn}\label{BCcat2} {\rm
Let $\hat \cI$ be an enhancement of an additive (symmetric) monoidal category $\cI$
as above, endowed with a Bost--Connes system given by endofunctors $\{ \sigma^\cI_n \}$
and $\{ \tilde\rho_n^\cI \}$ of $\hat \cI$ as in Definition~\ref{BCcat}, with
$\alpha_n$ the object in $\hat\cI$ with $\tilde\rho_n \circ \sigma_n (\alpha)= \alpha \otimes \alpha_n$.
Let $\{  \hat C_\alpha \} _{\alpha\in \hat\cI}$ be a collection of enhancements of  
additive (symmetric) monoidal categories $C_\alpha$, indexed by the 
objects of the auxiliary category $\hat \cI$, endowed with functors 
$f_n: \hat C_{\alpha^{\oplus n}}\to \hat C_\alpha$
and $h_n: \hat C_{\alpha \times \beta} \to \hat C_\alpha$. 
Let $\{ \sigma_n \}_{n\in \N}$
and $\{ \tilde\rho_n \}_{n\in\N}$ be two collections of functors
$$ \sigma_n : \hat C_\alpha \to \hat C_{\sigma^\cI_n(\alpha)} \ \ \ \text{ and } \ \ \
 \tilde\rho_n : \hat C_\alpha \to \hat C_{\tilde\rho_n^\cI(\alpha)} $$
satisfying the properties:
\begin{enumerate}
\item The functors $\sigma_n$ are compatible with both the
additive and the (symmetric) monoidal structure, while the
functors $\tilde\rho_n$ are functors of additive categories.
\item For all $n,m\in \N$ these functors satisfy
$$ \sigma_{nm}=\sigma_n \circ \sigma_m, \ \ \ \  \tilde\rho_{nm}=\tilde\rho_n \circ \tilde\rho_m. $$
\item The compositions 
$$ \sigma_n\circ \tilde\rho_n : \hat C_\alpha \to \hat C_{\alpha^{\oplus n}} \ \ \  \text{ and } \ \ \ 
\tilde\rho_n \circ \sigma_n: \hat C_\alpha \to \hat C_{\alpha \otimes \alpha_n} $$
satisfy
\begin{equation}\label{endofunctorBC}
\begin{array}{rcl}
f_n\circ \sigma_n\circ \tilde\rho_n (X,v_X)_\alpha & = & (X,v_X)_\alpha^{\oplus n} \ \ \ \text{ and } \\[3mm]  
h_n\circ \tilde\rho_n \circ \sigma_n (X,v_X)_\alpha & = & (X,v_X)_\alpha \otimes (Z_n,v_n)_\alpha, 
\end{array}
\end{equation}
for an object $(Z_n,v_n)_\alpha$ in $\hat C_\alpha$ that depends on $n$ and $\alpha$, but not on $(X,v_X)$. 
\end{enumerate} }
\end{defn} 

\smallskip

We will first focus on the case of assembler categories, as those were at
the basis of our constructions of Bost--Connes systems in \cite{ManMar2}, but 
we will also consider in \S \ref{NoriSec} a different categorical setting that will allow us to
identify analogous structures at a motivic level, following the formalism of
geometric diagrams and Nori motives.

\smallskip
\subsection{Assemblers for the relative Grothendieck ring}\label{RelativeAssSec} 

We consider the relative Grothendieck ring $K_0(\cV_S)$ of varieties over a base variety 
$S$ over a field $\K$, as in Definition~\ref{RelGrK0}.

\smallskip

An assembler $\cC_S$ such that the associated spectrum $K(\cC_S)$ has
$K_0(\cC_S)=\pi_0 K(\cC_S)$ given by the relative Grothendieck ring $K_0(\cV_S)$
can be obtained as a slight modification of the construction given in \cite{Zak3} for the
ordinary Grothendieck ring $K_0(\cV_\K)$.

\smallskip

\begin{defn}\label{assK0Sdef} {\rm 
The assembler $\cC_S$ for the relative Grothendieck ring $K_0(\cV_S)$ has 
objects $f: X \to S$ that are varieties over $S$ and morphisms that are locally
closed embeddings of varieties over $S$. }
\end{defn}

\smallskip

\begin{lem}\label{assK0S}
The category $\cC_S$ of Definition~\ref{assK0Sdef}
is indeed as assembler, with the Grothendieck topology on $\cC_S$
is generated by the covering families $$\{ Y \hookrightarrow X, X\smallsetminus Y \hookrightarrow X \}$$ 
with compatible maps \eqref{overSmaps}
\begin{equation}\label{overSm}
 \xymatrix{ 
Y \ar@{^{(}->}[r] \ar[rd]_{f|_Y} & X \ar[d]^f & X\smallsetminus Y \ar@{_{(}->}[l] \ar[ld]^{f|_{X\smallsetminus Y}} \\
& S & 
} \end{equation}
\end{lem}

\proof The argument is the same as in \cite{Zak1}, \cite{Zak3} and in \cite{ManMar2}.
In this setting finite disjoint covering families are maps
$$ \xymatrix{ X_i \ar@{^{(}->}[r] \ar[rd]_{f_i} & X \ar[d]^f \\ & S } $$
where $X_i=Y_i\smallsetminus Y_{i-1}$ with commutative diagrams
$$ \xymatrix{ Y_0 \ar@{^{(}->}[r] \ar[rrrd]_{f_0} & Y_1 \ar@{^{(}->}[r] \ar[rrd]^{f_1} & \cdots  \ar@{^{(}->}[r] & Y_n=X
\ar[d]^f \\ & & & S } $$
The category has pullbacks, hence as shown in \cite{Zak1} (Remark after Definition~2.4) this suffices to obtain that
any two finite disjoint covering families have a common refinement. Morphisms
are embeddings compatible with the structure maps as in \eqref{overSm}
hence in particular monomorphisms. Theorem~2.3 of \cite{Zak1} then shows that
the spectrum $K(\cC_S)$ associated to this assembler category has
$\pi_0 K(\cC_S)=K_0(\cV_S)$.
\endproof

\smallskip

In a similar way we obtain an assembler category and spectrum for the
equivariant version $K_0^{\hat\Z}(\cV_S)$. The argument is as in the
previous case and in Lemma~4.5.1 of \cite{ManMar2}, using the
inclusion-exclusion relations \eqref{overSmaps2}.

\begin{cor}\label{assK0ShatZ}
An assembler category $\cC^{\hat\Z}_{(S,\alpha)}$ for $K_0^{\hat\Z}(\cV_{(S,\alpha)})$ is constructed 
as in Lemma~\ref{assK0S} with objects the $\hat\Z$-equivariant $f: X \to S$, morphisms given by
$\hat\Z$-equivariant locally closed embeddings of varieties over $S$
and with Grothendieck topology generated by the covering families given by
$\hat\Z$-equivariant maps as in \eqref{overSmaps} and \eqref{overSmaps2}.
\end{cor}

\smallskip

As in Proposition~4.2 of \cite{ManMar2}, we show that the lifting of the integral Bost--Connes algebra
obtained in Proposition~\ref{BCliftSZ} and Theorem~\ref{liftBCmaps} further lifts to 
functors of the associated assembler categories,  with the $\sigma_n$ compatible with 
the monoidal structure, but not the $\tilde\rho_n$. 

\begin{thm}\label{liftAssSZ}
The maps $\sigma_n: (f: (X,\alpha_X) \to (S,\alpha)) \mapsto (f: (X,\alpha_X\circ \sigma_n) \to (S,\alpha\circ \sigma_n))$ and $\tilde\rho_n: (f: (X,\alpha_X) \to (S,\alpha)) \mapsto (f\times {\rm id}: (X\times Z_n, \Phi_n(\alpha_X)) \to (S\times Z_n, \Phi_n(\alpha)))$ define functors of the assembler categories 
$\sigma_n : \cC^{\hat\Z}_{(S,\alpha)} \to \cC^{\hat\Z}_{(S,\alpha\circ \sigma_n)}$ and
$\tilde \rho_n: \cC^{\hat\Z}_{(S,\alpha)} \to \cC^{\hat\Z}_{(S\times Z_n , \Phi_n(\alpha))}$. The
functors $\sigma_n$ are compatible with the monoidal structure.
\end{thm}

\proof The functors $\sigma_n$ defined as above on objects are compatibly defined on
morphisms by assigning to a locally closed embedding 
$$ \sigma_n : \ \ \ \xymatrix{ (Y,\alpha_Y) \ar[r]^{j} \ar[dr]_{f_Y} & (X,\alpha_X) \ar[d]^{f_X} \\ & (S,\alpha)
} \ \ \ \mapsto  \ \ \ 
\xymatrix{ (Y,\alpha_Y\circ \sigma_n) \ar[r]^{j} \ar[dr]_{f_Y} & (X,\alpha_X\circ \sigma_n) \ar[d]^{f_X} \\ & (S,\alpha\circ \sigma_n) }
$$
Similarly, we define the $\tilde\rho_n$ on morphisms by
$$ \tilde\rho_n : \ \ \ \xymatrix{ (Y,\alpha_Y) \ar[r]^{j} \ar[dr]_{f_Y} & (X,\alpha_X) \ar[d]^{f_X} \\ & (S,\alpha)
} \ \ \ \mapsto  \ \ \ 
\xymatrix{ (Y \times Z_n,\Phi_n(\alpha_Y)) \ar[r]^{j} \ar[dr]_{f_Y} & (X\times Z_n, \Phi_n(\alpha_X)) 
\ar[d]^{f_X} \\ & (S \times Z_n, \Phi_n(\alpha)) }
$$
The functors $\sigma_n$ are compatible with the monoidal structure since
$\sigma_n (X,\alpha_X) \times \sigma_n(X',\alpha_{X'})=(X\times X', (\alpha\times \alpha')\circ \Delta \circ \sigma_n)=\sigma_n ((X,\alpha_X) \times (X',\alpha_{X'}))$. 
\endproof

The functor of assembler categories determines an induced map of spectra and in turn
an induced map of homotopy groups. By construction the induced maps on the $\pi_0$
homotopy agree with the maps \eqref{sigmaSn} and \eqref{rhoSn} of Proposition~\ref{BCliftSZ}.
\endproof

\smallskip

\begin{rem}\label{asshatZrem}{\rm  We can associate to the assembler category 
$\cC^{\hat\Z}_{(S,\alpha)}$ of Corollary~\ref{assK0ShatZ} with the endofunctors $\sigma_n$
and $\tilde\rho_n$ a categorical Bost--Connes structure in the sence of Definition~\ref{BCcat2},
where the objects are $f: X \to S$ as above with the automorphisms given by elements 
$g\in\hat\Z$ acting on $f: X \to S$ through the action by $\alpha_X(g)$ on $X$ and by
$\alpha_S(g)$ on $S$, intertwined by the equivariant map $f$. }
\end{rem}

\medskip
\section{Torifications, $\F_1$-points, zeta functions, and spectra}\label{TorifiedSec}   

In this section we relate the point of view developed in \cite{ManMar2}, with lifts of the
Bost--Connes system to Grothendieck rings and spectra, to
the approach to $\F_1$-geometry based on torifications. This was first introduced
in \cite{LoLo}. Weaker forms of torification were also considered in \cite{ManMar},
which allow for the development of a form of $\F_1$-geometry suitable for the
treatment of certain classical moduli spaces. 

\smallskip

The approach we follow here, in order the relate the case of torified geometry
with the Bost--Connes systems on Grothendieck rings, assemblers, and spectra
discussed in \cite{ManMar2}, is based on the following simple setting.
Instead of working with the equivariant Grothendieck rings $K_0^{\hat\Z}(\cV)$ and
$K_0^{\hat\Z}(\cV_S)$, where one assumes the varieties have a good effectively finite $\hat\Z$-action, 
we consider here a variant that connects to
the torifications point of view on $\F_1$-geometry of \cite{LoLo}. We replace varieties
with $\Z/N\Z$-effectively finite $\hat\Z$-actions with varieties with a 
$\Q/\Z$-action induced by a torification, where the group schemes 
$\fm_n$ of $n$-th roots of unity, given by the kernels 
$$ 1\to  \fm_n  \to \bG_m \stackrel{\lambda\mapsto \lambda^n}{\longrightarrow} \bG_m \to 1 $$
determine a diagonal embedding in each torus and an action by multiplication.
This is a very restrictive class of varieties, because the existence of a torification on a variety 
implies that the Grothendieck class is a sum of classes of tori with non-negative coefficients.
The resulting construction will be more restrictive than the one
considered in \cite{ManMar2}. We will see, however, that one can still see in this context
several interesting phenomena, especially in connection with the ``dynamical" approach
to $\F_1$-geometry proposed in \cite{ManMar2}.

\smallskip
\subsection{Torifications} \label{TorifSec}  

A torification of an algebraic variety $X$ defined over $\Z$ is a decomposition $X=\sqcup_{i\in \cI} T_i$
into algebraic tori $T_i=\bG_m^{d_i}$. Weaker to stronger forms of torification \cite{ManMar} 
include
\begin{enumerate}
\item {\em torification of the Grothendieck class}:
$[X]=\sum_{i\in \cI} (\bL-1)^{d_i}$ with $\bL$ the Lefschetz motive; 
\item {\em geometric torification}: $X=\sqcup_{i\in \cI} T_i$ with $T_i=\bG_m^{d_i}$;
\item {\em affine torification}: the existence of an affine covering compatible with the geometric 
torification, \cite{LoLo};
\item {\em regular torification}: the closure of each torus in the geometric torification is also
a union of tori of the torification, \cite{LoLo}.
\end{enumerate}

\smallskip

Similarly, there are different possibilities when one considers morphisms of torified 
varieties, see \cite{ManMar}. In view of describing associated Grothendieck rings,
we review the different notions of morphisms. The Grothendieck classes are then
defined with respect to the corresponding type of isomorphism. 

A torified morphism of geometric torifications in the sense of \cite{LoLo} between torified
varieties $f: (X,T) \to (Y,T')$ is a morphism $f: X \to Y$ of varieties
together with a map $h: I \to J$ of the indexing sets of
the torifications $X=\sqcup_{i\in I} T_i$ and $Y=\sqcup_{j\in J} T'_j$
such that the restriction of $f$ to tori $T_i$ is a morphism of
algebraic groups $f_i : T_i \to T'_{h(i)}$. There are then three different classes
of morphisms of torified varieties that were introduced in \cite{ManMar}: strong, ordinary,
and weak morphisms. To describe them, one first defines strong, ordinary, and weak
equivalences of torifications, and one then uses
these to define the respective class of morphisms. 

\smallskip

Let $T$ and $T'$ be two geometric torifications of a variety $X$. 
\begin{enumerate}
\item The torifications $(X,T)$ and $(X,T')$ are {\em strongly equivalent} if the identity map on $X$ is
a torified morphism as above.
\item The torifications $(X,T)$ and $(X,T')$ are {\em ordinarily equivalent} if there exists an
automorphism $\phi: X \to X$ that is a torified morphism.
\item The torifications $(X,T)$ and $(X,T')$ are {\em weakly equivalent} if $X$ has two decompositions
$X=\cup_i X_i$ and $X=\cup_j X'_j$ into a disjoint union of subvarieties, compatible with the torifications,
such that there are isomorphisms of varieties $\phi_i: X_i \to X'_{j(i)}$ that are torified. 
\end{enumerate}

In the weak case a ``decomposition compatible with torifications" means that the intersections 
$T_i \cap X_j$ of the tori of $T$ with the pieces of the decomposition (when non-empty) are 
tori of the torification of $X_j$, and similarly for $T'_i \cap X'_j$.
In general weakly equivalent torification are not ordinarily equivalent because the maps 
$\phi_i$ need not glue together to define a single map $\phi$ on all of $X$. 

\smallskip

We then have the following classes of morphisms of torified varieties from \cite{ManMar}:
\begin{enumerate}
\item {\em strong morphisms}: these are torified morphisms in the sense of \cite{LoLo},
namely morphisms that restrict to morphisms of tori of the respective torifications. 
\item {\em ordinary morphisms}: an ordinary morphism of torified varieties  $(X,T)$ and $(Y,T')$
is a morphism $f: X \to Y$ such that becomes a torified morphism after composing with strong isomorphisms,
that is, $\phi_Y\circ f \circ \phi_X: (X,T)\to (Y,T')$ is a strong morphism of torified varieties,
for some isomorphisms $\phi_X: X\to X$ and $\phi_Y: Y \to Y$. In other words, if we denote by $T_\phi$
and $T'_\phi$ the torifications such that $\phi_X: (X,T)\to (X,T_\phi)$ and $\phi_Y: (Y,T'_\phi)\to (Y,T')$
are torified, then $f: (X,T_\phi) \to (Y,T'_\phi)$ is torified. 
\item {\em weak morphisms}: the torified varieties $(X,T)$ and $(Y,T')$ admit decompositions 
$X=\sqcup_i X_i$ and $Y=\sqcup_j Y_i$, compatible with the torifications, such that there exist 
ordinary morphisms $f_i: (X_i,T_i) \to (Y_{f(i)}, T'_{f(i)})$ of these subvarieties.
\end{enumerate}
Note that the strong isomorphisms $\phi_X: (X,T)\to (X,T_\phi)$ and $\phi_Y: (Y,T'_\phi)\to (Y,T')$
used in the definition of ordinary morphisms are ordinary equivalences of the torifications  $T$ and $T_\phi$,
respectively $T'$ and $T'_\phi$. 

\smallskip

Given these notions of morphisms, we can correspondingly 
construct Grothendieck rings $K_0(\cT)^s$, $K_0(\cT)^o$,
and $K_0(\cT)^w$ in the following way. 

\smallskip

As an abelian group $K_0(\cT)^s$ is generated by isomorphism
classes $[X,T]_s$ of pairs of a torifiable variety $X$ and a torification $T$ modulo strong
isomorphisms, with the inclusion-exclusion relations $[X,T]_s=[Y,T_Y]_s +[X\smallsetminus Y, T_{X\smallsetminus Y}]_s$ whenever $(Y, T_Y)\hookrightarrow (X,T)$ is a strong morphism
(that is, the inclusion of $Y$ in $X$ is compatible with the torification: $Y$ is a union of
tori of the torification of $X$) and $(Y,T_Y)$ is a {\em complemented subvariety} in $(X,T)$,
which means that the complement $X\smallsetminus Y$ is also a union of tori of the torification
so that the inclusion of $(X\smallsetminus Y, T_{X\smallsetminus Y})$ in $(X,T)$ is also a strong morphism. This 
complemented condition is very strong. Indeed, one can see that, for example, there are
in general very few complemented points in a torified variety. The product operation is
$[X,T]_s \cdot [Y,T']_s=[X\times Y, T\times T']_s$ with the torification $T\times T'$ given
by the product tori $T_{ij}=T_i\times T'_j=\bG_m^{d_i+d_j}$.

\smallskip

The abelian group $K_0(\cT)^o$ is generated by isomorphism classes $[X]_o$ varieties that admit 
a torification with respect to ordinary isomorphisms, with the inclusion-exclusion relations
$[X]_o=[Y]_o+[X\smallsetminus Y]_o$ whenever the inclusions $Y\hookrightarrow X$
and $X\smallsetminus Y \hookrightarrow X$ are ordinary morphisms. The product is
the class of the Cartesian product $[X]_o \cdot [Y]_o =[X\times Y]_o$. 

\smallskip

The abelian group $K_0(\cT)^w$ is generated by the isomorphism classes $[X]_w$
of torifiable varieties $X$ with respect to weak morphisms, with the inclusion-exclusion
relations $[X]_w=[Y]_w +[X\smallsetminus Y]_w$ whenever the inclusions $Y\hookrightarrow X$
and $X\smallsetminus Y \hookrightarrow X$ are weak morphisms. The product structure
is again given by $[X]_w\cdot [Y]_w=[X\times Y]_w$.

\smallskip

The reader can consult the explicit examples given in \cite{ManMar} to see how
these notions (and the resulting Grothendieck rings) can be different. For example,
as mentioned in \S 2.2 of \cite{ManMar}, consider the variety $X=\P^1\times \P^1$
and consider on it two torifications $T$ and $T'$, where $T$ is the standard
torification given by the decomposition of each $\P^1$ into cells $\A^0\cup \A^1$, with the cell
$\A^1$ torified as $\A^0\cup \bG_m$, while $T'$ is the torification where in the big cell $\A^2$ of
$\P^1\times \P^1$ we take a torification of the diagonal $\A^1$ and a torification of the 
complement of the diagonal, and we use the same torification of the lower dimensional 
cells as in $T$. These two torifications are related by a weak isomorphism, 
hence the elements $(\P^1\times \P^1,T)$
and $(\P^1\times \P^1,T')$ define the same class in $K_0(\cT)^w$, but they are not
related by an ordinary isomorphism so they define different classes in $K_0(\cT)^o$.

\smallskip

Note however that, in all these cases, the Grothendieck classes $[X]_a$ with $a=s,o,w$
have the form $[X]_a=\sum_{n\geq 0} a_n \bT^n$ with $a_n\in \Z_+$ and
$\bT^n=[\bG_m^n]$. 

\smallskip

In the following, whenever we simply write $a =  s,o,w $ without specifying one of the
three choices of morphisms, it means that the stated property holds for all of these
choices. 

\smallskip
\subsubsection{Relative case} \label{TorifRelSec}  

In a similar way, we can construct relative Grothendieck rings $K_S(\cT)^a$
with $a=s,o,w$ where in the strong case $S=(S,T_S)$ is a choice of a 
variety with an assigned torification, with $K_S(\cT)^s$ generated as an
abelian group by isomorphisms classes $[f: (X,T) \to (S,T_S)]$ where $f$ is a strong
morphism of torified varieties and the isomorphism class is taken with respect
to strong isomorphisms $\phi, \phi_S$ such that the diagram commutes
$$  \xymatrix{ (X,T) \ar[r]^{\phi}\ar[d]_f & (X',T') \ar[d]^{f'} \\
(S,T_S) \ar[r]^{\phi_S} & (S,T_S) } $$
with inclusion-exclusion relations
$$ [f: (X,T) \to (S,T_S)] = $$ $$ [f|_{(Y,T_Y)}: (Y,T_Y)\to (S,T_S)] + [
f|_{(X\smallsetminus Y, T_{X\smallsetminus Y})} : 
(X\smallsetminus Y, T_{X\smallsetminus Y}) \to (S,T_S)] $$
where $\iota_Y: (Y,T_Y) \hookrightarrow (X,T)$ is a strong morphism and
$(Y,T_Y)$ is complemented with $\iota_{X\smallsetminus Y}: (X\smallsetminus Y, T_{X\smallsetminus Y})\hookrightarrow (X,T)$
also a strong morphism and both these inclusions are compatible with the map
$f:(X,T) \to (S,T_S)$, so that $f_Y=f\circ \iota_Y$ and 
$f|_{(X\smallsetminus Y, T_{X\smallsetminus Y})} = f\circ \iota _{X\smallsetminus Y}$
are strong morphisms. The construction for ordinary and weak morphism is similar,
with the appropriate changes in the definition.

\smallskip
\subsection{Group actions}\label{QZactSec}  

In order to operate on Grothendieck classes with Bost--Connes type
endomorphisms, we introduce appropriate group actions.  

\smallskip

Torified varieties carry natural $\Q/\Z$ actions,
since the roots of unity embed diagonally in each torus of the torification
and act on it by multiplication. However, we will also be interested in
considering good effectively finite $\hat\Z$-actions, in the sense already discussed in
\cite{ManMar2}, that is, actions of $\hat\Z$ as in Definition~\ref{efinact}. 

\smallskip

\begin{rem}\label{hatZvsQZ}{\rm 
The main reason for working with $\hat\Z$-actions rather than
with $\Q/\Z$ actions lies in the fact that, in the construction of the
geometric Vershiebung action discussed in \S \ref{VerschSec} we need
to be able to describe the cyclic permutation action of $\Z/n\Z$ on the
finite set $Z_n$ as an action factoring through $\Z/n\Z$. This cannot
be done in the case of $\Q/\Z$-actions because 
there are no nontrivial group homomorphisms $\Q/\Z \to \Z/n\Z$
since $\Q/\Z$ is infinitely divisible.
}\end{rem}

\smallskip

In the case of the natural $\Q/\Z$-actions on torifications, 
we consider objects of the form $(X,T,\alpha)$ where $X$ is a torifiable
variety, $T$ a choice of a torification, and $\alpha: \Q/\Z \times X \to X$ 
an action of $\Q/\Z$ determined by an embedding of $\Q/\Z$ as roots of
unity in $\bG_m(\C)=\C^*$, which act on each torus $T_i=\bG_m^{k_i}$
diagonally by multiplication. An embedding of $\Q/\Z$ in $\bG_m$ is
determined by an invertible element in $\Hom(\Q/\Z,\bG_m)=\hat\Z$,
hence the action $\alpha$ is uniquely determined by the torification $T$
and by the choice of an element in $\hat\Z^*$. 

\smallskip

The corresponding morphisms are, respectively, strong, ordinary, or
weak morphisms of torified varieties compatible with the $\Q/\Z$-actions,
in the sense that the resulting torified morphism (after composing
with isomorphisms or with local isomorphisms in the ordinary and
weak case) are $\Q/\Z$-equivariant.  
We can then proceed as above and obtain equivariant Grothendieck rings
$K_0^{\Q/\Z}(\cT)^s$, $K_0^{\Q/\Z}(\cT)^o$,
and $K_0^{\Q/\Z}(\cT)^w$ of torified varieties. 

\smallskip

In the case of good $\Z/N\Z$-effectively finite $\hat\Z$-actions, the setting is essentially 
the same. We consider objects of the form $(X,T,\alpha)$ where $X$ is a torifiable
variety, $T$ a choice of a torification, and $\alpha: \Z/N\Z \times X \to X$ is given
by the action of the $N$-th roots of unity on the tori $T_i=\bG_m^{k_i}$ by
multiplication. Thus, a good $\hat\Z$-action is determined by $T$, by the
choice of an embedding of roots of unity in $\bG_m$ (an element of $\hat\Z^*$)
as above, and by the choice of $N\in \N$ that determines which subgroup of
roots of unity is acting. 

\smallskip

This choice of good $\Z/N\Z$-effectively finite $\hat\Z$-actions, with strong, ordinary, or
weak morphisms whose associated torified morphisms are 
$\Z/N\Z$-equivariant, determine equivariant Grothendieck rings
$K_0^{\hat\Z}(\cT)^s$, $K_0^{\hat\Z}(\cT)^o$,
and $K_0^{\hat\Z}(\cT)^w$ of torified varieties with good effectively finite $\hat\Z$-actions.

\smallskip
\subsection{Assembler and spectrum of torified varieties} \label{AssTorifiedSec}   

As in the previous cases of $K_0^{\hat Z}(\cV)$ of \cite{ManMar2} and in the case of $K_0^{\hat \Z}(\cV_S)$
discussed above, we consider the Grothendieck rings $K_0(\cT)^s$, $K_0(\cT)^o$,
and $K_0(\cT)^w$ and their corresponding equivariant versions 
$K_0^{\Q/\Z}(\cT)^s$, $K_0^{\Q/\Z}(\cT)^o$, $K_0^{\Q/\Z}(\cT)^w$, and
$K_0^{\hat\Z}(\cT)^s$, $K_0^{\hat\Z}(\cT)^o$, $K_0^{\hat\Z}(\cT)^w$
from the point of view of assemblers and spectra developed in
\cite{Zak1}, \cite{Zak2}, \cite{Zak3}.

\smallskip

\begin{prop}\label{assK0Ta}
For $a=s,o,w$, the category $\cC^a_\cT$ has objects that are pairs $(X,T)$ of a torifiable
variety and a torification, with morphisms the locally closed embeddings that are,
respectively, strong, ordinary, or weak morphisms of torified varieties.
The Grothendieck topology is generated by the covering families
\begin{equation}\label{covfamT}
 \{ (Y,T_Y) \hookrightarrow (X,T_X), (X\smallsetminus Y,T_{X\smallsetminus Y}) \hookrightarrow (X,T_X) \}
\end{equation} 
 where both embeddings are strong, ordinary, or weak morphisms, respectively. The category $\cC^a_\cT$ is an
assembler with spectrum $K(\cC^a_\cT)$ satisfying $\pi_0 K(\cC^a_\cT)=K_0(\cT)^a$. Similarly, for $G=\Q/\Z$
or $G=\hat\Z$ let $\cC^{G,a}_\cT$ be the category with objects $(X,T,\alpha)$ given by 
a torifiable variety $X$ with a torification $T$ and a $G$-action $\alpha$ of the kind discussed in \S \ref{QZactSec} and morphisms the 
locally closed embeddings that are $G$-equivariant strong, ordinary, or weak morphisms.
The Grothendieck topology is generated by covering families \eqref{covfamT} with $G$-equivariant
embeddings. The category $\cC^{G,a}_\cT$ is also an assembler, whose associated
spectrum $K(\cC^{G,a}_\cT)$ satisfies $\pi_0 K(\cC^{G,a}_\cT)=K_0^G(\cT)^a$. 
\end{prop}

\proof The argument is again as in \cite{Zak1}, see Lemma~\ref{assK0S}.
We check that the category admits pullbacks. In the strong case, if $(Y,T_Y)$ and
$(Y',T_{Y'})$ are objects with morphisms $f: (Y,T_Y) \hookrightarrow (X,T_X)$
and $f': (Y',T_{Y'}) \hookrightarrow (X,T_X)$ given by embeddings that are strong 
morphisms of torified varieties. This means that the tori of the torification $T_Y$
are restrictions to $Y$ of tori of the torification $T_X$ of $X$. Thus, both $Y$ and $Y'$
are unions of subcollections of tori of $T_X$. Their intersection $Y\cap Y'$ will then also
inherit a torification consisting of a subcollection of tori of $T_X$ and the resulting
embedding $(Y\cap Y', T_{Y\cap Y'})\hookrightarrow (X,T_X)$ is a strong morphism
of torified varieties.  
In the ordinary case, we consider embeddings $f: Y \hookrightarrow X$ 
and $f': Y' \hookrightarrow X$ that are ordinary morphisms of torified
varieties, which means that, for isomorphisms 
$\phi_X$, $\phi'_X$, $\phi_Y$, $\phi_{Y'}$, the compositions 
$$ \phi_X \circ f \circ \phi_Y: (Y,T_Y) \hookrightarrow (X,T_X), \ \ \ \ \ 
\phi'_X\circ f'\circ \phi_{Y'}: (Y',T_{Y'}) \hookrightarrow (X,T_X) $$
are (strong) torified morphisms. Thus, the tori of the torifications $T_Y$
and $T_{Y'}$ are subcollections of tori of $X$, under the embeddings
$\phi_X \circ f \circ \phi_Y$ and $\phi'_X\circ f'\circ \phi_{Y'}$. The intersection
$\phi_X \circ f \circ \phi_Y (Y)\cap \phi'_X\circ f'\circ \phi_{Y'}(Y') \subset X$
is isomorphic to a copy of $Y\cap Y'$ and has an induced torification $T_{Y\cap Y'}$ 
by a subcollection of tori of $T_X$. The embedding of $Y\cap Y'$ in $X$ with
this image is an ordinary morphism with respect to this torification. 
The weak case is constructed similarly to the ordinary case on the pieces of
the decomposition. The equivariant cases are constructed analogously, as
discussed in the case of equivariant Grothendieck rings of varieties in \cite{ManMar2}.
\endproof

\smallskip

\begin{cor}\label{AssTmaps}
There are inclusions of assemblers $\cC^s_\cT \hookrightarrow \cC^o_\cT \hookrightarrow  \cC^w_\cT$
that induce maps of $K$-theory, in particular $K_0(\cT)^s \to K_0(\cT)^o$ and $K_0(\cT)^o\to K_0(\cT)^w$.
Similarly, for the $G$-equivariant cases of Proposition~\ref{assK0Ta}.
\end{cor}

\proof Since for morphisms strong implies ordinary and ordinary implies weak,
one obtains inclusions of assemblers as stated. 
\endproof

\smallskip
\subsection{Lifting of the Bost--Connes system for torifications} \label{BCTorifSec}   

We consider here lifts of the integral Bost--Connes algebra to the
Grothendieck rings $K_0^{\hat\Z}(\cT)^s$, $K_0^{\hat\Z}(\cT)^o$,
and $K_0^{\hat\Z}(\cT)^w$ and to the assemblers and spectra 
$K^{\hat\Z}(\cC^s_\cT)$, $K^{\hat\Z}(\cC^o_\cT)$, and $K^{\hat\Z}(\cC^w_\cT)$. 

\smallskip

\begin{defn}\label{BCliftTdef}{\rm
We regard the zero-dimensional variety $Z_n$ as a torified
variety with the torification consisting of $n$ zero dimensional tori and with a
good $\hat\Z$ action factoring through $\Z/n\Z$ that cyclically permutes the
points of $Z_n$. We write $(Z_n,T_0,\gamma)$ for this object. For $(X,T,\alpha)$ a triple of 
a torifiable variety $X$, a given torification $T$, and an effectively finite action $\alpha$ 
of $\hat\Z$, we then set, for all $n\in\N$,
\begin{equation}\label{sigmarhonT}
\sigma_n(X,T,\alpha)=(X,T,\alpha\circ \sigma_n) \ \ \ \text{ and } \ \ \ \tilde\rho_n(X,T,\alpha)=(X\times Z_n, \sqcup_{a\in Z_n} T, \Phi_n(\alpha)) \, .
\end{equation}
}\end{defn}

\smallskip 

\begin{prop}\label{BCliftT}
The $\sigma_n$ and $\tilde\rho_n$ defined as in \eqref{sigmarhonT} 
determine endofunctors of the assembler categories $\cC^{\hat\Z,a}_\cT$ that
induce, respectively, ring homomorphisms $\sigma_n: K^{\hat\Z}(\cC^a_\cT) \to
K^{\hat\Z}(\cC^a_\cT)$ and group homomorphisms $\tilde\rho_n: K^{\hat\Z}(\cC^a_\cT) \to
K^{\hat\Z}(\cC^a_\cT)$ with the Bost--Connes relations 
$$ \tilde\rho_n \circ \sigma_n (X,T,\alpha) =(X,T, \alpha)\times (Z_n, T_0, \gamma)
\ \ \ \ \  \sigma_n\circ \tilde\rho_n (X,T,\alpha) = (X,T,\alpha)^{\oplus n}.$$ 
\end{prop}

\proof The proof is completely analogous to the case discussed in Theorem~\ref{liftAssSZ}
and to the similar cases discussed in \cite{ManMar2}.
\endproof

\begin{rem}\label{BCliftTcat} {\rm The $\sigma_n$ and $\tilde\rho_n$ defined as in \eqref{sigmarhonT} 
determine a categorical Bost--Connes system as in Definition~\ref{BCcat},
where the objects are pairs $(X,T)$ and the automorphisms are elements $g\in \hat\Z$
acting through the effectively finite action $\alpha(g)$.
}\end{rem}

\begin{rem}\label{endotoric} {\rm 
Bost--Connes type quantum statistical mechanical systems associated to
individual toric varieties (and more generally to varieties admitting torifications)
were constructed in \cite{JinMar}. Here instead of Bost--Connes endomorphisms
of individual varieties we are interested in a Bost--Connes system over the
entire Grothendieck ring and its associated spectrum.}
\end{rem}

\begin{rem}\label{multiBC}{\rm 
Variants of the construction above can be obtained by considering the
multivariable versions of the Bost--Connes system discussed in \cite{Mar},
with actions of subsemigroups of $M_N(\Z)^+$ on $\Q[\Q/\Z]^{\otimes N}$, that is,
subalgebras of the crossed product algebra
$$ \Q[\Q/\Z]^{\otimes N}\rtimes_\rho M_N(\Z)^+  $$
generated by $e(\underline{r})$ and $\mu_\alpha$, $\mu_\alpha^*$ with
$$ \rho_\alpha(e(\underline{r}))= \mu_\alpha e(\underline{r}) \mu_\alpha^*=
\frac{1}{\det\alpha} \sum_{\alpha(\underline{s})=\underline{r}} e(\underline{s}) $$
$$ \sigma_\alpha(e(\underline{r}))=\mu_\alpha^* e(\underline{r}) \mu_\alpha =
e(\alpha(\underline{r})) .$$
The relevance of this more general setting to $\F_1$-geometries lies in a
result of Borger and de Smit \cite{BorgerdeSmit} showing that every 
torsion free finite rank $\Lambda$-ring embeds in some $\Z[\Q/\Z]^{\otimes N}$
with the action of $\N$ determined by the $\Lambda$-ring structure compatible 
with the diagonal subsemigroup of $M_N(\Z)^+$. 
}\end{rem}

\medskip
\section{Torified varieties and zeta functions}\label{TorifiedZetaSec}   

We discuss in this section the connection between the dynamical
point of view on $\F_1$-geometry proposed in \cite{ManMar2} and
the point of view based on torifications. 

\smallskip

We first discuss in \S \ref{CountF1Sec} and \S \ref{BBsec} the notion of $\F_1$-points 
of a torified variety and its relation to the torification of the Grothendieck class, with some
explicit examples. We then introduce the $\F_1$-zeta function in \S \ref{CountF1zetaSec}
and we show its main properties in Proposition~\ref{zetaF1motmeas}, 
while in \S \ref{F1HWSec} we explain the relation
between the $\F_1$-zeta function and the Hasse--Weil zeta function.

\smallskip

In \S \ref{DynZetaSec} we consider the point of view on $\F_1$-structures proposed
in \cite{ManMar2} based on dynamical systems inducing quasi-uniponent endomorphisms
on homology, in the particular case of torified varieties with dynamical systems compatible
with the torification. We focus on the associated dynamical zeta functions, the Lefschetz
zeta function and the Artin--Mazur zeta function, whose properties we 
recall in \S \ref{DynZetaProSec}. We then prove in Proposition~\ref{zetadynmotmeas}
that the resulting dynamical zeta function have similar properties to the $\F_1$-zeta
function in the sense that both define exponentiable motivic measures from the
Grothendieck rings of torified varieties to the Witt ring. 

\smallskip
\subsection{Counting $\F_1$-points}\label{CountF1Sec}  

Assuming that a variety $X$ over $\Z$ admits an $\F_1$-structure, regarded here as one
of several possible forms of torified structure recalled above, \cite{LoLo}, \cite{ManMar},  
the number of points of $X$ over $\F_1$ is computed as the $q\to 1$ limit of the counting function
$N_X(q)$ of points over $\F_q$ of the mod $p$ reduction of $X$, for $q$ a power of $p$.
Any form of torified structure in particular implies that the variety is polynomially countable,
hence that the counting function $N_X(q)$ is a polynomial in $q$ with $\Z$-coefficients. 
The limit $\lim_{q\to 1} N_X(q)$, possibly normalized by a power of $q-1$, is interpreted 
as the number of $\F_1$-points of $X$, see \cite{Soule}. Similarly, one can define
``extensions" $\F_{1^m}$ of $\F_1$, in the sense of \cite{KapSmi} (see also \cite{CCM}).
These corresponds to actions of the groups $\fm_m$ of $m$-th roots of unity. In terms of
a torified structure, the points over $\F_{1^m}$ count $m$-th roots of unity in each torus 
of the decomposition. In terms of the counting function $N_X(q)$ the counting of points
of $X$ over the extension  $\F_{1^m}$ is obtained as the value $N_X(m+1)$, see
Theorem~4.10 of \cite{CoCo} and Theorem~1 of \cite{Deit}). Summarizing, we have
the following.

\begin{lem}\label{lemF1mpts}
Let $X$ be a variety over $\Z$ with torified Grothendieck class 
\begin{equation}\label{torclass}
 [X] = \sum_{i=0}^N a_i \bT^i 
\end{equation} 
with coefficients $a_i \in \Z_+$ and $\bT=[\bG_m]=\bL-1$.
Then the number of points over $\F_{1^m}$ of $X$ is given by
\begin{equation}\label{F1mpts}
\# X(\F_{1^m}) = \sum_{i=0}^N a_i \, m^i . 
\end{equation}
In particular, $\# X (\F_1)=a_0=\chi(X)$ the Euler characteristic. 
\end{lem}

\smallskip
\subsection{Bia\/{l}ynicki-Birula decompositions and torified geometries} \label{BBsec}   

As shown in \cite{Bano}, \cite{Brosnan}, the motive of a smooth projective 
variety with action of the multiplicative group
admits a decomposition, obtained via the method of Bia\/{l}ynicki-Birula,
\cite{BiBi1}, \cite{BiBi2}, \cite{BiBi3}. We recall the result here, in a particular case 
which gives rise to examples of torified varieties. 

\smallskip

\begin{lem}\label{torifBB}
Let $X$ be a smooth projective $k$-variety $X$ endowed with a $\bG_m$ action
such that the fixed point locus $X^{\bG_m}$ admits a torification of the Grothendieck class.
Then $X$ also admits a torification of the Grothendieck class. Consider the filtration
$X=X_n\supset X_{n-1}\supset\cdots\supset X_0\supset \emptyset$ with affine
fibrations $\phi_i: X_i \smallsetminus X_{i-1} \to Z_i$ over the components 
$X^{\bG_m}=\sqcup_i Z_i$, associated to the Bia\/{l}ynicki-Birula decomposition.
If the fixed point locus $X^{\bG_m}$ admits a geometric torification such that
the restrictions of the fibrations $\phi_i$ to the individual tori of the torification of 
$Z_i$ are trivializable, then $X$ also admits a geometric torification.
\end{lem}

\proof
The Bia\/{l}ynicki-Birula decomposition, \cite{BiBi1},
\cite{BiBi2}, \cite{BiBi3}, see also \cite{Hessel}, shows that 
a smooth projective $k$-variety $X$ endowed with a $\bG_m$ action
has smooth closed fixed point locus $X^{\bG_m}$ which decomposes
into a finite union of components $X^{\bG_m}=\sqcup_i Z_i$, of dimensions
$\dim Z_i$ the dimension of $TX_z^{\bG_m}$ at $z\in Z_i$.
The variety $X$ has a filtration $X=X_n\supset X_{n-1}\supset\cdots\supset X_0\supset \emptyset$
with affine fibrations $\phi_i: X_i \smallsetminus X_{i-1} \to Z_i$ of relative
dimension $d_i$ equal to the dimension of the positive eigenspace of the $\bG_m$-action on
the tangent space of $X$ at points of $Z_i$. The scheme $X_i \smallsetminus X_{i-1}$ is
identified with $\{ x\in X \,:\, \lim_{t\to 0} tx \in Z_i \}$ under the $\bG_m$-action $t: x\mapsto tx$,
with $\phi_i(x)=\lim_{t\to 0} tx$. As shown in \cite{Brosnan}, the object $M(X)$
in the category of correspondences ${\rm Corr}_k$ with integral coefficients (and in the
category of Chow motives) decomposes as
\begin{equation}\label{BBdecomp}
 M(X)=\bigoplus_i M(Z_i)(d_i),
\end{equation}
where $M(Z_i)$ are the motives of the components of the fixed point set and
$M(Z_i)(d_i)$ are Tate twists. The class in the Grothendieck ring $K_0(\cV_k)$
decomposes then as
\begin{equation}\label{BBK0}
[X]=\sum_i [Z_i] \bL^{d_i}. 
\end{equation}
It is then immediate that, if the components $Z_i$ admit a geometric torification
(respectively, a torification of the Grothendieck class) then the variety $X$ also does.
If $Z_i=\cup_{j=1}^{n_i} T_{ij}$ with $T_{ij}= \bG_m^{a_{ij}}$ or, respectively 
$[Z_i]=\sum_{j=1}^{n_i} (\bL-1)^{a_{ij}}$,
then $X=\cup_{i=0}^n (X_i\smallsetminus X_{i-1}) = \cup_{i=0}^n   \cF^{d_i}(Z_i)$, where
$\cF^{d_i}(Z_i)$ denotes the total space of the affine fibration $\phi_i : X_i\smallsetminus X_{i-1} \to Z_i$
with fibers $\A^{d_i}$. The Grothendieck class is then torified by 
$$ [X]=\sum_{i=1}^n \sum_{j=1}^{n_i} \bT^{a_{ij}}  (1+\sum_{k=1}^{d_i} \binom{d_i}{k}  \bT^k),   $$
with $\bT=\bL-1$ the class of the multiplicative group $\bT=[\bG_m]$, 
and where the affine spaces are torified by 
$$ \bL^n-1 =\sum_{k=1}^n \binom{n}{k}  \bT^k. $$
If the restriction of the fibration $\cF^{d_i}(Z_i)$ to the individual tori $T_{ij}$ of the torification of 
$Z_i$ is trivial, then it can
be torified by a products $T_{ij} \times T_k$ of the torus $T_{ij}$ and the tori $T_k$ of a torification of
the fiber affine space $\A^{d_i}$. This determines a 
a geometric torification of the affine fibrations $\cF^{d_i}(Z_i)$, hence of $X$.
\endproof

\smallskip
\subsection{An example of torified varieties} \label{BBexSec}   

A physically significant example of torified varieties of the type
described in Lemma~\ref{torifBB} arises in the context of BPS state counting of \cite{CKK}.
Refined BPS state counting computes the multiplicities of BPS particles
with charges in a lattice ($K$-theory changes of even $D$-branes) for assigned spin
quantum numbers of a ${\rm Spin}(4)=SU(2)\times SU(2)$ representation, see
\cite{CKK}, \cite{ChoiMai}, \cite{DreMai}.

\smallskip

We mention here the following explicit example from \cite{ChoiMai}, namely
the case of the moduli space $\cM_{\P^2}(4,1)$ of Gieseker semi-stable shaved on $\P^2$
with Hilbert polynomial equal to $4m+1$. In this case, it is proved in \cite{ChoiMai}
that $\cM_{\P^2}(4,1)$ has a torus action of $\bG_M^2$ for which the fixed point locus
consists of $180$ isolated points and $6$ components isomorphic to $\P^1$. The
Grothendieck class, obtained through the Bia\/{l}ynicki-Birula decomposition \cite{ChoiMai}
is given by 
$$ [ \cM_{\P^2}(4,1) ]=1+2\bL +6 \bL^2 + 10 \bL^3 + 14 \bL^4 + 15 \bL^5 $$
$$ + 16 \bL^6 + 16 \bL^7 + 16 \bL^8 + 16 \bL^9 + 16\bL^{10} + 16 \bL^{11} $$
$$ + 15 \bL^{12} + 14 \bL^{13} + 10 \bL^{14} + 6 \bL^{15} + 2 \bL^{16} + \bL^{17}. $$
Note that, for a smooth projective variety with Grothendieck class that is a polynomial in the
Lefschetz motive $\bL$, the Poincar\'e polynomial and the Grothendieck class are related
by replacing $x^2$ with $\bL$, since the variety is Hodge--Tate. In torified form the above gives
$$ [ \cM_{\P^2}(4,1) ]=
\bT^{17}+19\, \bT^{16} +174\, \bT^{15} +1020\, \bT^{14} +4284\, \bT^{13} +13665\, \bT^{12} 
+34230\, \bT^{11} $$ $$ +68678\, \bT^{10} +111606\, \bT^9 +147653\, \bT^8+ 159082\, \bT^7 
+139008\, \bT^6 $$ $$ + 97643\, \bT^5 +54320\,  \bT^4 +23370\, \bT^3 +7468\, \bT^2+ 1632\, \bT+192, $$
where $192= \chi(\cM_{\P^2}(4,1))$ is the Euler characteristics, which is also the number of points
over $\F_1$. The number of points over $\F_{1^m}$ gives $864045$ for $m=1$ (the number of tori in the torification), $383699680$ for $m=2$ (roots of unity of order two),
$36177267945$ for $m=3$ (roots of unity of order three), etc. 

\smallskip

In this example, the Euler characteristic $\chi(\cM_{\P^2}(4,1))$, which can also
be seen as the number of $\F_1$-points, is interpreted physically as determining the
BPS counting. It is natural to ask whether the counting of $\F_{1^m}$-points, which
corresponds to the counting of roots of unity in the tori of the torification, can also
carry physically significant information. 

\smallskip

Other examples of torified varieties relevant to physics can be found in the
context of quantum field theory, see \cite{BejMar} and \cite{Mu}.

\smallskip
\subsubsection{BPS counting and the virtual motive}\label{BPSsec}  

The formulation
of the refined BPS counting given in \cite{CKK} can be summarized as follows.
The virtual motive $[X]_{\rm vir}=\bL^{-n/2} [X]$, with $n=\dim(X)$, is a class in the ring of
 motivic weights $K_0(\cV)[\bL^{-1/2}]$, see \cite{BeBrSz}. When $X$ admits a $\bG_m$
 action and a Bia\/{l}ynicki-Birula decomposition as discussed in the previous section, 
 where all the components $Z_i$ of the fixed point locus of the $\bG_m$-action have
 Tate classes $[Z_i]=\sum_j c_{ij} \bL^{b_{ij}} \in K_0(\cV)$, with $c_{ij}\in \Z$ and $b_{ij}\in \Z_+$,
 the virtual motive $[X]_{\rm vir}$ is a Laurent polynomial in the square root $\bL^{1/2}$ of the
 Lefschetz motive,
 \begin{equation}\label{Xvirtclass}
  [X]_{\rm vir} = \sum_{i,j} c_{ij}\, \bL^{b_{ij} + d_i-1/2} , 
 \end{equation} 
 where, as before, $d_i$ is the dimension of the positive eigenspace of the $\bG_m$-action on
the tangent space of $X$ at points of $Z_i$. In applications to BPS counting, one
considers the virtual motive of a moduli space $M$ that admits a perfect obstruction theory,
so that it has virtual dimension zero and an associated invariant $\#_{\rm vir} M$ which
is computed by a virtual index
$$  \#_{\rm vir} M = \chi_{\rm vir}(M, K_{M,{\rm vir}}^{1/2}) 
=\chi(M,K_{M,{\rm vir}}^{1/2}\otimes \cO_{M,{\rm vir}}),$$
where $\cO_{M,{\rm vir}}$ is the virtual structure sheaf and 
$K_{M,{\rm vir}}^{1/2}$ is a square root of the
virtual canonical bundle, see \cite{FaGo}.

\smallskip
\subsubsection{The formal square root of the Leftschetz motive} \label{SqrtSec}  

The formal square root $\bL^{1/2}$ of the Leftschetz motive that occurs in \eqref{Xvirtclass} as
Grothendieck class can be introduced, at the
level of the category of motives, as shown in \S 3.4 of \cite{KoSo}, using the Tannakian
formalism, \cite{Deligne}. Let $\cC={\rm Num}_\Q^\dagger$ be the Tannakian category of pure motives 
with the numerical equivalence relation and the Koszul sign rule twist $\dagger$ in the 
tensor structure, with motivic Galois group $G={\rm Gal}(\cC)$. The inclusion of the
Tate motives (with motivic Galois group $\bG_m$) determines a group homomorphism
$t: G \to \bG_m$, which satisfies $t\circ w=2$ with the weight homomorphism $w:\bG_m \to G$
(see \S 5 of \cite{DeMi}).
The category $\cC(\Q(\frac{1}{2}))$ obtained by adjoining a square root of the Tate motive to $\cC$
is then obtained as the Tannakian category whose Galois group is the fibered product
$$ G^{(2)}=\{ (g,\lambda)\in G\times \bG_m\,:\, t(g)=\lambda^2 \}. $$

The construction of square roots of Tate motives described in \cite{KoSo} was generalized
in \cite{LoMa} to arbitrary $n$-th roots of Tate motives, obtained via the same Tannakian
construction, with the category $\cC(\Q(\frac{1}{n}))$ obtained by adjoining an $n$-th
root of the Tate motive determined by its Tannakian Galois group
$$ G^{(n)}=\{ (g,\lambda)\in G\times \bG_m\,:\, t(g)=\sigma_n(\lambda) \}, $$
with $\sigma_n:\bG_m \to \bG_m$, $\sigma_n(\lambda)=\lambda^n$. The category
$\hat\cC$ obtained by adjoining to $\cC={\rm Num}_\Q^\dagger$ arbitrary roots of the
Tate motives is the Tannakian category with Galois group $\hat G=\varprojlim_n G^{(n)}$.
The category $\hat\cC$ has an action of $\Q^*_+$ by automorphisms induced by the
endomorphisms $\sigma_n$ of $\bG_m$. 
These roots of Tate motives give rise to a good formalism of $\F_\zeta$-geometry,
with $\zeta$ a root of unity, lying ``below" $\F_1$-geometry and expressed at the
motivic level in terms of a Habiro ring type object associated to the Grothendieck ring
of orbit categories of $\hat\cC$, see \cite{LoMa}.

\smallskip
\subsection{Counting $\F_1$-points and zeta function}  \label{CountF1zetaSec}   

For a variety $X$ over $\Z$ that is polynomially countable (that is, the counting
functions $N_X(q)=\# X_p(\F_q)$ with $X_p$ the mod $p$ reduction is a polynomial
in $q$ with $\Z$ coefficients)
the counting of points over the ``extensions" $\F_{1^m}$ (in the sense
of \cite{KapSmi}) can be obtained as the values $N_X(m+1)$ (see Theorem~4.10 of
\cite{CoCo} and Theorem~1 of \cite{Deit}). As we discussed earlier, in the case of
a torified variety, with Grothendieck class $[X]=\sum_{i\geq 0} a_i \bT^i$ with
$a_i\in \Z_+$, this corresponds to the counting given in \eqref{F1mpts}.
This is the counting of the number of $m$-th roots of unity in each torus 
$\bT^i=[\bG_m^i]$ of the torification.

\smallskip

For a variety $X$ over a finite field $\F_q$ the Hasse--Weil zeta function is given,
in logarithmic form by
\begin{equation}\label{HasseWeil}
\log Z_{\F_q}(X,t) = \sum_{m\geq 1} \frac{\# X(\F_{q^m})}{m} t^m .
\end{equation}
In the case of torified varieties, there is an analogous zeta function
over $\F_1$. We think of this $\F_1$-zeta function as defined on torified 
Grothendieck classes, $Z_{\F_1}([X],t)$. In the case of geometric torifications, 
we can regard it as a function of the variety and the torification,
$Z_{\F_1}((X,T),t)$. For simplicity of notation, we will simply write 
$Z_{\F_1}(X,t)$ by analogy to the Hasse--Weil zeta function, with
\begin{equation}\label{defF1zeta}
\log Z_{\F_1}(X,t) := \sum_{m\geq 1} \frac{\# X(\F_{1^m})}{m} t^m .
\end{equation}

\begin{lem}\label{lemF1HasseWeil}
Let $X$ be a variety over $\Z$ with a torified Grothendieck class 
$[X]=\sum_{k\geq 0} a_k \bT^k$ with $a_k\in \Z_+$. Then the $\F_1$-zeta function
is given by
\begin{equation}\label{HasseWeilF1}
\log Z_{\F_1}(X,t) = \sum_{k=0}^N a_k \, {\rm Li}_{1-k}(t) ,
\end{equation}
where ${\rm Li}_s(t)$ is the polylogarithm function with ${\rm Li}_1(t)=-\log(1-t)$ and for $k\geq 1$
$$ {\rm Li}_{1-k}(t) =  (t \frac{d}{dt})^{k-1} \frac{t}{1-t}. $$
\end{lem}

\proof
For $[X]=\sum_{k\geq 0} a_k \bT^k$ with $a_k\in \Z_+$ as above, we can consider
a similar zeta function based on the counting of $\F_{1^m}$-points described
above. Using \eqref{F1mpts}, we obtain an expression of the form
$$
\log Z_{\F_1}(X,t) = \sum_{m\geq 1} \frac{\# X(\F_{1^m})}{m} t^m =
\sum_{k=0}^N a_k \sum_{m\geq 1} m^{k-1} t^m = \sum_{k=0}^N a_k \, {\rm Li}_{1-k}(t) ,
$$
given by a linear combination of polylogarithm functions ${\rm Li}_s(t)$ at integer values $s\leq 1$.
\endproof

Such polylogarithm functions can be expressed explicitly in the form ${\rm Li}_1(t)=-\log(1-t)$ 
and for $k\geq 1$
$$ {\rm Li}_{1-k}(t)= (t \frac{d}{dt})^{k-1} \frac{t}{1-t} =  \sum_{\ell=0}^{k-1} \ell ! \, S(k,\ell+1) \left( \frac{t}{1-t} \right)^{\ell+1}, $$
with $S(k,r)$ the Stirling numbers of the second kind
$$ S(k,r)=\frac{1}{r!} \sum_{j=0}^r (-1)^{r-j} \binom{r}{j} j^k. $$

\smallskip

As in the case of the Hasse--Weil zeta function over $\F_q$ (see \cite{Ram}), the
$\F_1$-zeta function gives an exponentiable motivic measure.

\begin{prop}\label{zetaF1motmeas}
The $\F_1$-zeta function is an exponentiable motivic measure, that is, a ring
homomorphism $Z_{\F_1}: K_0(\cT)^a \to W(\Z)$ from the Grothendieck ring of
torified varieties (with either $a=w,o,s$) to the Witt ring.
\end{prop}

\proof 
Clearly with respect to addition in the Grothendieck ring of torified varieties we have
$[X]+[X'] =\sum_{i\geq 0} a_i \bT^i + \sum_{j \geq 0} a'_j \bT^j=\sum_{k\geq 0} b_k \bT^k$
with $b_k=a_k+a'_k$, hence
$$ \log Z_{\F_1}([X]+[X'],t) = \sum_{k=0}^N b_k \, {\rm Li}_{1-k}(t) = \log Z_{\F_1}([X],t)+\log Z_{\F_1}([X'],t). $$
The behavior with respect to products $[X]\cdot [Y]$ in the Grothendieck ring of torified varieties
can be analyzed as in \cite{Ram} for the Hasse--Weil zeta function. We view the $\F_1$-zeta function
$$ Z_{\F_1}(X,t) =\exp( \sum_{k=0}^N a_k \, {\rm Li}_{1-k}(t) ) $$
as the element in the Witt ring $W(\Z)$ with ghost components $\# X(\F_{1^m})=\sum_{k=0}^N m^k$,
by writing the ghost map ${\rm gh}: W(\Z)\to \Z^\N$ as
$$ {\rm gh}: Z(t)=\exp(\sum_{m\geq 1} \frac{N_m}{m}\, t^m) \mapsto t\frac{d}{dt}\, \log Z(t) =\sum_{m\geq 1} N_m t^m \mapsto (N_m)_{m\geq 1}. $$
The ghost map is an injective ring homomorphism. Thus, it suffices to see that
on the ghost components $N_m(X\times Y)=N_m(X)\cdot N_m(Y)$. If $[X]=\sum_{k\geq 0} a_k \bT^k$
and $[Y]=\sum_{\ell\geq 0} b_\ell \bT^\ell$ then $[X\times Y]=\sum_{n\geq 0}
\sum_{k+\ell=n} a_k b_\ell \bT^n$ and $N_m(X\times Y)=\sum_{n\geq 0}
\sum_{k+\ell=n} a_k b_\ell m^n=N_m(X)\cdot N_m(Y)$.
\endproof

\smallskip
\subsection{Relation to the Hasse--Weil zeta function}\label{F1HWSec}  

We discuss here the relation between the $\F_1$-zeta function $Z_{\F_1}(X,t)$ introduced
in \eqref{defF1zeta}
above, for a variety $X$ over $\Z$ with torified Grothendieck class $[X]=\sum_{k\geq 0} a_k \bT^k$,
and the Hesse--Weil zeta function $Z_{\F_q}(X,t)$, defined as in \eqref{HasseWeil}.

\smallskip

\begin{defn}\label{Z0Z1WZ}{\rm 
Consider the following elements in the Witt ring $W(\Z)$, for $k\geq 0$:
\begin{equation}\label{Z0Witt}
Z_{0,k,q}(t):= \exp( \sum_{m\geq 1} (q-1)^k\, \frac{t^m}{m} ) =\frac{1}{(1-t)^{(q-1)^k}}
\end{equation}
\begin{equation}\label{Z1Witt}
Z_{1,k,q}(t):=\exp( \sum_{m\geq 1} (\#\P^{m-1}(\F_q))^k\,\, \frac{t^m}{m} )
= \exp( \sum_{m\geq 1} (\sum_{\ell=0}^{m-1} q^\ell)^k \frac{t^m}{m} ) .
\end{equation}
}\end{defn}

\smallskip

\begin{lem}\label{ZF1HWlem}
Let $Z_{\F_q}(\bT^k,t)$ be the Hasse-Weil zeta function of a torus $\bT^k$.
The function $Z_{0,k,q}(t)$ of \eqref{Z0Witt} divides $Z_{\F_q}(\bT^k,t)$
in the Witt ring with quotient the function $Z_{1,k,q}(t)$ of \eqref{Z1Witt}.
\end{lem}

\proof
Given elements $Q=Q(t)$ and $P=P(t)$ in the Witt ring $W(\Z)$, we have that $Q$ divides $P$
iff the ghost components $q_m$ of $Q$ in $\Z^\N$ divide the corresponding ghost components $p_m$
of $P$. There is then an element $S=S(t)$ in $W(\Z)$, with ghost components $s_m=p_m/q_m$,
such that the Witt product gives $S \star_W Q =P$. The $m$-th ghost components of $Z_{\F_q}(\bT^k,t)$
is $(q^m-1)^k =\# \bT^k(\F_{q^m})$, and we have $(q^m-1)^k/(q-1)^k =(1+q+\cdots + q^{m-1})^k$.
\endproof

\smallskip

Given elements $Q,P\in W(\Z)$ such that $Q|P$ as above, we write $S=P/_W Q$ for the resulting
element $S\in W(\Z)$ with $S\star_W Q=P$.

\smallskip

The $\F_1$-zeta function of \eqref{defF1zeta} is obtained from the Hasse--Weil zeta function of 
\eqref{HasseWeil} in the following way.

\begin{prop}\label{ZF1HWprop}
Let $X$ be a variety $X$ over $\Z$ with torified Grothendieck class $[X]=\sum_{k\geq 0} a_k \bT^k$.
The $\F_1$-zeta function is given by
\begin{equation}\label{ZF1HWlim}
Z_{\F_1}(X,t)=\lim_{q\to 1}\,\, {}^W\sum_{k\geq 0} \left(  Z_{\F_q}(\bT^k,t)\,/_W\,\, Z_{0,k,q}(t) \right)^{a_k} =
\lim_{q\to 1} \,\, {}^W\sum_{k\geq 0} Z_{1,k,q}(t)^{a_k} ,
\end{equation}
while the Hasse--Weil zeta function is given by 
\begin{equation}\label{ZHWeq}
 Z_{\F_q}(X,t)=\,\, {}^W\sum_{k\geq 0} Z_{\F_q}(\bT^k,t)^{a_k}, 
\end{equation} 
where ${}^W\sum$ denotes the sum in the Witt ring.
\end{prop}

\proof For the Hasse--Weil zeta function we have 
$$ Z_{\F_q}(X,t)=\exp( \sum_{m\geq 1} \#X(\F_{q^m}) \frac{t^m}{m}) =
\exp(\sum_{k\geq 0} a_k  \sum_{m\geq 1} (q^m-1)^k \frac{t^m}{m}) $$
$$ =\prod_{k\geq 0} \exp(a_k \log Z_{\F_q}(\bT^k,t)), $$
hence we get \eqref{ZHWeq}. To obtain the $\F_1$-zeta function we then 
use Lemma~\ref{ZF1HWlem} and the fact that $(q^m-1)^k/(q-1)^k =(1+q+\cdots + q^{m-1})^k$, 
with $\lim_{q\to 1} (1+q+\cdots + q^{m-1})^k=m^k$.
\endproof

\smallskip
\subsection{Dynamical zeta functions} \label{DynZetaSec}  

The dynamical approach to $\F_1$-structures proposed in \cite{ManMar2} is
based on the existence of an endomorphism $f: X \to X$ that induces a quasi-unipotent
morphism $f_*$ on the homology $H_*(X,\Z)$. In particular, this means that the
map $f_*$ has eigenvalues that are roots of unity. 

\smallskip

In the case of a variety $X$ endowed with a torification $X=\sqcup_i T^{d_i}$,
one can consider in particular endomorphisms $f: X \to X$ that preserve the
torification and that restrict to endomorphisms of each torus $T^{d_i}$.

\smallskip

We recall the definition and main properties of the relevant dynamical zeta functions,
which we will consider in Proposition~\ref{zetadynmotmeas}.

\smallskip
\subsubsection{Properties of dynamical zeta functions}\label{DynZetaProSec}  

In general to a self-map $f: X \to X$, one can associate the dynamical
Artin--Mazur zeta function and the homological Lefschetz zeta function.
A particular class of maps with the property that they induce quasi-unipotent
morphisms in homology is given by the Morse--Smale diffeomorphisms of
smooth manifolds, see \cite{ShuSul}. These are diffeomorphisms characterized
by the properties that the set of nonwandering points is finite and hyperbolic, 
consisting of a finite number of periodic points, and for any pair of these points $x,y$ the
stable and unstable manifolds $W^s(x)$ and $W^u(y)$ intersect transversely.
Morse--Smale diffeomorphisms are structurally stable among all diffeomorphisms,
\cite{Franks}, \cite{ShuSul}. 

\smallskip

The Lefschetz zeta function is given by
\begin{equation}\label{zetaLef}
\zeta_{\cL,f}(t) =\exp\left( \sum_{m\geq 1} \frac{L(f^m)}{m} t^m \right),
\end{equation}
with $L(f^m)$ the Lefschetz number of the $m$-th iterate $f^m$, 
$$ L(f^m)=\sum_{k\geq 0} (-1)^k \Tr((f^m)_*\,|\, H_k(X,\Q)). $$
For a function with finite sets of 
fixed points ${\rm Fix}(f^m)$ this is also equal to $$L(f^m)=\sum_{x\in {\rm Fix}(f^m)} \cI(f^m,x),$$
with $\cI(f^m,x)$ the index of the fixed point. The zeta function can then be written as a rational function of the form
$$ \zeta_{\cL,f}(t)= \prod_k \det(1-t \, f_* | H_k(X,\Q))^{(-1)^{k+1}}. $$
In the case of a map $f$ with finitely many 
periodic points, all hyperbolic, the Lefschetz zeta function can be equivalently written (see \cite{Franks})
as the rational function 
$$ \zeta_{\cL,f}(t) = \prod_\gamma (1-\Delta_\gamma t^{p(\gamma)})^{(-1)^{u(\gamma)+1} }. $$
Here the product is over periodic orbits $\gamma$ with least period $p(\gamma)$ and  
$u(\gamma)=\dim E^u_x$, for $x\in \gamma$, is the dimension of the span of eigenvectors 
of $D f^{p(\gamma)}_x: T_x M \to T_x M$ with eigenvalues $\lambda$ with $|\lambda|>1$. 
One has $\Delta_\gamma=\pm 1$
according to whether $D f^{p(\gamma)}_x$ is orientation preserving or reversing. 
The relation comes from the identity $\cI(f^m,x)=(-1)^{u(\gamma)}\Delta_\gamma$.
The Artin--Mazur zeta function is given by
\begin{equation}\label{zetaAM}
\zeta_{AM,f}(t) =\exp\left( \sum_{m\geq 1} \frac{\# {\rm Fix}(f^m)}{m} t^m \right). 
\end{equation}
The case of Morse--Smale diffeomorphisms can be treated as in \cite{Franks2}
to obtain rationality and a description in terms of the homological zeta functions.

\smallskip

In the setting of real tori $\R^d/\Z^d$, one can considers the case of a toral 
endomorphism specified by a matrix $M \in M_d(\Z)$.
In the hyperbolic case, the counting of isolated fixed points of $M^m$ is given
by $|\det(1-M^m)|$ and the dynamical Artin--Mazur 
zeta function is expressible in terms of the Lefschetz zeta function,
associated to the signed counting of fixed points, through the fact that
the Lefschetz zeta function agrees with the zeta function 
\begin{equation}\label{AMzetasign}
\zeta_M(t) = \exp( \sum_{n\geq 1} \frac{t^n}{n}  a_n) , \ \ \ \text {with } \ \ \ a_n=\det(1-M^n),
\end{equation}
where $a_n=\det(1-M^n)$ is a signed fixed point counting. The general
relation between the zeta functions for the signed $\det(1-M^n)$ and for
$|\det(1-M^m)|$ is shown in \cite{Baake} for arbitrary toral endomorphisms,
with $M\in M_d(\Z)$.

\smallskip

In the case of complex algebraic tori $T^d=\bG_m^d(\C)$, one can
similarly consider the endomorphisms action of the semigroup of matrices $M\in M_d(\Z)^+$
by the linear action on $\C^d$ preserving $\Z^d$ and the exponential map
$0 \to \Z \to \C \to \C^* \to 1$ so that, for $M=(m_{ab})$ and $\lambda_a =\exp(2\pi i u_a)$,
with the action given by
$$  \lambda=(\lambda_a) \mapsto M(\lambda)= \exp(2\pi i \sum_b m_{ab} u_b). $$
The subgroup $\SL_n(\Z)\subset M_n(\Z)^+$ acts by automorphisms.
These generalize the Bost--Connes endomorphisms $\sigma_n: \bG_m \to \bG_m$,
which correspond to the ring homomorphisms of $\Z[t,t^{-1}]$ given by $\sigma_n: P(t)\mapsto P(t^n)$ 
and determine multivariable versions of the Bost--Connes algebra, see \cite{Mar}. We can
consider in this way maps of complex algebraic tori $T_\C^d=\bG_m^d(\C)$ that induce maps
of the real tori obtained as the subgroup $T_\R^d=U(1)^d\subset \bG_m^d(\C)$, and associate
to these maps the Lefschetz and Artin--Mazur zeta functions of the induced map of real tori. 

\smallskip
\subsubsection{Torifications and dynamical zeta functions} \label{TorDynZetaSec}  

In the case of a variety with a torification, we consider endomorphisms $f: X \to X$
that preserves the tori of the torification and restricts to each torus $T^{d_i}$ 
to a diffeomorphism $f_i: T_\R^{d_i} \to T_\R^{d_i}$. In particular, we consider
toral endomorphism with a matrix $M_i\in M_{d_i}(\Z)$,
we can associate to the pair $(X,f)$ a zeta function of the form
\begin{equation}\label{zetafi}
 \zeta_{\cL,f}(X,t) = \prod_i \zeta_{\cL,f_i}(t) , \ \ \ \ \   \zeta_{AM,f}(X,t) = \prod_i \zeta_{AM,f_i}(t). 
\end{equation}

\medskip

\begin{prop}\label{zetadynmotmeas}
The zeta functions \eqref{zetaLef} and \eqref{zetaAM} define exponentiable motivic measures 
on the Grothendieck
ring $K_0^\Z(\cV_\C)$ of \S 6 of \cite{ManMar2} with values in the Witt ring $W(\Z)$. The zeta functions
\eqref{zetafi} define exponentiable motivic measures on the Grothendieck ring $K_0(\cT)^a$ of torified varieties
with values in $W(\Z)$.
\end{prop}

\proof The Grothendieck ring $K_0^\Z(\cV_\C)$ considered in \S 6 of \cite{ManMar2} consists
of pairs $(X,f)$ of a complex quasi-projective variety and an automorphism $f: X \to X$ that induces
a quasi-uniponent map $f_*$ in homology. The addition is simply given by the disjoint union, and
both the counting of periodic points $\#{\rm Fix}(f^m)$ and the Lefschetz numbers
$L(f^m)$ behave additively under disjoint unions. Thus, the zeta functions $\zeta_{\cL,f}(t)$
and $\zeta_{AM,f}(t)$, seen as elements in the Witt ring $W(\Z)$ add
$$ \zeta_{\cL,f_1 \sqcup f_2}(t) =\exp\left( \sum_{m\geq 1} \frac{L((f_1 \sqcup f_2)^m)}{m} t^m \right) = $$ $$ 
\exp\left( \sum_{m\geq 1} \frac{L(f_1^m)}{m} t^m \right) \cdot \exp\left( \sum_{m\geq 1} \frac{L(f_2^m)}{m} t^m \right) = \zeta_{\cL,f_1}(t)  +_W \zeta_{\cL,f_2}(t) $$
and similarly for $\zeta_{AM, f_1 \sqcup f_2}(t) =\zeta_{AM, f_1}(t) +_W \zeta_{AM,f_2}(t)$. The product
is given by the Cartesian product $(X_1,f_1)\times (X_2,f_2)$. Since ${\rm Fix}((f_1\times f_2)^m)
={\rm Fix}(f_1^m)\times {\rm Fix}(f_2^m)$ and the same holds for Lefschetz numbers since
$$ L((f_1\times f_2)^m)=\sum_{k\geq 0} \sum_{\ell+r=k} (-1)^{\ell+r} \Tr((f_1^m)_* \otimes (f_2^m)_*\,|\, H_\ell(X_1,\Q)\otimes H_r(X_2,\Q)) $$ which gives $L(f_1^m)\cdot L(f_2^m)$.
Thus, we can use as in Proposition~\ref{zetaF1motmeas} the fact that the ghost map
${\rm gh}: W(\Z)\to \Z^\N$ 
$$ {\rm gh}: \exp(\sum_{m\geq 1} \frac{N_m}{m}\, t^m) \mapsto \sum_{m\geq 1} N_m t^m 
\mapsto (N_m)_{m\geq 1} $$
is an injective ring homomorphism to obtain the multiplicative property. The case of the torified
varieties and the zeta functions \eqref{zetafi} is analogous, combining the additive and
multiplicative behavior of the fixed point counting and the Lefschetz numbers on each torus
and of the decomposition into tori as in Proposition~\ref{zetaF1motmeas}.
\endproof

In the case of quasi-unipotent maps of tori the Lefschetz zeta function can
be computed completely explicitly. Indeed, 
it is shown in \cite{Berr1}, \cite{Berr2} that, for a quasi-unipotent self map
$f: T_\R^n \to T_\R^n$, the Lefschetz zeta function has an explicit form
that is completely determined by the map on the first homology. Under the
quasi-unipotent assumption all the eigenvalues of the induced map on $H_1$
are roots of unity, hence the characteristic polynomial $\det(1-t \, f_* | H_1(X))$
is a product of cyclotomic polynomials $\Phi_{m_1}(t) \cdots \Phi_{m_N}(t)$ where
$$ \Phi_m(t) =\prod_{d|m} (1-t^d)^{\mu(m/d)}, $$
with $\mu(n)$ the M\"obius function. It is shown in \cite{Berr2} that the Lefschetz
zeta function has the form
\begin{equation}\label{LefTnquasiuni}
\zeta_{\cL,f}(t)=\prod_{d|m} (1-t^d)^{-s_d},
\end{equation}
where $m={\rm lcm}\{m_1,\ldots, m_N\}$ and 
$$ s_d=\frac{1}{d} \sum_{k|d} F_k\, \mu(d/k) $$
$$ F_k = \prod_{i=1}^N (\Phi_{m_i/(k,m_i)}(1))^{\varphi(m_i)/\varphi(m_i/(k,m_i))} $$
where the Euler function 
$$ \varphi(m)= m \prod_{p|m, \, p \text{ prime}} (1-p^{-1}) $$
is the degree of $\Phi_m(t)$.

\smallskip

\begin{rem}\label{poscharzeta} {\rm
The properties of dynamical Artin--Mazur zeta functions change significantly when,
instead of considering varieties over $\C$ one considers varieties in positive
characteristic, \cite{Bridy}, \cite{BysCor}. The prototype model of this phenomenon 
is illustrated by considering the Bost--Connes endomorphisms 
$\sigma_n : \lambda \mapsto \lambda^n$ of $\bG_m(\bar\F_p)$. In this case, the
dynamical zeta function of $\sigma_n$ is rational or transcendental depending on whether $p$
divides $n$ (Theorem~1.2 and 1.3 and \S 3 of \cite{Bridy} and Theorem~1 of \cite{Bridy2}).
Similar phenomena in the more general case of endomorphisms of Abelian varieties in
positive characteristic have been investigated in \cite{BysCor}. In the positive
characteristic setting, where one is considering the characteristic $p$ version of the Bost--Connes 
system of \cite{CCM}, one should then replace the dynamical zeta function by the
tame zeta function considered in \cite{BysCor}.  }
\end{rem}

\medskip
\section{Spectra and zeta functions}\label{ZetaSec}   

We have already discussed in \S \ref{CountF1zetaSec} and \S \ref{DynZetaSec} 
zeta functions arising from certain counting functions that define 
ring homomorphisms from suitable Grothendieck rings to the Witt ring $W(\Z)$. 
We consider here a more general setting of exponentiable motivic measures.

\smallskip

A motivic measure is a ring homomorphism $\mu: K_0(\cV) \to R$, from
the Grothendieck ring of varieties $K_0(\cV)$ to a commutative ring $R$.
Examples include the counting measure, for varieties defined over finite fields,
which counts the number of algebraic points over $\F_q$, the topological
Euler characteristic or the Hodge--Deligne polynomials for complex algebraic 
varieties. 

\smallskip

The Kapranov motivic zeta function \cite{Kapr} is defined as $\zeta(X,t)=\sum_{n=0}^\infty [S^n(X)] t^n$,
where $S^n(X)=X^n/S_n$ are the symmetric products of $X$ and $[S^n(X)]$ are the classes in
$K_0(\cV)$. Similarly, the zeta function of a motivic measure is defined as
\begin{equation}\label{zetamu}
 \zeta_\mu (X, t) =\sum_{n=0}^\infty \mu(S^n(X))\, t^n .
\end{equation}
It is viewed as an element in the Witt ring $W(R)$. The addition in $K_0(\cV)$
is mapped by the zeta function to the addition in $W(R)$, which is the usual 
product of the power series,
\begin{equation}\label{zetasum}
 \zeta_\mu (X \sqcup Y, t)= \zeta_\mu (X, t) \cdot \zeta_\mu (Y, t) =\zeta_\mu (X, t) +_{W(R)} \zeta_\mu (Y, t). 
\end{equation} 
The motivic measure $\mu: K_0(\cV) \to R$ is said to be exponentiable (see \cite{Ram}, \cite{RamTab})
if the zeta function \eqref{zetamu} defines a ring homomorphism
$$ \zeta_\mu : K_0(\cV) \to W(R), $$
that is, if in addition to \eqref{zetasum} one also has 
\begin{equation}\label{zetaprod}
 \zeta_\mu (X \times Y, t)= \zeta_\mu (X, t) \star_{W(R)} \zeta_\mu (Y, t) . 
\end{equation} 

\smallskip

We investigate here how to lift the zeta functions of exponentiable motivic measures
to the level of spectra. To this purpose, we first investigate how to construct a spectrum
whose $\pi_0$ is a dense subring $W_0(R)$ of the Witt ring $W(R)$ and then we
consider how to lift the ring homomorphisms given by zeta functions $\zeta_\mu$
of exponentiable measures with a rationality and a factorization condition. 

\smallskip
\subsection{The Endomorphism Category}\label{EndoCatSec}  

Let $R$ be a commutative ring. We denote by $\cE_R$ the
endomorphism category of $R$, which is defined as follows
(see \cite{Alm1}, \cite{Alm2}, \cite{DreSie}).

\begin{defn}\label{ERcat}  {\rm 
The category $\cE_R$ has objects given by the pairs
$(E,f)$ of a finite projective module $E$ over $R$ and
an endomorphism $f\in {\rm End}_R(E)$, and morphisms
given by morphisms $\phi: E \to E'$ of finite projective 
modules that commute with the endomorphisms,
$f' \circ \phi = \phi \circ f$. 
The endomorphism category has direct sum 
$(E,f)\oplus (E',f')=(E\oplus E', f\oplus f')$ and tensor product
$(E,f)\otimes (E',f') =(E\otimes E', f\otimes f')$.  }
\end{defn}

\smallskip

The category of finite
projective modules over $R$ is identified with the subcategory
corresponding to the objects $(E,0)$ with trivial endomorphism.

\smallskip

An exact sequence in $\cE_R$ is a sequence of
objects and morphisms in $\cE_R$ which is 
exact as a sequence of finite projective modules over $R$
(forgetting the endomorphisms). This determines a collection
of admissible short exact sequence (and of admissible
monomorphisms and epimorphisms). 
The endomorphism category
$\cE_R$ is then an exact category, hence it has an associated
$K$-theory defined via the Quillen $Q$-construction, \cite{Quillen}. 
This assigns to the exact category $\cE_R$ the category $\cQ\cE_R$
with the same objects and morphisms $\Hom_{\cQ\cE_R}((E,f),(E',f'))$
given by diagrams 
$$ \xymatrix{ & (E'',f'') \ar[dl] \ar[dr] & \\
(E,f) & & (E',f'),
} $$
where the first arrow is an admissible epimorphism and the
second an admissible monomorphism, with composition
given by pullback. By the Quillen construction $K$-theory 
of $\cE_R$ is then $K_{n-1}(\cE_R)=\pi_n(\cN(\cQ\cE_R))$, 
with $\cN(\cQ\cE_R)$ the nerve of $\cQ\cE_R$. 

\smallskip

The forgetful functor
$(E,f) \mapsto E$ induces a map on $K$-theory
$$ K_n(\cE_R) \to K_n(\cP_R)=K_n(R), $$ 
which is a split surjection. Let 
$$ \cE_n(R):={\rm Ker}(K_n(\cE_R) \to K_n(R)). $$

\smallskip

In the case of $K_0$ an explicit description is given by the
following, \cite{Alm1}, \cite{Alm2}.
Let $K_0(\cE_R)$ denote the $K_0$ of the endomorphism
category $\cE_R$. It is a ring with the product structure induced
by the tensor product. It is proved in \cite{Alm1}, \cite{Alm2}
that the quotient 
\begin{equation}\label{W0R}
W_0(R)=K_0(\cE_R)/K_0(R)
\end{equation}
embeds as a dense subring of the big Witt ring $W(R)$
via the map 
\begin{equation}\label{Lmap}
L: (E,f) \mapsto \det(1-t \, M(f))^{-1}, 
\end{equation}
with $M(f)$ the matrix associated to $f\in {\rm End}_R(E)$, where 
$\det(1-t \, M(f))^{-1}$ is viewed as an element in $\Lambda(R)=1+ t R[[t]]$.
As a subring $W_0(R) \hookrightarrow W(R)$ of the big Witt ring, $W_0(R)$
consists of the rational Witt vectors
$$ W_0(R)= \left\{ \frac{1+a_1 t + \cdots + a_n t^n}{1+b_1 t+ \cdots + b_m t^m} \, | \,\, a_i,b_i\in R,\,\, n,m\geq 0 \right\}. $$
Equivalently, one can consider the ring $\cR = (1+t R[t])^{-1} R[t]$ and identify the
above with $1+t \cR$, where the multiplication in $1+t \cR$ corresponds to 
the addition in the Witt ring, and the Witt product is determined by the identity $(1-at)\star (1-b t) = (1-ab t)$.

\smallskip

This description of Witt rings in terms of endomorphism categories was
applied to investigate the arithmetic structures of the Bost--Connes
quantum statistical mechanical system, see \cite{CoCo}, \cite{MaRe}, \cite{MaTa}.

\smallskip

This relation between the Grothendieck ring and Witt vectors was 
extended to the higher $K$-theory in \cite{Gray}, where an explicit
description for the kernels $\cE_n(R)$ is obtained, by showing that
$$ \cE_{n-1}(R) = {\rm Coker}(K_n(R) \to K_n(\cR)), $$
where $\cR= (1+t R[t])^{-1} R[t]$ and $K_n(R) \to K_n(\cR)$ is a split
injection. The identification above is obtained in \cite{Gray} by
showing that there is an exact sequence
\begin{equation}\label{KnE}
 0 \to K_n(R) \to K_n(\cR) \to K_{n-1}(\cE_R) \to K_{n-1}(R) \to 0. 
\end{equation} 
The identification \eqref{W0R} for $K_0$ is then recovered as
the case with $n=0$ that gives an identification $\cE_0(R)\simeq 1+t\cR$.

\smallskip
\subsection{Spectrum of the Endomorphism Category and Witt vectors} \label{WittSpSec} 

Let $\cP_R$ denote the category
of finite projective modules over a commutative ring $R$ 
with unit. Also let $\cE_R$ be the endomorphism
category recalled in \S \ref{GammaSpSec}. By the Segal construction
described in \S \ref{GammaSpSec}, we obtain associated $\Gamma$-spaces $F_{\cP_R}$ and $F_{\cE_R}$
and spectra $F_{\cP_R}(\bS)=K(R)$, the $K$-theory spectrum of $R$, and
$F_{\cE_R}(\bS)$, the spectrum of the endomorphism category. 

\smallskip

We obtain in the following way a functorial ``spectrification" of the Witt ring $W_0(R)$, namely
a spectrum $\bW(R)$ with $\pi_0 \bW(R) =W_0(R)$.

\begin{defn}\label{defWittSpec}{\rm
For a commutative ring $R$, with $\cP_R$ the category of finite projective modules and
$\cE_R$ the category of endomorphisms, the spectrum $\bW(R)$ is defined as the cofiber
$\bW(R):=F_{\cE_R}(\bS)/F_{\cP_R}(\bS)$ obtained from the $\Gamma$-spaces 
$F_{\cP_R} : \Gamma^0 \to \Delta_*$ and $F_{\cE_R} : \Gamma^0 \to \Delta_*$
associated to the categories $\cP_R$ and $\cE_R$.
}\end{defn}

\smallskip

\begin{lem}\label{WittSpec}
For a commutative ring $R$, 
the inclusion of the category $\cP_R$ of finite projective modules as the subcategory of the
endomorphism category $\cE_R$ determines a long exact sequence
$$ \cdots \to \pi_n (F_{\cP_R}(\bS)) \to \pi_n(F_{\cE_R}(\bS)) \to \pi_n(F_{\cE_R}(\bS)/F_{\cP_R}(\bS)) \to \pi_{n-1} (F_{\cP_R}(\bS)) \to 
\cdots $$ $$ \cdots \to \pi_0 (F_{\cP_R}(\bS)) \to \pi_0 (F_{\cE_R}(\bS)) 
\to \pi_0(F_{\cE_R}(\bS)/F_{\cP_R}(\bS)) $$
of the homotopy groups of the spectra $F_{\cP_R}(\bS)$, $F_{\cE_R}(\bS)$ with cofiber
$\bW(R)$ as in Definition~\ref{defWittSpec}.
The spectrum $\bW(R)$ satisfies $\pi_0 \bW(R) =W_0(R)$. 
\end{lem}

\proof
The functoriality of the Segal construction implies that the 
inclusion of $\cP_R$ as the subcategory of $\cE_R$ 
given by objects $(E,0)$ with trivial endomorphism 
determines a map of $\Gamma$-spaces $F_{\cP_R} \to F_{\cE_R}$, which is
a natural transformation of the functors $F_{\cP_R} : \Gamma^0 \to \Delta_*$
and $F_{\cE_R} : \Gamma^0 \to \Delta_*$.
After passing to endofunctors
$F_{\cP_R} : \Delta_* \to \Delta_*$
and $F_{\cE_R} : \Delta_* \to \Delta_*$ we obtain a map of spectra $K(R) \to F_{\cE_R}(\bS)$, 
induced by the inclusion of $\cP_R$ as subcategory of $\cE_R$.
The category $\Delta_*$ of simplicial sets has products and equalizers, hence pullbacks. Thus,
given two functors $F,F':\Gamma^0\to \Delta_*$, a natural transformation $\alpha: F \to F'$
is mono if and only if for all objects $X\in \Gamma^0$ the morphism $\alpha_X : F(X)\to F'(X)$
is a monomorphism in $\Delta_*$. 
An embedding $\cC \hookrightarrow \cC'$  determines by composition an embedding
$\Sigma_{\cC}(X) \hookrightarrow \Sigma_{\cC'}(X)$ of the categories of summing functors,
for each object $X\in \Gamma^0$. This gives a monomorphism $F_{\cC}(X) = \cN \Sigma_{\cC}(X)
\to F_{\cC'}(X)=\cN \Sigma_{\cC'}(X)$, hence a monomorphism $F_{\cC} \to F_{\cC'}$ of
$\Gamma$-spaces. 
Arguing as in Lemma~1.3 of \cite{Schwede} we then obtain from such a map $F_{\cC} \to F_{\cC'}$
of $\Gamma$-spaces a long exact sequence of homotopy groups of the associated spectra
$$ \cdots \to \pi_n (F_{\cC}(\bS)) \to \pi_n( F_{\cC'}(\bS)) \to \pi_n(F_{\cC'}(\bS)/F_{\cC}(\bS) ) \to 
\pi_{n-1} (F_{\cC}(\bS)) \to \cdots $$   $$ \cdots  \to \pi_0 (F_{\cC}(\bS)) \to \pi_0( F_{\cC'}(\bS)) \to \pi_0(F_{\cC'}(\bS)/F_{\cC}(\bS) ), $$
where $F_{\cC'}(\bS)/F_{\cC}(\bS)$ is the cofiber.  When applied to the subcategory 
$\cP_R \hookrightarrow \cE_R$ this gives the long exact sequence
$$ \cdots \to \pi_n (F_{\cP_R}(\bS)) \to \pi_n(F_{\cE_R}(\bS)) \to \pi_n(F_{\cE_R}(\bS)/F_{\cP_R}(\bS)) \to \pi_{n-1} (F_{\cP_R}(\bS)) \to 
\cdots $$ $$ \cdots \to \pi_0 (F_{\cP_R}(\bS)) \to \pi_0 (F_{\cE_R}(\bS)) 
\to \pi_0(F_{\cE_R}(\bS)/F_{\cP_R}(\bS)). $$
Here we have $\pi_n (F_{\cP_R}(\bS))=K_n(R)$. Moreover, by construction we
have $\pi_0 (F_{\cE_R}(\bS))=K_0(\cE_R)$ so that we identify
$$ \pi_0(F_{\cE_R}(\bS)/F_{\cP_R}(\bS)) = W_0(R) = K_0(\cE_R)/K_0(R). $$
Thus, the spectrum $\bW(R):=F_{\cE_R}(\bS)/F_{\cP_R}(\bS)$ given by the cofiber
of $F_{\cP_R}(\bS) \to F_{\cE_R}(\bS)$ provides a spectrum whose zeroth homotopy
group is the Witt ring $W_0(R)$.
\endproof

\smallskip

The forgetful functor $\cE_R\to \cP_R$ also induces a corresponding 
map of $\Gamma$-spaces $F_{\cE_R} \to F_{\cP_R}$. Moreover, one
can also construct a spectrum with $\pi_0$ equal to $W_0(R)$ using the
characterization given in \cite{Gray}, that we recalled above, in terms of
the map on $K$-theory (and on $K$-theory spectra) $K(R)\to K(\cR)$
with $\cR=(1+rR[t])^{-1} R[t]$. One can obtain in this way a reformulation
in terms of spectra of the result of \cite{Gray}. However, for our purposes
here, it is preferable to work with the spectrum constructed in Lemma~\ref{WittSpec}. 

\smallskip

We give a variant of Lemma~\ref{WittSpec} that will be useful in the following.
We denote by $\cP_R^\pm$ and $\cE_R^\pm$, respectively, the categories of $\Z/2\Z$-graded
finite projective $R$-modules and the $\Z/2\Z$-graded endomorphism category with objects
given by pairs $\{ (E_+,f_+), (E_-,f_-) \}$, which we write simply as $(E_\pm, f_\pm)$ and
with morphisms $\phi: E_\pm \to E'_\pm$ of $\Z/2\Z$-graded finite projective modules that
commute with $f_\pm$. The sum in $\cE_R^\pm$ is given by
$$ (E_\pm, f_\pm)\oplus (E'_\pm, f'_\pm)=((E_+\oplus E'_+,E_-\oplus E'_-), (f_+\oplus f'_+,f_-\oplus f'_-)) $$
while the tensor product $(E_\pm, f_\pm)\otimes (E'_\pm, f'_\pm)$ is given by
$$ ((E_+\otimes E'_+ \oplus E_-\otimes E'_-, f_+\otimes f'_+ \oplus f_-\otimes f'_-),(E_+\otimes E'_- \oplus E_-\otimes E'_+, f_+\otimes f'_- \oplus f_-\otimes f'_+)). $$
Again we consider $\cP_R^\pm$ as a subcategory of $\cE_R^\pm$
with trivial endomorphisms. 

\begin{lem}\label{K0pmW0}
The map $\delta: K_0(\cE_R^\pm) \to K_0(\cE_R)$ given by $[E_\pm,f_\pm]\mapsto [E_+,f_+]-[E_-,f_-]$
is a ring homomorphism and it descends to a ring homomorphism
$$ K_0(\cE_R^\pm) /K_0(\cP_R^\pm) \to K_0(\cE_R)/K_0(R) \simeq W_0(R). $$
\end{lem}

\proof The map is clearly compatible with sums. Compatibility with product also holds since 
$[E_\pm,f_\pm] \cdot [E'_\pm,f'_\pm]  \mapsto 
 ([E_+,f_+]-[E_-,f_-]) \cdot ([E'_+,f'_+]-[E'_-,f'_-]).$ Moreover, it maps $K_0(\cP_R^\pm)$
 to $K_0(\cP_R)$.
\endproof

As before, the categories $\cP_R^\pm$ and $\cE_R^\pm$ have associated $\Gamma$-spaces 
$F_{\cP_R^\pm}: \Gamma^0\to \Delta_*$ and $F_{\cE_R^\pm}: \Gamma^0 \to \Delta_*$
and spectra $F_{\cP_R^\pm}(\bS)$ and $F_{\cE_R^\pm}(\bS)$. The following result follows
as in Lemma~\ref{WittSpec}. 

\begin{lem}\label{GammaWpm}
The inclusion of $\cP_R^\pm$ as a subcategory of $\cE_R^\pm$ induces a long exact sequence
$$ \cdots \to \pi_n (F_{\cP_R^\pm}(\bS)) \to \pi_n(F_{\cE_R^\pm}(\bS)) \to \pi_n(F_{\cE_R^\pm}(\bS)/F_{\cP_R^\pm}(\bS)) \to \pi_{n-1} (F_{\cP_R^\pm}(\bS)) \to 
\cdots $$ $$ \cdots \to \pi_0 (F_{\cP_R^\pm}(\bS)) \to \pi_0 (F_{\cE_R^\pm}(\bS)) 
\to \pi_0(F_{\cE_R^\pm}(\bS)/F_{\cP_R^\pm}(\bS)) $$
of the homotopy groups of the spectra $F_{\cP_R^\pm}(\bS)$ and $F_{\cE_R^\pm}(\bS)$,
which at the level of $\pi_0$ gives $K_0(\cP_R^\pm)\to K_0(\cE_R^\pm)\to K_0(\cE_R^\pm)/K_0(\cP_R^\pm)$.
\end{lem}

We denote by $\bW^\pm(R)=F_{\cE_R^\pm}(\bS)/F_{\cP_R^\pm}(\bS)$ the cofiber
of $F_{\cP_R^\pm}(\bS) \to F_{\cE_R^\pm}(\bS)$.

\smallskip

\begin{rem}\label{SpWRrem}{\rm
It is important to point out that our treatment of Witt vectors and their spectrification, as
presented in this section, differs from the one in \cite{Hess} (see especially Theorem~2.2.9 
and equation (2.2.11) in that paper), and in \cite{Camp}. Nonetheless, the circle action on THH
that is used to obtain the spectrum TR is closely related to the Bost--Connes structure
investigated in the present paper. A more direct relation between Bost--Connes structures and
topological Hochschild and cyclic homology will also relate naturally to the point of view on
$\F_1$-geometry developed in \cite{CC16}. We will leave this topic for future work.
}\end{rem}

\smallskip
\subsection{Exponentiable measures and maps of $\Gamma$-spaces}\label{zetaGammaSec} 

The problem of lifting to the level of spectra the Hasse--Weil zeta function associated 
to the counting motivic measure for varieties over finite fields was discussed in \cite{CaWoZa}. 
We consider here a very similar setting and procedure, where we want to lift a zeta function 
$\zeta_\mu: K_0(\cV) \to W(R)$ associated to an exponentiable motivic measure
to the level of spectra. To this purpose, we make some assumptions of rationality
and the existence of a factorization for our zeta functions of exponentiable motivic measures.
We then consider the spectrum $K(\cV)$ of \cite{Zak1}, \cite{Zak3} with $\pi_0K(\cV)=K_0(\cV)$
and a spectrum, obtained from a $\Gamma$-space, associated to the subring $W_0(R)$ 
of the big Witt ring $W(R)$. 

\smallskip

\begin{defn}\label{factormotmeas}{\rm 
A motivic measure, that is, a ring homomorphism $\mu: K_0(\cV)\to R$
of the Grothendieck ring of varieties to a commutative ring $R$, is called {\rm factorizable} 
is it satisfies the following three properties:
\begin{enumerate}
\item {\bf exponentiability}: the associated zeta function $\zeta_\mu(X,t)$ is a {\rm ring}
homomorphism $\zeta_\mu: K_0(\cV)\to W(R)$ to the Witt ring of $R$;
\item {\bf rationality}: the homomorphism $\zeta_\mu$ factors through the inclusion of
the subring $W_0(R)$ of the Witt ring, $\zeta_\mu: K_0(\cV)\to W_0(R) \hookrightarrow W(R)$.
\item {\bf factorization}: the rational functions $\zeta_\mu(X,t)$ admit a factorization into linear factors
$$ \zeta_\mu(X,t)= \frac{\prod_i (1-\alpha_i t)}{\prod_j (1-\beta_j t)} = \zeta_{\mu,+}(X,t)\, 
 -_W\,\, \zeta_{\mu,-}(X,t) $$
where $\zeta_{\mu,+}(X,t)=\prod_j (1-\beta_j t)^{-1}$ and $\zeta_{\mu,-}(X,t)=\prod_i (1-\alpha_i t)^{-1}$ and
$-_W$ is the difference in the Witt ring, that is the ratio of the two polynomials. 
\end{enumerate} }
\end{defn}

\begin{lem}\label{AssEnd}
A factorizable motivic measure $\mu: K_0(\cV)\to R$, as in Definition~\ref{factormotmeas}, determines a
functor $\Phi_\mu:  \cC_\cV \to \cE_R^\pm$ where $\cC_\cV$ is the assembler category encoding
the scissor-congruence relations of the Grothendieck ring $K_0(\cV)$ and $\cE_R^\pm$ is the
$\Z/2\Z$-graded endomorphism category.
\end{lem}

\proof The objects of $\cC_\cV$ are varieties $X$ and the morphisms are locally closed embeddings,
\cite{Zak1}, \cite{Zak3}. To an object $X$ we assign an object of $\cE_R$ obtained in the following
way. Consider a factorization 
$$ \zeta_\mu(X,t)= \frac{\prod_{i=1}^n (1-\alpha_i t)}{\prod_{j=1}^m (1-\beta_j t)}
 = \zeta_{\mu,+}(X,t)\,  -_W\,\, \zeta_{\mu,-}(X,t) $$
as above of the zeta function of $X$. Let $E_+^{X,\mu} =R^{\oplus m}$ and $E_-^{X,\mu}=R^{\oplus n}$ with
endomorphisms $f_\pm^{X,\mu}$ respectively given in matrix form by 
$M(f_+^{X,\mu})={\rm diag}(\beta_j)_{j=1}^m$
and $M(f_-^{X,\mu})={\rm diag}(\alpha_i)_{i=1}^n$. 
The pair $(E_\pm^{X,\mu}, f_\pm^{X,\mu})$ is an object of the endomorphism 
category $\cE_R^\pm$. Given an embedding
$Y \hookrightarrow X$, the zeta function satisfies
$$ \zeta_\mu(X,t)=\zeta_\mu(Y,t) \cdot \zeta_\mu(X\smallsetminus Y,t) = \zeta_\mu(Y,t) \, +_W \,\zeta_\mu(X\smallsetminus Y,t). $$
Using the factorizations of each term, this gives 
$$ (E_\pm^{X,\mu},f_\pm^{X,\mu}) = (E_\pm^{Y,\mu},f_\pm^{Y,\mu})\oplus (E_\pm^{X\smallsetminus Y,\mu},f_\pm^{X\smallsetminus Y,\mu}), $$
hence a morphism in $\cE_R^\pm$ given by the canonical morphism to the direct sum
$$ (E_\pm^{Y,\mu},f_\pm^{Y,\mu} ) \to  (E_\pm^{X,\mu},f_\pm^{X,\mu}) . $$
\endproof

\begin{prop}\label{PhimuGamma}
The functor $\Phi_\mu:  \cC_\cV \to \cE_R^\pm$ of Lemma~\ref{AssEnd} induces a map of
$\Gamma$-spaces and of the associated spectra $\Phi_\mu: K(\cV) \to F_{\cE_R^\pm}(\bS)$.
The induced maps on the homotopy groups has the property that the composition
\begin{equation}\label{Phidelta}
K_0(\cV) \stackrel{\Phi_\mu}{\to} K_0(\cE_R^\pm)  \stackrel{\delta}{\to} K_0(\cE_R) \to K_0(\cE_R)/K_0(R)=W_0(R)   
\end{equation}
with $\delta$ as in Lemma~\ref{K0pmW0},  is given by the zeta function $\zeta_\mu: K_0(\cV)\to W_0(R)$. 
\end{prop}

\proof The $\Gamma$-space associated to the assembler category $\cC_\cV$ is obtained in the
following way, \cite{Zak1}, \cite{Zak3}.  One first associates to the assembler category $\cC_\cV$ another
category $\cW(\cC_\cV)$ whose objects are finite collections $\{ X_i \}_{i\in I}$ of non-initial
objects of $\cC_\cV$ with morphisms $\varphi=(f,f_i): \{ X_i \}_{i\in I} \to \{ X'_j \}_{j\in J}$ given 
by a map of the indexing sets $f: I\to J$ and morphisms $f_i: X_i \to X'_{f(i)}$ in $\cC_\cV$, such
that, for every fixed $j\in J$ the collection $\{ f_i: X_i \to X'_j\, :\, i\in f^{-1}(j)\}$ is a disjoint covering
family of the assembler $\cC_\cV$. This means, in the case of the assembler $\cC_\cV$ underlying the
Grothendieck ring of varieties, that the $f_i$ are closed embeddings of the varieties $X_i$ in the given
$X'_j$ with disjoint images. We first show that the functor $\Phi_\mu:  \cC_\cV \to \cE_R^\pm$ of Lemma~\ref{AssEnd} extends to a functor (for which we still use the same notation)
$\Phi_\mu:  \cW(\cC_\cV) \to \cE_R^\pm$. We define $\Phi_\mu( \{ X_i \}_{i\in I})= \oplus_{i\in I} \Phi_\mu(X_i)=
\oplus_{i\in I} (E_\pm^{X_i,\mu}, f_\pm^{X_i,\mu})$. 
Given a covering family $\{ f_i: X_i \to X'_j\, :\, i\in f^{-1}(j)\}$ as above, each morphism
$f_i: X_i \to X'_j$ determines a morphism $\Phi_\mu(f_i) : (E_\pm^{X_i,\mu},f_\pm^{X_i,\mu}) \to
(E_\pm^{X'_j,\mu},f_\pm^{X'_j,\mu})$ given by the canonical morphism to the direct sum
$(E_\pm^{X_i,\mu},f_\pm^{X_i,\mu}) \to (E_\pm^{X_i,\mu},f_\pm^{X_i,\mu}) \oplus (E_\pm^{X'_j\smallsetminus X_i,\mu},f_\pm^{X'_j\smallsetminus X_i,\mu})$. This determines a morphism 
$\Phi_\mu(\varphi): \oplus_{i\in I}   (E_\pm^{X_i,\mu}, f_\pm^{X_i,\mu}) \to \oplus_{j\in J}  (E_\pm^{X'_j,\mu}, f_\pm^{X'_j,\mu})$. We then show that the functor $\Phi_\mu:  \cW(\cC_\cV) \to \cE_R^\pm$ constructed in
this way determines a map of the associated $\Gamma$-spaces. The $\Gamma$-space associated to
$\cW(\cC_\cV)$ is constructed in \cite{Zak1}, \cite{Zak3} as the functor that assigns to a finite pointed
set $S\in \Gamma^0$ the simplicial set given by the nerve $\cN \cW(S\wedge \cC_\cV)$, where the
coproduct of assemblers $S\wedge \cC_\cV=\bigvee_{s\in S\smallsetminus \{ s_0 \}} \cC_\cV$ 
has an initial object and a copy of the non-initial objects of $\cC_\cV$ for each point 
$s\in S\smallsetminus \{ s_0 \}$ and morphisms induced by those of $\cC_\cV$. This means that
we can regard objects of $\cW(S\wedge \cC_\cV)$ as collections $\{ X_{s,i} \}_{i\in I}$, for some 
$s\in S\smallsetminus \{ s_0 \}$ and morphisms $\varphi_s=(f_s,f_{s,i}): \{ X_{s,i} \}_{i\in I}\to 
\{ X'_{s,j} \}_{j\in J}$ as above. In order to obtain a map of $\Gamma$-spaces between 
$F_\cV: S \mapsto \cN \cW(S\wedge \cC_\cV)$ and $F_{\cE_R^\pm}: S\mapsto \cN \Sigma_{\cE_R^\pm}(S)$,
we construct a functor $\cW(S\wedge \cC_\cV)\to \Sigma_{\cE_R^\pm}(S)$ from the category 
$\cW(S\wedge \cC_\cV)$ described above to the category of summing functors $\Sigma_{\cE_R^\pm}(S)$.
To an object $X_{S,I}:=\{ X_{s,i} \}_{i\in I}$ in $\cW(S\wedge \cC_\cV)$ we associate a functor
$\Phi_{X_{S,I}}: \cP(S)\to \cE_R^\pm$ that maps a subset $A_+=\{s_0\}\sqcup A\in \cP(X)$ to
$\Phi_{X_{S,I}}(A_+)=\oplus_{a\in A} \Phi_\mu(\{ X_{a,i}\}_{i\in I})$ where $\Phi_\mu:  \cW(\cC_\cV) \to \cE_R^\pm$ is the functor constructed above. It is a summing functor since $\Phi_{X_{S,I}}(A_+ \cup B_+)=
\Phi_{X_{S,I}}(A_+)\oplus \Phi_{X_{S,I}}(B_+)$ for $A_+\cap B_+=\{ s_0 \}$. This induces a map of simplicial
sets $\cN\cW(S\wedge \cC_\cV) \to \cN\Sigma_{\cE_R^\pm}(S)$ which determines a natural transformation
of the functors $F_\cV: S \mapsto \cN \cW(S\wedge \cC_\cV)$ and $F_{\cE_R^\pm}: S\mapsto \cN \Sigma_{\cE_R^\pm}(S)$. This map of $\Gamma$-spaces in turn determines a map of the associated
spectra and an induced map of their homotopy groups. It remains to check that the induced
map at the level of $\pi_0$ agrees with the expected map of Grothendieck rings 
$K_0(\cV)\to K_0(\cE_R^\pm)$, hence with the zeta function when further mapped to $K_0(\cE_R)$
and to the quotient $K_0(\cE_R)/K_0(R)$. This is the case since by construction the induced map 
$\pi_0 K(\cV)=K_0(\cV)\to K_0(\cE_R^\pm)=\pi_0 F_{\cE_R^\pm}(\bS)$ is given by the
assignment $[X]\mapsto [ E_\pm^{X,\mu},f_\pm^{X,\mu} ]$. 
\endproof

\smallskip

\begin{cor}\label{PhideltaSp}
The map of Grothendieck rings given by the composition \eqref{Phidelta} also 
lifts to a map of spectra. 
\end{cor}

\proof
It is possible to realize the map $\delta: K_0(\cE_R^\pm)\to K_0(\cE_R)$ of Lemma~\ref{K0pmW0}
at the level of spectra. The $K$-theory spectrum of an abelian category $\cA$ is weakly equivalent to 
the $K$-theory spectrum of the category of bounded chain complexes over $\cA$.  In fact, this holds
more generally for $\cA$ an exact category closed under kernels. Thus, in the case
of the category $\cE_R$, there is a weak equivalence $K({\rm Ch}^\flat(\cE_R))\tilde\to K(\cE_R)$ which descends on the level $\pi_0$ to the map $K_0({\rm Ch}^\flat(\cE_R))\tilde\to K_0(\cE_R)$ given by $[E^\cdot, f^\cdot]\mapsto\sum_k(-1)^k[E^k, f^k]$.  To an object $(E^\pm, f^\pm)$
of $\cE_R^\pm$ we can assign a chain complex in ${\rm Ch}^\flat(\cE_R)$ of the form 
$0\to(E^-, f^-)\overset{0}{\to}(E^+, f^+)\to 0$, where $(E^+,f^+)$ sits in degree $0$. 
This descends on the level of K-theory to a map $K(\cE_R^\pm)\to K({\rm Ch}^\flat(\cE_R))$,
which at the level of $\pi_0$ gives the map $[E^\pm, f^\pm]\mapsto[E^+, f^+]-[E^-, f^-]$. The
functor $\cE_R^\pm\to {\rm Ch}^\flat(\cE_R)$ used here does not respect tensor products, although
the induced map $\delta: K_0(\cE_R^\pm)\to K_0(\cE_R)$ at the level of $K_0$ is compatible with
products.  Thus, the composition \eqref{Phidelta} can also be lifted at the level of spectra.
\endproof

\smallskip

It should be noted that the construction of a derived motivic zeta function outlined above is not the first to appear in the literature.  In \cite{CaWoZa}, the authors describe a derived motivic measure 
$\zeta:K(\cV_k)\to K(\text{Rep}_{cts}(\text{Gal}(k^s/k);\bZ_\ell))$ from the Grothendieck spectrum of varieties to the $K$-theory spectrum of the category of continuous $\ell$-adic Galois representations.  
This map corresponds to the assignment $X\mapsto H^*_{\text{et}, c}(X\times_k k^s, \bZ_\ell)$. In particular, they show that when $k=\bF_q$ for $\ell$ coprime to $q$, on the level of $\pi_0$, $\zeta$ corresponds to the Hasse-Weil zeta function.  They then use $\zeta$ to prove that $K_1(\cV_{\bF_q})$ is not only nontrivial, but contains interesting algebro-geometric data.

\smallskip

Essentially, the approach in \cite{CaWoZa} was to start with a Weil Cohomology theory (in this case, $\ell$-adic cohomology) and then to construct a derived motivic measure realizing on the level of $K$-theory the assignment to a variety $X$ of its corresponding cohomology groups.  The methods used in the case of $\ell$-adic cohomology may not immediately generalize to other Weil cohomology theories.  This method has 
yielded deep insight into the world of algebraic geometry.
Our approach here, in contrast, is to take an interesting class of motivic measures, namely Kapranov motivic zeta functions (exponentiable motivic measures, \cite{Kapr}, \cite{Ram}, \cite{RamTab}), and to determine reasonable conditions under which such a motivic measure can be derived directly.  This method still needs to be studied further to yield additional insights into what it captures about the geometry of varieties.

\smallskip
\subsection{Bost--Connes type systems via motivic measures}  \label{BCMotMeasSec}  

The lifting of the integral Bost--Connes algebra to various Grothendieck rings, their
assembler categories, and the associated spectra, that we discussed in \cite{ManMar2}
and in the earlier sections of this paper, can be viewed as an instance of a more
general kind of operation. As discussed in \cite{CoCo2}, there is a close relation between
the endomorphisms $\sigma_n$ and the maps $\tilde\rho_n$ of the integral Bost--Connes
algebra and the operation of Frobenius and Verschiebung in the Witt ring. Thus, we can
formulate a more general form of the question investigated above, of lifting of the integral 
Bost--Connes algebra to a Grothendieck ring through an Euler characteristic map, in terms
of lifting the Frobenius and Verschiebung operations of a Witt ring to a Grothendieck ring
through the zeta function $\zeta_\mu$ of an exponentiable motivic measure.  A
prototype example of this more general setting is provided by the Hasse--Weil zeta
function $Z: K_0(\cV_{\F_q})\to W(\Z)$, which has the properties that the action of the
Frobenius $F_n$ on the Witt ring $W(\Z)$ corresponds to passing to a field extension,
$F_n Z(X_{\F_q},t) = Z(X_{\F_{q^n}},t)$ and the action of the Verschiebung $V_n$ on
the Witt ring $W(\Z)$ is related to the Weil restriction of scalars from $\F_{q^n}$ to $\F_q$
(see \cite{Ram} for a precise statement). 

\smallskip

Recall that, if one denotes by $[a]$ the elements $[a]=(1-at)^{-1}$ in the Witt ring $W(R)$,
for $a\in R$, then the Frobenius ring homomorphisms $F_n: W(R)\to W(R)$ of the Witt ring
are determined by $F_n([a])=[a^n]$ and the Verschiebung group homomorphisms
$V_n: W(R)\to W(R)$ are defined on an arbitrary $P(t)\in W(R)$ as $F_n: P(t)\mapsto P(t^n)$.
These operations satisfy an analog of the Bost--Connes relations 
\begin{equation}\label{relsVnFn}
 F_n \circ F_m =F_{nm}, \ \  V_n\circ V_m = V_{nm}, \ \  \  F_n\circ V_n = n \cdot {\rm id}, \ \ \
F_n \circ V_m = V_m F_n \text{ if } (n,m)=1. 
\end{equation}
These correspond, respectively, to the semigroup structure of the $\sigma_n$ and $\tilde\rho_n$
of the integral Bost--Connes algebra and the relations $\sigma_n \circ \tilde\rho_n = n \cdot {\rm id}$,
while the last relation is determined in the Bost--Connes case by the commutation of the 
generators $\tilde\mu_n$ and $\mu_m^*$ for $(n,m)=1$. 

\smallskip

\begin{defn}\label{BCmotmeas} {\rm 
A factorizable motivic measure $\mu: K_0(\cV)\to R$, in the sense of Definition~\ref{factormotmeas},  
is of {\em Bost--Connes type}
if there is a lift to $K_0(\cV)$ of the Frobenius $F_n$ and Verschiebung $V_n$ 
of the Witt ring $W(R)$ to $K_0(\cV)$ such that the diagrams commute
$$ \xymatrix{  K_0(\cV) \ar[r]^{\zeta_\mu} \ar[d]^{\sigma_n} & W(R) \ar[d]^{F_n} \\
K_0(\cV) \ar[r]^{\zeta_\mu} & W(R)  } \ \ \  \ \ \
\xymatrix{  K_0(\cV) \ar[r]^{\zeta_\mu} \ar[d]^{\tilde\rho_n} & W(R) \ar[d]^{V_n} \\
K_0(\cV) \ar[r]^{\zeta_\mu} & W(R)  } $$
Such a motivic measure $\mu: K_0(\cV)\to R$ 
is of {\em homotopic Bost--Connes type} if the maps $\sigma_n$ and $\tilde\rho_n$
in the diagrams above also lift to endofunctors of the assembler category $\cC_\cV$ 
of the Grothendieck ring $K_0(\cV)$ with the endofunctors $\sigma_n$ compatible
with the monoidal structure. }
\end{defn}

\smallskip

\begin{defn}\label{EndFVdef} {\rm
The Frobenius and Verschiebung on the category $\cE^\pm_R$ are defined as the endofunctors
$F_n (E,f)=(E,f^n)$ and $V_n (E_\pm, f_\pm)=(E_\pm^{\oplus n}, V_n(f_\pm))$ 
with $V_n(f)$ defined by
\begin{equation}\label{VnEf}
 V_n: (E,f) \mapsto (E^{\oplus n}, V_n(f)), \ \ \ \ V_n(f)=\left( \begin{array}{ccccc}
0 & 0 & \cdots & 0 & f \\
1 & 0 & \cdots & 0 & 0 \\
0 & 1 & \cdots & 0 & 0 \\
\vdots & \vdots & \cdots & \vdots & \vdots \\
0 & 0 & \cdots & 1 & 0  \end{array}\right). 
\end{equation} }
\end{defn}

\smallskip

It is worth noting that the endofunctors of Definition~\ref{EndFVdef}
are akin to those used in the definitions of topological cyclic and topological 
restriction homology, \cite{HesMa}.

\smallskip

\begin{lem}\label{lemVnFn}
The Frobenius and Verschiebung $F_n$ and $V_n$ of Definition~\ref{EndFVdef} are 
endofunctors of the category $\cE^\pm_R$ with the property that
the maps they induce on $W_0(R)=K_0(\cE_R)/K_0(R)$ agree with the restrictions 
to $W_0(R)\subset W(R)$ of the Frobenius and Verschiebung maps. These
endofunctors determine natural transformations (still denoted $F_n$ and $V_n$) 
of the $\Gamma$-space $F_{\cE^\pm_R}: \Gamma^0\to \Delta_*$.
\end{lem} 

\proof The homomorphism $K_0(\cE_R) \to W_0(R)$ given by $$(E,f) \mapsto L(E,f)=\det(1-t M(f))^{-1}$$
sends the pair $(R,f_a)$ with $f_a$ acting on $R$ as multiplication by $a\in R$ to the element
$[a]=(1-at)^{-1}$ in the Witt ring. The action of the Frobenius $F_n([a])=[a^n]$ is induced from
the Frobenius $F_n (E,f)=(E,f^n)$ which is an endofunctor of $\cE_R$. This extends to a compatible
endofunctor of $\cE^\pm_R$ by $F_n(E_\pm,f_\pm)=(E_\pm,f_\pm^n)$. Similarly, the Verschiebung
map that sends $\det(1-t M(f))^{-1}\mapsto \det(1-t^n M(f))^{-1}$ is induced from the Verschiebung on 
$\cE_R$ given by \eqref{VnEf}, 
since we have $L(E^{\oplus n},V_n(f))=\det(1-t^n M(f))^{-1}$, with compatible endofunctors
$V_n (E_\pm, f_\pm)=(E_\pm^{\oplus n}, V_n(f_\pm))$ on $\cE^\pm_R$. The Frobeniius and
Verschiebung on $\cE^\pm_R$ induce natural transformations of the $\Gamma$-space
$F_{\cE^\pm_R}: \Gamma^0\to \Delta_*$ by composition of the summing functors $\Phi: \cP(X)\to \cE^\pm_R$
in $\Sigma_{\cE^\pm_R}(X)$ with the endofunctors $F_n$ and $V_n$ of $\cE^\pm_R$.
\endproof

\smallskip

\begin{prop}\label{zetamuBC}
Let $\mu: K_0(\cV)\to R$ be a factorizable motivic measure, as in Definition~\ref{factormotmeas}, 
that is of homotopical Bost--Connes type. Then the endofunctors $\sigma_n$ and $\tilde\rho_n$
of the assembler category $\cC_\cV$ determine natural transformations (still denoted by
$\sigma_n$ and $\tilde\rho_n$) of the associated $\Gamma$-space $F_\cV: \Gamma^0 \to \Delta_*$
that fit in the commutative diagrams
$$ \xymatrix{  F_\cV \ar[r]^{\Phi_\mu} \ar[d]^{\sigma_n} & F_{\cE^\pm_R} \ar[d]^{F_n} \\
F_\cV \ar[r]^{\Phi_\mu} & F_{\cE^\pm_R}  } \ \ \  \ \ \
\xymatrix{  F_\cV \ar[r]^{\Phi_\mu} \ar[d]^{\tilde\rho_n} & F_{\cE^\pm_R} \ar[d]^{V_n} \\
F_\cV \ar[r]^{\Phi_\mu} & F_{\cE^\pm_R}  } $$
where $\Phi_\mu: F_\cV \to F_{\cE^\pm_R}$ is the natural transformation of $\Gamma$-spaces
of \eqref{PhimuGamma} and $F_n$ and $V_n$ are the natural transformations of Lemma~\ref{lemVnFn}.
\end{prop}

\proof  The natural transformation $\Phi_\nu$ is determined as in Proposition~\ref{PhimuGamma} by
the functor $\Phi_\mu: \cC_\cV\to \cE^\pm_R$  that assigns $\Phi_\mu: X \mapsto
(E_\pm^X,f_\pm^X)$ constructed as in Lemma~\ref{AssEnd}. 
Suppose we have endofunctors $\sigma_n$ and $\tilde\rho_n$ of the assembler category $\cC_\cV$
that induce maps $\sigma_n$ and $\tilde\rho_n$ on $K_0(\cV)$ that lift the Frobenius and Verschiebung
maps of $W(R)$ through the zeta function $\zeta_\mu: K_0(\cV)\to W(R)$. This means that $\zeta_\mu(\sigma_n(X),t)=F_n\zeta_\mu(X,t)$ and $\zeta_\mu(\tilde\rho_n(X),t)=V_n\zeta_\mu(X,t)=\zeta_\mu(X,t^n)$. By Lemma~\ref{lemVnFn}, we
have $F_n\zeta_\mu(X,t) =L(F_n(E^X_\pm, f^X_\pm))=L(E_\pm^X,(f_\pm^X)^n)$ and 
$V_n \zeta_\mu(X,t) =L(V_n(E^X_\pm, f^X_\pm))=L((E^X_\pm)^{\oplus n}, V_n(f^X_\pm))$. 
This shows the compatibilities of the natural transformations in the diagrams above. 
\endproof

\medskip
\subsection{Spectra and spectra}\label{SpSpSec}  

We apply a construction similar to the one discussed in the previous subsections
to the case of the map $(X,f) \mapsto \sum_{\lambda\in \Sp(f_*)} m_\lambda \lambda$
that assigns to a variety over $\C$ with a quasi-unipotent map the spectrum of the
induced map $f_*$ in homology, seen as an element in $\Z[\Q/\Z]$, as in \S 6 of \cite{ManMar2}.

\smallskip

In this section the term spectrum will appear both in its homotopy theoretic sense
and in its operator sense. Indeed, we consider here a lift to the level of spectra 
(in the homotopy theoretic sense) of the construction
described in \S 6 of \cite{ManMar2}, based on the spectrum (in the operator sense) 
Euler characteristic.

\smallskip

We consider here a setting as in \cite{EbGuZa}, \cite{GuZa}, where $(X,f)$ is a pair 
of a variety over $\C$ and an endomorphism $f: X \to X$ such that the induced 
map $f_*$ on $H_*(X,\Z)$ has spectrum consisting of roots of unity. As discussed
in \cite{ManMar2} and in a related form in \cite{EbGuZa} the
spectrum determines a ring homomorphism (an Euler characteristic)
\begin{equation}\label{SpEul}
\sigma: K_0^{\Z}(\cV_\C) \to \Z[\Q/\Z] 
\end{equation}
where $K_0^{\Z}(\cV_\C)$ denotes the Grothendieck ring of pairs $(X,f)$
with the operations defined by the disjoint union and the Cartesian product. 
It is shown in \cite{ManMar2} that one can lift the operations $\sigma_n$
and $\tilde\rho_n$ of the integral Bost--Connes algebra from $\Z[\Q/\Z]$
to $K_0^{\Z}(\cV_\C)$ via the ``spectral Euler characteristic" \eqref{SpEul},
and that the operations can further be lifted from $K_0^{\Z}(\cV_\C)$ to
a (homotopy theoretic) spectrum with $\pi_0$ equal to $K_0^{\Z}(\cV_\C)$
via the assembler category construction of \cite{Zak1}. 

\smallskip

In the next sub section we discuss how to lift the right hand side of \eqref{SpEul},
namely the original Bost--Connes algebra $\Z[\Q/\Z]$ with the operations
$\sigma_n$ and $\tilde\rho_n$ to the level of a homotopy theoretic spectrum,
so that the spectral Euler characteristic \eqref{SpEul} becomes induced
by a map of spectra. 

\smallskip
\subsubsection{Bost--Connes Tannakian categorification and lifting of the
spectral Euler characteristic}\label{TannakaSec}   

To construct a categorification of the map \eqref{SpEul} compatible with the
Bost--Connes structure, we use the lift of the left-hand-side of \eqref{SpEul}
to an assembler category, as in Proposition~6.6 of \cite{ManMar2}, while for
the right-hand-side of \eqref{SpEul} we use the categorification of Bost--Connes 
system constructed in \cite{MaTa}. 

\smallskip

We begin by recalling the categorification of the Bost--Connes algebra of \cite{MaTa}.
Let ${\rm Vect}^{\bar\Q}_{\Q/\Z}(\Q)$ be the category of pairs $(W,\oplus_{r\in \Q/\Z}\bar W_r)$
with $W$ a finite dimensional $\Q$-vector space and $\oplus_r \bar W_r$ a $\Q/\Z$-graded vector
space with $\bar W =W\otimes \bar\Q$. This is a neutral Tannakian category with fiber functor
the forgetful functor $\omega: {\rm Vect}^{\bar\Q}_{\Q/\Z}(\Q) \to {\rm Vect}(\Q)$ and with
${\rm Aut}^\otimes(\omega)=\Sp(\bar\Q[\Q/\Z]^G)$ and $G={\rm Gal}(\bar\Q/\Q)$, see Theorem~3.2
of \cite{MaTa}. The category ${\rm Vect}^{\bar\Q}_{\Q/\Z}(\Q)$ is endowed with 
additive symmetric monoidal functors $\sigma_n(W)=W$ and $\overline{\sigma_n(W)}_r=\oplus_{r'\,:\, \sigma_n(r')=r} \bar W_{r'}$ if $r$ is in the range of $\sigma_n$ and zero otherwise and 
additive functors $\tilde\rho_n(W) = W^{\oplus n}$ and $\overline{\tilde \rho_n(W)}_r=\bar W_{\sigma_n(r)}$
satisfying $\sigma_n \circ \tilde\rho_n = n\cdot {\rm id}$ that induce the Bost--Connes maps on
$\Q[\Q/\Z]$. 

\smallskip

As shown in Theorem~3.18 of \cite{MaTa}, this category can be equivalently
described as a category of automorphisms ${\rm Aut}^{\bar\Q}_{\Q/\Z}(\Q)$ with objects pairs
$(W,\phi)$ of a $\Q$-vector space $V$ and a $G$-equivariant diagonalizable automorphism of 
$\bar W$ with eigenvalues that are roots of unity (seen as elements in $\Q/\Z$). There is an
equivalence of categories between ${\rm Vect}^{\bar\Q}_{\Q/\Z}(\Q)$ and ${\rm Aut}^{\bar\Q}_{\Q/\Z}(\Q)$ 
under which the functors $\sigma_n$ and $\tilde\rho_n$ correspond, respectively, to the Frobenius
and Verschiebung on ${\rm Aut}^{\bar\Q}_{\Q/\Z}(\Q)$, given by 
\begin{equation}\label{AutbarQFnVn}
 F_n : (W,\phi) \mapsto (W,\phi^n), \ \ \ \  V_n: (W,\phi)\mapsto (W^{\oplus n}, V_n(\phi)),
\end{equation} 
with
\begin{equation}\label{Vnphi}
V_n(\phi)=\left( \begin{array}{cccccc} 0 & 0 & \cdots & 0 & \phi \\
1 & 0 & \cdots & 0 & 0 \\
0 & 1 & \cdots & 0 & 0  \\
\vdots & & \vdots & & \vdots \\
0 & 0 & \cdots & 1 & 0 \\
\end{array}\right).  
\end{equation}
The equivalence is realized by mapping $(W,\phi)\mapsto (W,\oplus_r \bar W_r)$ where $\bar W_r$
are the eigenspaces of $\phi$ with eigenvalue $r\in \Q/\Z$.

\smallskip

\begin{rem}\label{2categorify}{\rm 
Conceptually, the first description of the categorification in terms of the Tannakian
category ${\rm Vect}^{\bar\Q}_{\Q/\Z}(\Q)$ is closer to the integral Bost--Connes algebra 
as introduced in \cite{CCM}, while its equivalent description in terms of ${\rm Aut}^{\bar\Q}_{\Q/\Z}(\Q)$ 
is closer to the reinterpretation of the Bost--Connes algebra in terms of Frobenius and Verschiebung 
operators, as in \cite{CoCo2}. Since we have introduced here the Bost--Connes algebra in the form
of \cite{CCM}, we are recalling both of these descriptions of the categorification, even though
in the following we will be using only the one in terms of ${\rm Aut}^{\bar\Q}_{\Q/\Z}(\Q)$. }
\end{rem}

\smallskip

\begin{prop}\label{AssTann}
Let $\cC^\Z_\C$ be the assembler category underlying $K_0^{\Z}(\cV_\C)$, as in 
Proposition~6.6 of \cite{ManMar2}.
The assignment $\Phi(X,f) = (H_*(X,\Q), \oplus_r E_r(f_*))$, where $E_r(f_*)$ is the eigenspace with
eigenvalue $r\in \Q/\Z$, determines a functor $\Phi: \cC^\Z_\C \to {\rm Aut}^{\bar\Q}_{\Q/\Z}(\Q)$
that lifts the Frobenius and Vershiebung functors on ${\rm Aut}^{\bar\Q}_{\Q/\Z}(\Q)$ to the
endofunctors $\sigma_n$ and $\tilde\rho_n$ of $\cC^\Z_\C$ implementing the Bost--Connes
structure.
\end{prop}

\proof
We can construct the functor from the assembler category $\cC^\Z_\C$ of \S 6 of \cite{ManMar2},
underlying $K_0^{\Z}(\cV_\C)$ to ${\rm Aut}^{\bar\Q}_{\Q/\Z}(\Q)$ by following along the lines of 
Lemma~\ref{AssEnd} and Proposition~\ref{PhimuGamma}, where we assign 
$\Phi(X,f) = (H_*(X,\Q), \oplus_r E_r(f_*))$ where $E_r(f_*)$ is the eigenspace with
eigenvalue $r\in \Q/\Z$. The Bost--Connes algebra then lifts to the Frobenius and Vershiebung
functors on ${\rm Aut}^{\bar\Q}_{\Q/\Z}(\Q)$ and the latter lift to geometric Frobenius and
Verschiebung operations on the pairs $(X,f)$ mapping to $(X,f^n)$ and to $(X\times Z_n,\Phi_n(f))$.
\endproof

\smallskip

This point of view, that replaces the Bost--Connes algebra with it categorification
in terms of the Tannakian category  ${\rm Aut}^{\bar\Q}_{\Q/\Z}(\Q)$ as in \cite{MaTa}
will also be useful in Section~\ref{NoriSec}, where we reformulate our categorical setting,
by passing from Grothendieck rings, assemblers and spectra, to Tannakian categories
of Nori motives, and we compare in Lemma~\ref{TDV} and Theorem~\ref{BCNori}
the categorification of the Bost--Connes algebra obtained via Nori motives with
the one of \cite{MaTa} recalled here.

\section{Bost--Connes systems in categories of Nori motives}\label{NoriSec}  

We introduce in this section a motivic framework, with Bost--Connes type systems that
on Tannakian categories of motives.  The main result in this part of
the paper will be Theorem~\ref{BCNori}, showing the existence of
a fiber functor from the Tannakian category of Nori motives with good effectively
finite $\hat\Z$-action to the Tannakian category ${\rm Aut}^{\bar \Q}_{\Q/\Z}(\Q)$
that lifts the Bost--Connes system given by Frobenius and Verschiebung on
the target category to a Bost--Connes system on Nori motives. 
Proposition~\ref{motsheavesprop} then extends this Bost--Connes structure to
the relative case of motivic sheaves. 

\smallskip

This is a natural generalization of the approach to Grothendieck rings 
via assemblers, which can be extended in an interesting way to the domain
of motives, namely, Nori motives.

\smallskip

Roughly speaking, the theory of Nori motives starts with 
lifting the relations   
\begin{equation}\label{relK0S}
 [ f: X \to S ] = [f|_Y : Y \to S ] + [ f|_{X\smallsetminus Y}: X\smallsetminus Y \to S] 
\end{equation} 
of (relative) Grothendieck rings $K_0(\cV_S)$ to the level of ``diagrams'', which intuitively
can be imagined as ``categories without multiplication of morphisms.''

\smallskip

\smallskip
\subsection{Nori diagrams}\label{NoriDiagSubsec}  

More precisely, (cf.~Definition~7.1.1 of \cite{HuM-S17}, p.~137), we have the following
definitions. 

\begin{defn}\label{DiagDef}{\rm
A {\rm diagram} (also called {\rm a quiver}) $D$ is a family consisting of a
set of vertices $V(D)$ and a set of oriented edges, $E(D)$.
Each edge $e$ either connects  two different vertices, going, say, from a vertex
$\partial_{\text{out}}e=v_1$ to a vertex $\partial_{\text{in}}e=v_2$,
or else is ``an identity'', starting and ending
with one and the same  vertex $v$. We will consider only diagrams
with one identity for each vertex.  }
\end{defn}

\smallskip

Diagrams can be considered
as objects of  a category, with obvious morphisms.

\smallskip

\begin{defn}\label{DiagRepDef}{\rm 
Each small category $\cC$ can be considered as a diagram $D(\cC)$, with 
$V(D(\cC))= {\rm Ob}\, \cC$, $E(D(\cC))= {\rm Mor}\, \cC$,
so that each morphism $X\to Y$ ``is'' an oriented edge from $X$ to $Y$.
More generally, {\it a representation}  $T$ of a diagram $D$ in a (small) category 
$\cC$ is a morphism of directed graphs $T:\, D\to D(\cC)$. }
\end{defn}

\smallskip

Notice that a considerably more general treatment of graphs with markings,
including diagrams etc. in the operadic environment, can be found in \cite{MaBo07}. 
We do not use it here, although it might be highly relevant.

\smallskip
\subsection{From geometric diagrams to Nori motives}\label{DiagNoriSec}   

We recall the main idea in the construction of Nori motives from geometric diagrams.
For more details, see \cite{HuM-S17}, pp.~140--144.

\smallskip

\begin{enumerate}

\item {\it Start} with  the  following data:

\begin{itemize}
\item[a)] a diagram $D$;

\item[b)] a noetherian commutative ring with unit $R$ and the category of finitely generated $R$--modules
$R\hbox{-Mod}$;

\item[c)] a representation $T$ of $D$ in $R\hbox{-Mod}$, in the sense of Definition~\ref{DiagRepDef}.
\end{itemize}

\medskip

\item {\it Produce}  from them the category $C(D,T)$ defined in the following way:

\begin{itemize}
\item[d1)] If $D$ is finite, then $C(D,T)$ is the category of finitely generated $R$--modules
equipped with an $R$--linear action of $End (T)$.
\smallskip
\item[d2)] If $D$ is infinite, first consider its  all finite subdiagrams $F$.

\medskip

\item[d3)] For each $F$ construct $C(F, T|_F)$ as in d1). Then apply the following limiting procedure:
$$ C(D,T):={\rm colim}_{F\subseteq D\, \text{finite}} \,\, C(F, T|_F)\, $$
Thus, the category $C(D,T)$ has the following structure:

\begin{itemize}

\item Objects of $C(D,T)$ will be all objects of the categories $C(F,T|_F)$. If $F\subset F^{\prime}$,
then each object $X_F$ of $C(F, T|_F)$ can be canonically  extended to an object of 
$C(F^{\prime},T|_{F^{\prime}})$.

\item Morphisms from $X$ to $Y$ in $C(D,T)$ will be  defined as colimits over $F$
of morphisms from $X_F$ to $Y_F$ with respect to these extensions.
\end{itemize}

\item[d4)] The fact that $C(D,T)$ has a functor to $R\hbox{-Mod}$ follows directly from
the definition and the finite case.

\end{itemize}

\medskip

\item {\it The result} is called {\it the diagram category $C(D,T)$}.

It is an $R$--linear abelian category which
is endowed with $R$--linear faithful exact forgetful functor 
$$
f_T:\,C(D,T)\to R\hbox{-Mod}.
$$  
\end{enumerate}

\medskip

\subsubsection{Universal diagram category} 
The following results explain why abstract diagram
categories play a central role in the formalism
of Nori motives: they formalise the Grothendieck intuition of
motives as objects of the universal cohomology theory.

\smallskip

\begin{thm}\label{DiagCat} {\rm {\bf \cite{HuM-S17}}}
\begin{itemize}
  \item[(i)]  Any representation $T:\,D\to R\hbox{-Mod}$
can be presented as post-composition of the forgetful functor $f_T$
with an appropriate representation $\tilde{T}:\,D\to C(D,T)$:
$$
T = f_T \circ \tilde{T} .
$$
with the following universal property:

\smallskip

Given any $R$--linear abelian category $A$ with a representation $F:\,D\to A$
and $R$--linear faithful exact  functor $f:\,A\to R\hbox{-Mod}$
with $T=f\circ F$, it factorizes through a faithful exact functor
$L(F):\, C(D,T) \to A$ compatibly with the decomposition
$$
T = f_T \circ \tilde{T} .
$$

\smallskip

\item[(ii)] The functor $L(F)$ is unique up to unique isomorphism of exact additive functors.
\end{itemize}
\end{thm}

For proofs, cf.~\cite{HuM-S17}, pp. 140--141 and p.~167.

\medskip

\subsubsection{Nori geometric diagrams} \label{NoriGDiagSec}
If we start not with an abstract category but with  a ``geometric" category $\cC$ (in the sense that
its objects are spaces/varieties/schemes, possibly endowed with additional
structures), in which one can define morphisms of closed 
embeddings $Y\hookrightarrow X$ (or $Y\subset X$) and morphisms of complements to closed embeddings
$X\setminus Y \to X$, we can  define the Nori diagram of
{\it effective pairs} $D(\cC)$ in the following way (see \cite{HuM-S17}, 
pp.~207--208).

\smallskip

\begin{itemize}
\item[a)] One vertex of $D(\cC)$ is a triple $(X,Y,i)$ where $Y\hookrightarrow X$ is a closed
embedding, and $i$ is an integer.

\item[b)] Besides obvious identities, there are edges of two types.

\smallskip

\item[b1)] Let $(X,Y)$ and $(X^{\prime}, Y^{\prime})$ be two pairs
of closed embeddings.
Every morphism $f: X\to X^{\prime}$ such that $f(Y)\subset Y^{\prime}$
produces  functoriality edges $f^*$ (or rather $(f^*,i)$) going from $(X^{\prime}, Y^{\prime},i)$
to $(X, Y, i)$.

\item[b2)] Let $(Z\subset Y\subset X)$ be a stair of closed embeddings. Then
it defines coboundary edges  $\partial$ from $(Y,Z,i)$ to $(X,Y, i+1)$.
\end{itemize}

\medskip

\subsubsection{(Co)homological representatons of Nori geometric diagrams} If we start not just from
the initial category of spaces $\cC$, but rather from a pair $(\cC, H)$ where $H$ is a cohomology
theory,  then assuming
reasonable properties of this pair, we can define the respective representation $T_H$ of
$D(\cC)$ that we will call  a {\it (co)homological representation of $D(\cC)$}.

\smallskip

For a survey of  such pairs $(\cC, H)$ that were studied in the context
of Grothendieck's motives, see \cite{HuM-S17}, pp. 31--133. The relevant cohomology
theories include, in particular, singular cohomology, and algebraic and holomorphic de Rham cohomologies.

\smallskip

Below we will consider the basic example of cohomological representations of Nori diagrams
that leads to Nori motives.

\subsubsection{Effective Nori motives}  We follow \cite{HuM-S17},  pp. 207--208.
Take as a category $\cC$, the 
starting object in the definition of Nori geometric diagrams above, the category $\cV_k$ of
varieties $X$ defined over a subfield $k\subset \bC$.

\smallskip

We can then define the Nori diagram $D(\cC)$ as above. This diagram will be
denoted ${\rm Pairs}^{eff}$ from now on, 
$$ {\rm Pairs}^{eff} = D(\cV_k). $$

\smallskip

The category of effective mixed Nori motives is the diagram category $C({\rm Pairs}^{eff}, H^*)$ where
 $H^i(X,\bZ)$ is the respective singular
cohomology of the analytic space $X^{an}$ (cf.~\cite{HuM-S17}, pp.~31--34 and further on).

\smallskip

It turns out (see~\cite{HuM-S17}, Proposition 9.1.2. p.~208) that the map
$$
H^*:\, {\rm Pairs}^{eff} \to \bZ\hbox{-Mod}
$$
sending $(X,Y,i)$ to the relative singular cohomology $H^i(X(\bC), Y(\bC); \bZ)$,
 naturally extends  to a representation
of the respective Nori diagram in the category of finitely generated abelian groups $\bZ\hbox{-Mod}$.

\smallskip
\subsection{Category of equivariant Nori motives}\label{NoriTannSec}     

We now introduce the specific category of Nori motives that we will be
using for the construction of the associated Bost--Connes system.

\smallskip

Let $D(\cV)$ the Nori geometric diagrams associated to the category $\cV$ of 
varieties over $\Q$, constructed as described in \S \ref{DiagNoriSec}. 

\smallskip

As in \cite{ManMar2} and in \S~\ref{RelGrothSec} of this paper, we
consider here the category $\cV^{\hat\Z}$ of varieties $X$ with a 
good effectively finite action of $\hat\Z$.
We can view $\cV^{\hat\Z}$ as an enhancement $\hat\cV$ of the category $\cV$,
in the sense described in \S \ref{EnrichSec}. 

\smallskip

Define the Nori diagram of {\it effective pairs} $D(\cV^{\hat\Z})$ as we recalled 
earlier in \S \ref{DiagNoriSec}:

\begin{itemize}
\item[a)]  One vertex of $D(\cV^{\hat\Z})$ is a triple $((X,\alpha_X), (Y,\alpha_Y),i)$, of varieties $X$ and $Y$
with good effectively finite $\hat\Z$ actions, $\alpha_X: \hat\Z\times X \to X$ and
$\alpha_Y: \hat\Z \times Y \to Y$, and an integer $i$, together with a
closed embedding $j:\,Y\hookrightarrow X$ that is equivariant with respect to the $\hat\Z$ actions.
For brevity, we will denote such a triple $(\hat{X}, \hat{Y}, i)$ and call it a closed embedding in 
the enhancement $\hat{\cV}$. 

\item[b)] Identity edges, functoriality edges, and coboundary edges
are obvious enhancements of the respective edges defined in \S \ref{DiagNoriSec},
with the requirement that all these maps are $\hat\Z$-equivariant.

\item[b1)] Let $(\hat{X},\hat{Y})$ and $(\hat{X}^{\prime}, \hat{Y}^{\prime})$ be two pairs
of closed embeddings in $\hat{\cV}$.
Every morphism $f: X\to X^{\prime}$ such that $f(Y)\subset Y^{\prime}$ and $f\circ \alpha_X = \alpha_{X^\prime}\circ f$ produces  functoriality edges $f^*$ (or rather $(f^*,i)$) 
going from $((X^{\prime}, \alpha_{X^{\prime}}), (Y^{\prime}, \alpha_{Y^{\prime|}}),i)$
to $(X, Y, i)$.

\smallskip

\item[b2)] Let $(Z\subset Y\subset X)$ be a stair of closed embeddings compatible with enhancements (equivariant with respect tot the $\hat\Z$-actions). Then
it defines coboundary edges  $\partial$  
$$
((Y,\alpha_Y),(Z,\alpha_Z),i) \to ((X,\alpha_X), (Y,\alpha_Y),  i+1).
$$
\end{itemize}

We have thus defined he Nori geometric diagram of enhanced effective pairs, which
we denote equivalently by $D(\hat{\cV})$ or $D(\cV^{\hat\Z})$. 
\smallskip

Notice that forgetting in this diagram all enhancements, we obtain the map
$D(\hat{\cV})\to D(\cV)$ which is {\it injective} both on vertices and edges.

\smallskip
\subsection{Bost--Connes system on Nori motives}\label{BCNorisubsec}  

We now construct a Bost--Connes system on a category of Nori motives
obtained from the diagram $D(\cV^{\hat\Z})$ described above, which lifts
to the level of motives the categorification of the Bost--Connes algebra
constructed in \cite{MaTa}. 

\smallskip

As we recalled in \S \ref{TannakaSec}, we can describe the
categorification of the Bost--Connes algebra of \cite{MaTa} in terms of the Tannakian
category ${\rm Vec}^{\bar\Q}_{\Q/\Z}(\Q)$ with suitable functors $\sigma_n$ and
$\tilde\rho_n$ constructed as in Theorem~3.7 of \cite{MaTa} or in terms of an
equivalent Tannakian category ${\rm Aut}^{\bar\Q}_{\Q/\Z}(\Q)$ endowed with
Frobenius and Verschiebung functors. We are going to use here the second description. 

\smallskip

\begin{lem}\label{TDV}
The assignment $T: ((X,\alpha_X), (Y,\alpha_Y),i)\mapsto H^i(X(\C), Y(\C),\Q)$
determines a representation $T: D(\cV^{\hat\Z}) \to {\rm Aut}^{\bar\Q}_{\Q/\Z}(\Q)$
of the diagram $D(\cV^{\hat\Z})$ constructed above.
\end{lem}

\proof As discussed in the previous subsection, we view elements $((X,\alpha_X), (Y,\alpha_Y),i)$
of $D(\cV^{\hat\Z})$ in terms of an enhancement $\hat\cV$ of the category $\cV$ defined as in
\S \ref{EnrichSec}, by choosing a primitive root of unity that generates the cyclic group $\Z/N\Z$,
so that the actions $\alpha_X$ and $\alpha_Y$ are determined by self maps $v_X$ and $v_Y$
as in \S \ref{EnrichSec}. We identify the element above with $((X,v_X), (Y,v_Y), i)$, which we
also denoted by $(\hat X, \hat Y, i)$ in the previous subsection. Since the embedding
$Y\hookrightarrow X$ is $\hat\Z$-equivariant, the map $v_Y$ is the restriction to $Y$ of the
map $v_X$ under this embedding. We denote by $\phi^i$ the induced map on the cohomology
$H^i(X(\C), Y(\C),\Q)$. The eigenspaces of $\phi^i$ are the subspaces of the decomposition of
$H^i(X(\C), Y(\C),\Q)$ according to characters of $\hat\Z$, that is, elements in 
$\Hom(\hat\Z, \C^*)=\nu^*\simeq \Q/\Z$. Thus, we obtain an object 
$(H^i(X(\C), Y(\C),\Q), \phi^i)$ in the category ${\rm Aut}^{\bar\Q}_{\Q/\Z}(\Q)$.  Edges in the diagram are
$\hat\Z$-equivariant maps so they induce morphisms between the corresponding objects 
in the  category ${\rm Aut}^{\bar\Q}_{\Q/\Z}(\Q)$.
\endproof

\smallskip

One can also see in a similar way that the fiber functor 
$T: ((X,\alpha_X), (Y,\alpha_Y),i)\mapsto H^i(X(\C), Y(\C),\Q)$ determines an object in the category ${\rm Vec}^{\bar\Q}_{\Q/\Z}(\Q)$. Indeed, the pair $(X,Y)$ with $Y\subset X$ is endowed with compatible 
good effectively finite $\hat\Z$-actions $\alpha_X$ and $\alpha_Y$, hence the singular 
cohomology $H^i(X(\C), Y(\C),\Q)$ 
carries a resulting $\hat\Z$-representation. Thus, the vector space $H^i(X(\C), Y(\C),\bar\Q)$
can be decomposed into eigenspaces of this representations according to characters
$\chi \in \Hom(\hat\Z, \bG_m)=\Q/\Z$. Thus, we obtain a decomposition of
$H^i(X(\C), Y(\C),\bar\Q) =\oplus_{r\in \Q/\Z} \bar V_r$ as a $\Q/\Z$-graded vector
space. We choose to work with the category ${\rm Aut}^{\bar\Q}_{\Q/\Z}(\Q)$ because the
Bost--Connes structure is more directly expressed in terms of Frobenius and Verschiebung,
which will make the lifting of this structure to the resulting category of Nori motives more 
immediately transparent, as we discuss below.

\smallskip

The representation $T: D(\cV^{\hat\Z}) \to {\rm Aut}^{\bar\Q}_{\Q/\Z}(\Q)$ replaces, at this
motivic level, our previous use in \cite{ManMar2} of the equivariant Euler characteristics 
$K_0^{\hat\Z}(\cV)\to \Z[\Q/\Z]$ (see \cite{Looij}) as a way to lift the Bost--Connes algebra. We proceed in
the following way to obtain the Bost--Connes structure in this setting. 

\smallskip

\begin{defn}\label{intertwineTann}{\rm
Let $D$ be a diagram, endowed with a representation $T: D \to {\rm Aut}^{\bar\Q}_{\Q/\Z}(\Q)$, and let 
$\cC(D,T)$ be the associated diagram category, obtained as in \S \ref{DiagNoriSec}, with the induced
functor $\tilde T: \cC(D,T) \to {\rm Aut}^{\bar\Q}_{\Q/\Z}(\Q)$. 
We say that the functor $\tilde T$ intertwines the Bost--Connes structure, if there are endofunctors
$\sigma_n$ and $\tilde\rho_n$ of $\cC(D,T)$ (where the $\sigma_n$ but not the $\tilde\rho_n$ are
compatible with the tensor product structure) such that the following diagrams commute,
$$ \xymatrix{ \cC(D,T) \ar[r]^{\tilde T} \ar[d]^{\sigma_n} & {\rm Aut}^{\bar\Q}_{\Q/\Z}(\Q) \ar[d]^{F_n} \\
\cC(D,T) \ar[r]^{\tilde T} & {\rm Aut}^{\bar\Q}_{\Q/\Z}(\Q) }  \ \ \  \ \ \  
 \xymatrix{ \cC(D,T) \ar[r]^{\tilde T}  & {\rm Aut}^{\bar\Q}_{\Q/\Z}(\Q)  \\
\cC(D,T) \ar[r]^{\tilde T}\ar[u]^{\tilde\rho_n} & {\rm Aut}^{\bar\Q}_{\Q/\Z}(\Q) \ar[u]^{V_n} } $$
where on the right-hand-side of the diagrams, the $F_n$ and $V_n$ are the Frobenius and Verschiebung
on ${\rm Aut}^{\bar\Q}_{\Q/\Z}(\Q)$, defined as in \eqref{AutbarQFnVn} and \eqref{Vnphi}.
}\end{defn}

\smallskip

\begin{defn}\label{sigmarhonNori}{\rm
For $((X,\alpha_X), (Y,\alpha_Y), i)$ in the category $\cC(D(\cV^{\hat\Z}),T)$ of Nori motives
associated to the diagram $D(\cV^{\hat\Z})$ define
\begin{equation}\label{sigmanNori}
\sigma_n: ((X,\alpha_X), (Y,\alpha_Y), i) \mapsto ((X,\alpha_X\circ \sigma_n),
(Y,\alpha_Y\circ \sigma_n), i)
\end{equation}
\begin{equation}\label{rhonNori}
\tilde\rho_n : ((X,\alpha_X), (Y,\alpha_Y), i) \mapsto (X\times Z_n,  
\Phi_n(\alpha_X), (Y\times Z_n, \Phi_n(\alpha_Y)), i) ,
\end{equation}
where $Z_n={\rm Spec}(\Q^n)$ and $\Phi_n(\alpha)$ is the geometric Verschiebung
defined as in \S \ref{VerschSec}. 
}\end{defn}

\smallskip

\begin{thm}\label{BCNori}
The $\sigma_n$ and $\tilde\rho_n$ of \eqref{sigmanNori} and \eqref{rhonNori}
determine a Bost--Connes system on the category $\cC(D(\cV^{\hat\Z}),T)$ of Nori motives
associated to the diagram $D(\cV^{\hat\Z})$. The representation 
$T: D(\cV^{\hat\Z}) \to {\rm Aut}^{\bar\Q}_{\Q/\Z}(\Q)$ constructed above has the property
that the induced functor $$\cC(D(\cV^{\hat\Z}),T) \to {\rm Aut}^{\bar\Q}_{\Q/\Z}(\Q)$$ intertwines
the endofunctors $\sigma_n$ and $\tilde\rho_n$ of the Bost--Connes system on 
$\cC(D(\cV^{\hat\Z}),T)$ and the Frobenius $F_n$ and Verschiebung $V_n$ of
the  Bost--Connes structure on ${\rm Aut}^{\bar\Q}_{\Q/\Z}(\Q)$.
\end{thm}

\proof 
Consider the mappings $\sigma_n$ and $\tilde\rho_n$ defined in \eqref{sigmanNori} and \eqref{rhonNori},
The effect of the transformation $\sigma_n$, when written
in terms of the data $((X,v_X), (Y,v_Y),i)$ is to send $v_X \mapsto v_X^n$ and $v_Y\mapsto v_Y^n$,
hence it induces the Frobenius map $F_n$ acting on $(H^i(X(\C), Y(\C),\Q), \phi^i)$ in 
${\rm Aut}^{\bar\Q}_{\Q/\Z}(\Q)$. Similarly, we have
$T(X\times Z_n,  \Phi_n(\alpha_X), (Y\times Z_n, \Phi_n(\alpha_Y)), i)=H^i(X\times Z_n, Y\times Z_n, \Q)$
where by the relative version of the K\"unneth formula we have 
$(H^i(X(\C)\times Z_n (\C), Y(\C)\times Z_n(\C), \Q)\simeq H^i(X(\C), Y(\C), \Q)^{\oplus n}$ with the
induced map $V_n(\phi^i)$. The maps $\sigma_n$ and $\tilde\rho_n$ defined as above determine
self maps of the diagram $D(\cV^{\hat\Z})$. By Lemma~7.2.6 of \cite{HuM-S17} given a map 
$F: D_1\to D_2$ of diagrams and a representation $T: D_2 \to R\hbox{-Mod}$, there is an $R$-linear
exact functor $\cF: C(D_1, T\circ F) \to C(D_2, T)$ such that the following diagram commutes:
$$ \xymatrix{ D_1 \ar[rr]^{F} \ar[d] &  & D_2 \ar[d] \\
C(D_1, T\circ F) \ar[rr]^{\cF} \ar[dr] & & C(D_2,T) \ar[dl] \\
& R\hbox{-Mod} & } $$
We still denote by $\sigma_n$ and $\tilde\rho_n$ the endofunctors induced in 
this way on $C(D(\cV^{\hat\Z}), T)$. To check the compatibility of the $\sigma_n$ functors
with the monoidal structure, we use the fact that for Nori motives the product structure is
constructed using ``good pairs" (see \S 9.2.1 of \cite{HuM-S17}), that is, elements
$(X,Y,i)$ with the property that $H^j(X,Y,\Z)=0$ for $j\neq i$. For such elements the
product is given by $(X,Y,i) \times (X',Y', j)= (X\times X', X\times Y' \cup Y \times X', i+j)$.
The diagram category $C({\rm Good}^{eff},T)$ obtained
by replacing effective pairs ${\rm Pairs}^{eff}$ with good effective pairs ${\rm Good}^{eff}$
is equivalent to $C({\rm Pairs}^{eff},T)$ (Theorem 9.2.22 of \cite{HuM-S17}), hence the
tensor structure defined in this way on  $C({\rm Good}^{eff},T)$ determines the tensor structure
of $C({\rm Pairs}^{eff},T)$ and on the resulting category of Nori motives, 
see \S 9.3 of \cite{HuM-S17}. Thus, to check the compatibility of the functors $\sigma_n$
with the tensor structure it suffices to see that on a product of good pairs, where indeed we have
$$ \sigma_n((X,\alpha_X), (Y,\alpha_Y), i) \times \sigma_n((X',\alpha_X'), (Y',\alpha_Y'), j) = $$
$$ ((X \times X', (\alpha_X\times \alpha_X')\circ\Delta \circ \sigma_n),
((X\times Y',(\alpha_X\times \alpha_Y')\circ\Delta\circ \sigma_n)    \cup (Y \times X', (\alpha_Y\times \alpha_X') \circ \Delta \circ \sigma_n)), i+j) $$
$$ = \sigma_n( ( (X,\alpha_X), (Y,\alpha_Y), i) \times  ((X',\alpha_X'), (Y',\alpha_Y'), j)). $$
The functors $\tilde\rho_n$ are not compatible with the tensor product structure, as expected.
\endproof

\smallskip

\begin{rem}\label{MotBC}{\rm
In \cite{MaTa} a motivic interpretation of the categorification of the Bost--Connes algebra is
given by identifying the Tannakian category ${\rm Vec}^{\bar\Q}_{\Q/\Z}(\Q)$ with a limit of
orbit categories of Tate motives. Here we presented a different motivic categorification of
the Bost--Connes algebra by lifting the Bost--Connes structure to the level of the
category of Nori motives. In \cite{MaTa} a motivic Bost--Connes structure was also
constructed using the category of motives over finite fields and the larger class of Weil 
numbers replacing the roots of unity of the Bost--Connes system. } \end{rem}

\smallskip
\subsection{Motivic sheaves and the relative case}\label{NoriSheavesSec}   

The argument presented in Theorem~\ref{BCNori} lifting the Bost--Connes
structure to the category of Nori motives, which provides a Tannakian category
version of the list to Grothendieck rings via the equivariant Euler characteristics 
$K_0^{\hat\Z}(\cV)\to \Z[\Q/\Z]$, can also be generalized to the relative setting,
where we considered the Euler characteristic
$$  \chi_S^{\hat\Z}: K_0^{\hat\Z}(\cV_S) \to K_0^{\hat\Z}(\Q_S) $$
with values in the Grothendieck ring of constructible sheaves, discussed 
in \S \ref{RelGrothSec} of this paper.
The categorical setting of Nori motives that is appropriate for this relative case is
the Nori category of motivic sheaves introduced in \cite{Arapura}. 

\smallskip

We recall here briefly the construction of the category of motivic sheaves of \cite{Arapura}
and we show that the Bost--Connes structure on the category of Nori motives described
in Theorem~\ref{BCNori} extends to this relative setting.

\smallskip

Consider pairs $(X\to S, Y)$ of varieties over a base $S$ 
with $Y \subset X$ endowed with the restriction $f_Y: Y \to S$.
Morphisms $f: (X\to S, Y) \to (X' \to S, Y')$ are morphisms of varieties $h: X\to X'$
satisfying the commutativity of
$$ \xymatrix{ X \ar[rr]^{h} \ar[dr]^f & & X' \ar[dl]_{f'} \\ & S & } $$
and such that $h(Y)\subset Y'$. As before, we consider varieties endowed
with good effectively finite $\hat\Z$-action.
We denote by $(S,\alpha)$ the base with its good effectively finite $\hat\Z$-action and
by $((X\alpha_X)\to (S,\alpha), (Y,\alpha_Y))$ the pairs as above where we
assume that the map $f: X \to S$ and the inclusion $Y\hookrightarrow X$ are 
$\hat\Z$-equivariant.

\smallskip

Following \cite{Arapura}, a diagram $D(\cV_S)$ is obtained by considering as
vertices elements of the form $(X \to S, Y, i,w)$ with $(X \to S, Y)$ a pair as above,
$i\in \bN$ and $w\in \Z$. The edges are given by the three types of edges
\begin{enumerate}
\item geometric morphisms $h: (X\to S, Y) \to (X'\to S, Y')$ as above
determine edges  $h^* : (X'\to S, Y', i,w) \to (X\to S, Y, i,w)$;
\item connecting morphisms $\partial: (Y\to S,Z,i,w) \to (X\to S, Y, i+1,w)$ for a chain 
of inclusions $Z\subset Y \subset X$;
\item twisted projections: 
$(X,Y,i,w) \to (X\times \P^1, Y\times \P^1 \cup X \times \{ 0 \}, i+2, w+1)$.
\end{enumerate}
For consistency with our previous notation we have here written the
morphisms in the contravariant (cohomological) way rather than in the
covariant (homological) way used in \S 3.3 of \cite{Arapura}. 

\smallskip

Note that in the previous section, following \cite{HuM-S17} we described
the effective Nori motives as $\cM\cN^{eff}=C({\rm Pairs}^{eff},T)$, with
the category of Nori motives $\cM\cN$ being then obtained as the localization
of $\cM\cN^{eff}$ at $(\bG_m,\{ 1 \},1)$ (inverting the Lefschetz motive). 
Here in the setting of \cite{Arapura} the Tate motives are accounted for in
the diagram construction by the presence of the twist $w$ and the last class of edges.

\smallskip

Given $f: X \to S$ and a sheaf $\cF$ on $X$ one has $H^i_S(X;\cF)=R^i f_* \cF$. 
In the case of a pair $(f:X\to S, Y)$, let $j: X\smallsetminus Y \hookrightarrow X$ be
the inclusion and consider $H^i_S(X,Y;\cF)=R^i f_* j_{!} \cF |_{X\smallsetminus Y}$. 
The diagram representation $T$ in this case maps $T(X \to S, Y, i,w) = H^i_S(X,Y, \cF)(w)$
to the (Tate twisted) constructible sheaf $H^i_S(X,Y; \cF)$. 
It is shown in \cite{Arapura} that the Nori formalism of geometric diagrams applies
to this setting and gives rise to a Tannakian category of motivic sheaves $\cM\cN_S$.
In particular one considers the case where $\cF$ is constant with $\cF =\Q$, so that
the diagram representation $T:  D(\cV_S) \to \Q_S$
and the induced functor on $\cM\cN_S$ replace at the motivic level the 
Euler characteristic map on the relative Grothendieck ring $K_0(\cV_S) \to K_0(\Q_S)$.

\smallskip

As in the previous cases, we consider an enhancement of this category of
motivic sheaves, in the sense of \S \ref{EnrichSec}, by introducing good effectively finite
$\hat\Z$-actions. We modify the construction of \cite{Arapura} in the following way. 

\smallskip

We consider a diagram $D(\cV_{(S,\alpha)}^{\hat\Z})$ where the vertices are
elements $$((X,\alpha_X) \to (S,\alpha), (Y,\alpha_Y), i,w)$$ so that
the maps $f: X\to S$ and the inclusion $Y\hookrightarrow X$ are
$\hat\Z$-equivariant, and with morphisms as above, where all the
maps are required to be compatible with the $\hat\Z$-actions. 
One obtains by the same procedure as in \cite{Arapura} a category of
equivariant motivic sheaves $\cM\cN_S^{\hat\Z}$.
The representation above maps $D(\cV_{(S,\alpha)}^{\hat\Z})$ to $\hat\Z$-equivariant
constructible sheaves over $(S,\alpha)$.  
Then the same argument we used \S \ref{RelGrothSec}
at the level of Grothendieck rings, assemblers and spectra applies to this setting and gives
the following result.

\smallskip

\begin{prop}\label{motsheavesprop}
The maps of diagrams 
$$  \sigma_n : D(\cV_{(S,\alpha)}^{\hat\Z}) \to D(\cV_{(S,\alpha\circ \sigma_n)}^{\hat\Z}) $$
$$ \tilde\rho_n : D(\cV_{(S,\alpha)}^{\hat\Z}) \to D(\cV_{(S \times Z_n, \Phi_n(\alpha))}^{\hat\Z}) $$
given by
$$ \sigma_n ((X,\alpha_X) \to (S,\alpha), (Y,\alpha_Y), i,w) = 
((X,\alpha_X \circ \sigma_n) \to (S,\alpha\circ \sigma_n), (Y,\alpha_Y\circ \sigma_n), i,w) $$
$$ \begin{array}{c} \tilde\rho_n ((X,\alpha_X) \to (S,\alpha), (Y,\alpha_Y), i,w) = \\[3mm]
((X \times Z_n, \Phi_n(\alpha_X)) \to (S \times Z_n,\Phi_n(\alpha)), (Y \times Z_n,\Phi_n(\alpha_Y)), i,w) 
\end{array} $$
determine functors of the resulting category of motivic sheaves $\cM\cN_S^{\hat\Z}$
such that $\sigma_n\circ \tilde\rho_n = n\, {\rm id}$ and $\tilde\rho_n \circ \sigma_n$
is a product with $(Z_n,\alpha_n)$.  Thus, one obtains on
the category $\cM\cN_S^{\hat\Z}$ a Bost--Connes system as in Definition~\ref{BCcat2}.
\end{prop}

\proof The argument is as in Proposition~\ref{BCliftSZ}, 
using again, as in Theorem~\ref{BCNori} 
the fact that maps of diagrams induce functors of the resulting categories of Nori motives.
\endproof

\medskip
\subsection{Nori geometric diagrams for assemblers, and a challenge}\label{NoriAssSec}   

We conclude this section on Bost--Connes systems and Nori motives by formulating 
a question about Nori diagrams and assembler categories.

\smallskip

According to the Nori formalism
as it is presented in \cite{HuM-S17}, we must start with
 a ``geometric'' category $C$ of spaces/varieties/schemes,
 possibly  endowed with additional
structures, in which one can define morphisms of closed 
embeddings $Y\hookrightarrow X$ (or $Y\subset X$) and morphisms of complements to closed embeddings
$X\smallsetminus Y \to X$,
Then the Nori diagram of
{\it effective pairs} $D(C)$ is defined as in \cite{HuM-S17}, 
pp.~207--208, see \S \ref{NoriGDiagSec}.

\smallskip

In the current context, {\it objects} of our category $C$ will be {\it assemblers}
$\cC$ (of course, described in terms of a category of lower level). In particular, each 
such $\cC$ is  endowed with a Grothendieck topology.

\medskip

A {\it vertex} of the Nori diagram $D(C)$ will be a triple $(\cC, \cC \setminus \cD, i)$
where its first two terms are taken from an abstract scissors congruence in $C$, and $i$
is an integer.
Intuitively, this means that we we are considering the canonical embedding 
$\cC \setminus \cD \hookrightarrow \cC$
as an analog of closed embedding. This intuition makes translation
of the remaining components of Nori's diagrams obvious, except for one: {\it what is the geometric meaning
of the integer $i$ in $(\cC, \cC \setminus \cD, i)$}\,? 

\medskip

The answer in the general context of assemblers, seemingly, was not yet
suggested, and already in the algebraic--geometric contexts is non--obvious
and non--trivial. Briefly, $i$ translates to the level of Nori geometric diagrams
the {\it weight filtration} of various cohomology theories
(cf.~\cite{HuM-S17},  10.2.2, pp.~238--241), and the existence of such translation
and its structure are encoded in  several versions of {\it Nori's Basic Lemma}
independently and earlier discovered by A.~Beilinson and K.~Vilonen
(cf.~\cite{HuM-S17},  2.5, pp.~45--59).

\medskip

The most transparent and least technical version of the Basic Lemma 
(\cite{HuM-S17},  Theorem 2.5.2 , p.~46) shows that in algebraic geometry
the existence of weight filtration is based upon special properties
of {\it affine schemes.} As we will see in the last section, lifts
of Bost--Connes algebras to the level of cohomology
 based upon the techniques  of {\it enhancement} also require
 a definition of {\it affine} assemblers. Since we do not know its combinatorial version,
 the enhancements that we can study now, force us to return to algebraic geometry.

\medskip

This challenge suggests to think about other possible 
geometric contexts in which dimensions/weights
of the relevant  objects may take, say, $p$--adic
values (as in the theory of $p$--adic weights of automorphic
forms  inaugurated by J.~P.~Serre), or rational values
(as it happens in some corners of ``geometries below ${\rm Spec}\,\bZ$''),
or even real values (as in various fractal geometries).

\medskip

{\it Can one transfer the scissors congruences imagery there?}

\medskip

See, for example, the formalism of Farey semi--intervals as base of
$\infty$--adic topology.

\bigskip

\subsection*{Acknowledgment} We thank the referee for many very detailed and useful 
comments and suggestions on how to improve the structure and presentation of the paper. 
The first and third authors were supported in part by the Perimeter
Institute for Theoretical Physics. The third author is also partially supported by
NSF grant DMS-1707882, and by NSERC Discovery Grant RGPIN-2018-04937 and Accelerator
Supplement grant RGPAS-2018-522593.

\end{document}